\documentclass[opre,nonblindrev]{informs3} % current default for manuscript submission

\usepackage{amssymb}
\usepackage{color}
\usepackage{rotating}
\usepackage{multirow}
\usepackage{latexsym}
\usepackage{secdot}
\usepackage{setspace}
\usepackage{wasysym}
\usepackage[us,12hr]{datetime}
\usepackage{bm}
\usepackage{pdfsync}
\usepackage{rccol}
\usepackage{mathrsfs}
\usepackage{epsfig}
\usepackage[margin=1in,paper=letterpaper]{geometry}
\usepackage{url}
\usepackage{multibib}
\newcites{appendix}{Online Appendix References}
\usepackage{enumitem}
\usepackage{stackrel}
\usepackage{soul}
\usepackage{multibib}

\newtheorem{thm}{Theorem}[section]

\newtheorem{lem}[thm]{Lemma}

{\theoremstyle{THkey}}

\OneAndAHalfSpacedXI % current default line spacing
\usepackage{natbib}
 \bibpunct[, ]{(}{)}{,}{a}{}{,}%
 \def\bibfont{\small}%
 \def\bibsep{\smallskipamount}%

\EquationsNumberedThrough    % Default: (1), (2), ...
\MANUSCRIPTNO{} 

\newcommand{\ts}{{\,}}
\newcommand{\Ccal}{{\mathcal C}}

\newcommand{\Fcal}{{\mathcal F}}
\newcommand{\Jcal}{{\mathcal J}}
\newcommand{\Kcal}{{\mathcal K}}
\newcommand{\Lcal}{{\mathcal L}}
\newcommand{\Pcal}{{\mathcal P}}
\newcommand{\Tcal}{{\mathcal T}}
\newcommand{\Wcal}{{\mathcal W}}
\newcommand{\Ccaltilde}{\widetilde{\mathcal C}}
\newcommand{\Fcaltilde}{\widetilde{\mathcal F}}
\newcommand{\Jtilde}{{\widetilde J}}
\newcommand{\Ztilde}{{\widetilde Z}}
\newcommand{\cbar}{{\overline c}}
\newcommand{\wbar}{{\overline w}}
\newcommand{\xbar}{{\overline x}}
\newcommand{\ybar}{{\overline y}}
\newcommand{\Dbar}{{\overline D}}
\newcommand{\Kbar}{{\overline K}}
\newcommand{\Zbar}{{\overline Z}}
\newcommand{\alphabar}{{\overline \alpha}}
\newcommand{\betabar}{{\overline \beta}}
\newcommand{\etabar}{{\overline \eta}}
\newcommand{\gammabar}{{\overline \gamma}}
\newcommand{\lambdabar}{{\overline \lambda}}

\newcommand{\Dhat}{{\widehat D}}
\newcommand{\Zhat}{{\widehat Z}}
\newcommand{\avec}{{\bm a}}
\newcommand{\cvec}{{\bm c}}
\newcommand{\uvec}{{\bm u}}
\newcommand{\xvec}{{\bm x}}
\newcommand{\yvec}{{\bm y}}
\newcommand{\wvec}{{\bm w}}
\newcommand{\Dvec}{{\bm D}}
\newcommand{\Jvec}{{\bm J}}
\newcommand{\alphavec}{{\bm \alpha}}
\newcommand{\betavec}{{\bm \beta}}
\newcommand{\etavec}{{\bm \eta}}
\newcommand{\muvec}{{\bm \mu}}
\newcommand{\wvecbar}{\overline{\bm w}}
\newcommand{\xvecbar}{\overline{\bm x}}
\newcommand{\yvecbar}{\overline{\bm y}}
\newcommand{\alphavecbar}{\overline{\bm \alpha}}
\newcommand{\betavecbar}{\overline{\bm \beta}}
\newcommand{\etavecbar}{\overline{\bm \eta}}
\newcommand{\lp}{{\text{\sf LP}}}
\newcommand{\apx}{{\text{\sf APX}}}
\newcommand{\opt}{{\text{\sf OPT}}}
\newcommand{\prf}{{\text{\small \sf PRF}}}
\newcommand{\exf}{{\text{\small \sf EXF}}}
\newcommand{\total}{{\text{\sf total}}}
\newcommand{\linear}{{\text{\sf linear}}}

\newcommand{\ru}[1]{\lceil #1 \rceil}
\newcommand{\ruBig}[1]{\Bigl \lceil #1 \Bigr\rceil}

\newcommand{\Ommega}{{O}}

\renewcommand{\qed}{\hfill \mbox{\raggedright \rule{0.1in}{0.1in}}}

\begin{document}
%%%%%%%%%%%%%%%%

\RUNAUTHOR{\normalfont Li, Rusmevichientong,  Topaloglu:}
\RUNTITLE{\normalfont{Revenue Management with Calendar-Aware and Dependent Demands}}

\TITLE{\vspace{-15mm}
\begin{minipage}{6.5in} \large \centering Technical Note -- Revenue Management with  Calendar-Aware and Dependent Demands: Asymptotically Tight Fluid Approximations \end{minipage}}

%\large \mbox{\!\!\!\!\!\!\!\! ~~~~~~~Revenue Management with Calendar\hspace{0.5mm}-Aware and Dependent Demands:~~~~~} \\ Asymptotically Tight Fluid Approximations \\ \normalsize \large(Technical Note)}

\ARTICLEAUTHORS{\vspace{-2mm}
\AUTHOR{Weiyuan Li$^1$, Paat Rusmevichientong$^2$, Huseyin Topaloglu$^1$}
\vspace{1mm}
\AFF{$^1$School of Operations Research and Information Engineering, Cornell Tech, New York, NY 10044
\\ 
$^2$Marshall School of Business, University of Southern California, Los Angeles, CA 90089
\\
\EMAIL{wl425@cornell.edu, rusmevic@marshall.usc.edu, topaloglu@orie.cornell.edu}}
\phantom{\small (Replies to the referee comments are appended to the end of the manuscript)}
%\today
\vspace{-3mm}
}
%% Enter all authors

\ABSTRACT{%
When modeling the demand in revenue management systems, a natural approach is to focus on a canonical interval of time, such as a week, so that we forecast the demand over each week in the selling horizon.~Ideally, we would like to use random variables with general distributions to model the demand over each week.~The current demand can give a signal for the future demand, so we also would like to capture the dependence between the demands over different weeks. Prevalent demand models in the literature, which are based on a discrete-time approximation to a Poisson process, are not compatible with these needs. In this paper, we focus on revenue management models that are compatible with a natural approach for forecasting the demand.~Building such models through dynamic programming is not difficult. We divide the selling horizon into multiple stages, each stage being a canonical interval of time on the calendar. We have random number of customer arrivals in each stage, whose distribution is arbitrary and depends on the number of arrivals in the previous stage.~The question we seek to answer is the form of the corresponding fluid approximation.~We give the correct fluid approximation in the sense that it yields asymptotically optimal policies. The form of our fluid approximation is surprising as its constraints use expected capacity consumption of a resource up to a certain time period, conditional on the demand in the stage just before the time period in question.~Letting~$K$~be the number of stages in the selling horizon, $c_{\min}$ be the smallest resource capacity and $\epsilon$ be a lower bound on the mass function~of the demand in a stage, we use the fluid approximation to give a policy with a performance guarantee of~$1 - \Ommega\Big( \frac{\sqrt{(c_{\min} +  K / \epsilon^6 ) \log c_{\min}}}{c_{\min}} \Big)$.~As the resource capacities and number of stages increase with the same rate, the performance guarantee converges to one. To our knowledge, this result gives the first asymptotically optimal policy under dependent demands with arbitrary distributions. When the demands in different stages are independent, letting $\sigma^2$ be the variance proxy for the demand in each stage, a similar performance guarantee holds by replacing $\frac{1}{\epsilon^6}$ with $\sigma^2$. Our computational work indicates that using the right fluid approximation can make a dramatic impact in practice. {\it Dated \today.}}

%\KEYWORDS{choice modeling, multinomial logit, social welfare optimization, totally unimodular constraints} 
%\HISTORY{This paper was first submitted on April 12, 1922 and has been with the authors for 83 years for 65 revisions.}

\maketitle

\vspace{-8mm}

\section{Introduction}

A natural approach for modeling the demand in revenue management systems is to focus on a canonical interval of time, such as a week, so that we forecast the demand over each week in the selling horizon.~Ideally, we would like to model the demand over each week through a random variable with an arbitrary distribution. Indeed, it is common for revenue management systems to produce forecasts stating that the demand, for example, during the week of July 17-23, 2023 has the log-normal distribution with mean 500 and standard deviation 250. These forecasts involve arbitrary demand distributions.~Also, these forecasts are aware of the calendar in the sense that they have a concept of when each week ends and the next one starts. Furthermore, the current demand usually gives a signal for the future demand, so we also would like to capture the dependence between the demands over different weeks. Prevalent demand models in the literature are not compatible with such a natural approach for forecasting the demand. In particular, a common demand model is based on dividing the selling horizon into a number of time periods such that there is at most one customer arrival at each time period and the arrivals at successive time periods are independent.~This model corresponds to a \mbox{discrete-time} approximation to a Poisson process, but it has shortcomings.~Under this model, the demand over an interval of time will always be approximately Poisson. To make matters worse, using $\lambda_t$ to denote the probability that we have a customer arrival at time period $t$, over an interval of $T$ time periods, the mean and standard deviation of the number of customer arrivals are, respectively,  $\sum_{t=1}^T \lambda_t$ and $\sqrt{\sum_{t=1}^T \lambda_t \ts (1-\lambda_t)}$. Noting that $\sqrt{\sum_{t=1}^T \lambda_t \ts (1-\lambda_t)} \leq \sqrt{\sum_{t=1}^T \lambda_t}$, the ratio between the standard deviation and mean of the demand is at most $1 / \sqrt{\sum_{t=1}^T \lambda_t}$. Thus, as the mean demand gets large, the coefficient of variation of the demand gets smaller. In other words, large demand variability and large demand volume cannot co-exist in this demand model. Finally, because this demand model is based on a Poisson process, the demands over different time intervals must be independent.

Motivated by the shortcomings discussed in the previous paragraph, we focus on revenue management models that are intrinsically compatible with a natural approach for forecasting the demand. Such a natural approach for forecasting the demand may specify the distribution of the demand over, for example, different weeks in the selling horizon, possibly along with the correlation structure between the demands in successive weeks. It is not too difficult to build these revenue management models through dynamic programming. We can divide the selling horizon into a number of stages, each stage representing a canonical interval of time on the calendar, such as a week. The number of customer arrivals in each stage is a random variable, whose distribution is arbitrary and depends on the number of customer arrivals in the previous stage. Therefore, along with the remaining capacities of the resources, the state variable in the dynamic program needs to keep track of the number of customer arrivals so far in the current stage and the number of customer arrivals in the previous stage, so that as a function of these two quantities, we can compute the probability of having one more customer arrival in the current stage. This dynamic program would give a precise specification of our model, but it is not useful for computing the optimal policy in practice because it involves a high-dimensional state variable. 

We focus on fluid approximations for our model so that we can construct practical policies with performance guarantees and compute upper bounds on the optimal total expected revenues.

{\bf \underline{Our Results and Contributions\phantom{p}\!\!\!}:} We start with a revenue management model based on a dynamic programming formulation that can handle arbitrary distributions for the number of customer arrivals in each stage and allow dependence between the number of customer arrivals in successive stages. A stage may correspond to a canonical interval of time over which forecasts are produced.~We capture the dependence between the demands through a Markov chain that specifies the distribution of the number of customer arrivals in one stage as a function of the number of customer arrivals in the previous stage. We give the correct fluid approximation for our revenue management model in the sense that the fluid approximation satisfies two properties. First, given that it is computationally difficult to find the optimal policy through a dynamic programming formulation, we can use our fluid approximation to construct approximate policies with performance guarantees. Second, we can use our fluid approximation to obtain an upper bound on the optimal total expected revenue, so that we can compare the total expected revenue of a heuristic policy with the upper bound to assess the optimality gap of the heuristic policy.

\vspace{-0.3mm}

{\it \underline{Structure of the Fluid Approximation}.} The structure of our fluid approximation turns out to be novel. In our fluid approximation, we use decision variables that allow the probability of accepting a customer request at a time period to depend on the number of customer arrivals in the previous stage. Due to the dependence between the demands in successive stages, this form of the decision variables is perhaps not surprising, but we are not aware of other fluid approximations with similar decision variables. More importantly, the constraints in our fluid approximation keep track of the expected capacity consumption of a resource up to a certain time period in a particular stage, conditional on the demand in the stage just before the time period in question. The conditional form of these constraints is surprising and does not appear in the literature. We show that the optimal objective value of our fluid approximation is an upper bound on the optimal total expected revenue. Thus, we can use this upper bound to assess the optimality gaps of heuristic policies.
 
\vspace{-0.3mm}
 
{\it \underline{Policies under Dependent Demands}.} Using our fluid approximation, we give an approximate policy. Letting $K$ be the number of stages in the selling horizon, $c_{\min}$ be the smallest resource capacity, $L$ be the maximum number of resources used by a product and $\epsilon$ be a lower bound on the mass function of the demand in a stage given the demand in the previous stage, we show that our approximate policy has a performance guarantee of $\max\Big\{ \frac{1}{4L} , \Big( 1 - 4 \ts \frac{\sqrt{(c_{\min} + (K-1) / \epsilon^6 ) \log c_{\min}}}{c_{\min}} - \frac{L}{c_{\min}} \Big) \Big\}$.  In many applications, the number of resources used by a particular product is usually small. In the airline setting, for example, the number of flight legs in an itinerary rarely exceeds two, corresponding to \mbox{$L=2$}. By the first term in the max operator, our approximate policy has a constant-factor performance guarantee when $L$ is uniformly bounded. By the second term in the max operator, as the number of stages and the capacities of the resources both increase linearly with rate~$\theta$, our performance guarantee converges to one with rate $1-\frac{1}{\sqrt \theta}$. 

\vspace{-0.3mm}

Thus, our approximate policy is asymptotically optimal for systems that command large demands for the products and involve large capacities for the resources. During the course of the proof of this result, we  also show that the upper bound provided by our fluid approximation is asymptotically tight in the same regime. To our knowledge, our approximate policy is the first to yield asymptotic optimality guarantee under dependent demands with arbitrary distributions. The proof for our performance guarantee uses techniques that have not been used in the related literature. In particular, under the standard demand model, to analyze policies from fluid approximations, we upper bound the total consumption of a resource  by using a random variable expressed as a sum of independent random variables. In this case, we can use a concentration inequality for sums of independent random variables to upper bound the tail probability of the total consumption of a resource, which in turn, yields a lower bound on the probability that the  policy has enough resource availabilities to accept the product request at each time period. The lower bound on the resource availability probabilities are used to lower bound the performance of the policy.

\vspace{-0.3mm}

Because the demands in successive stages are dependent in our setting, concentration inequalities for sums of independent random variables are not helpful to us, so we explicitly construct the concentration inequalities that we need. It is standard to use the moment generating function of a random variable to  bound its tail probabilities. In particular, if the moment generation function of the random variable $Z$ satisfies $\mathbb E\{ e^{\lambda Z} \} \leq f(\lambda)$ for all $\lambda \geq 0$, then we can bound its tail probabilities as $\mathbb P\{ Z \geq c\} = \mathbb P\{ e^{\lambda Z} \geq e^{\lambda c}\} \leq \frac 1{e^{\lambda c}} \mathbb E\{ e^{\lambda Z}\}\leq \frac{f(\lambda)}{e^{\lambda c}}$, where the first inequality is the Markov inequality.~We use martingales and the method of bounded differences to bound the moment generating functions of the capacity consumptions of the resources; see \cite{DuPa09}. This technique has been used in analyzing randomized algorithms, but their use in revenue management appears to be new. We hope that our use of this technique will further motivate other fluid approximations under even more sophisticated demand models.

\vspace{-0.3mm}

{\it \underline{Policies under Independent Demands}.} When the demands in different stages are independent, our revenue management model captures the case where the number of customer arrivals in, say, each week has an arbitrary distribution and the decision maker has a concept of when each week ends and the next one starts. Thus, this demand model is different from a demand model, where there is a single stage consisting of possibly multiple weeks and the number of customer arrivals in the single stage has an arbitrary distribution. In particular, the uncertainty in the demand in our model resolves sequentially over multiple stages. Recalling that the variance proxy of a sub-Gaussian random variable is an upper bound on its variance, under the assumption that the demand in each stage is \mbox{sub-Gaussian} with variance proxy $\sigma^2/200$, we give an approximate policy with a performance guarantee of \mbox{$\max\Big\{ \frac{1}{4L} , \Big( 1 - 4 \ts \frac{\sqrt{(c_{\min} + \sigma^2 (K-1)) \log c_{\min}}}{c_{\min}} - \frac{L}{c_{\min}} \Big) \Big\}$}. Variance may be a more intuitive statistic to work with than a lower bound on the mass function.

%Naturally, our earlier performance guarantee under dependent demands holds under independent demands, but variance may be a more intuitive statistic than a lower bound on the mass function.

\vspace{-0.3mm}

Even when the numbers of customer arrivals in different stages are independent, establishing the performance guarantee in the previous paragraph requires going one step beyond concentration inequalities for sums of independent random variables. In particular, the random variable that captures the demand in a particular stage creates dependence between the capacity consumptions of a resource at different time periods in the stage. Thus, as far as we can see, our performance guarantee with independent demands in each stage does not  follow from the proof techniques used for the existing asymptotic optimality results under a discrete-time approximation to a Poisson process.~We end up constructing the concentration inequalities that we need by exploiting the assumption of \mbox{sub-Gaussian} demand random variables. The assumption of sub-Gaussian demand random variables is relatively mild, as this class is rather general; see Section 2.1.2 in \cite{Wa19}.~Any bounded random variable, for example, is sub-Gaussian.

\vspace{-0.3mm}

{\it \underline{Computational Performance}.} To our knowledge, there is no work on asymptotically tight fluid approximations and asymptotically optimal policies for revenue management problems in which the demands over different time intervals are dependent and have arbitrary distributions. Building such fluid approximations and approximate policies is theoretically interesting, but fluid approximations with a sound theoretical footing can also make a significant impact in practice. In our computational experiments, we make comparisons with existing fluid approximations. While we can show that these fluid approximations do provide upper bounds, they do not provide asymptotically optimal policies. Our fluid approximation, owing to its sound theoretical footing, provides significantly tighter upper bounds and better approximate policies for a range of test problems.

%\vspace{-0.3mm}

{\bf \underline{Related Literature\phantom{p}\!\!\!}:} There is significant work on fluid approximations in revenue management, but this work is under demand models that use a discrete-time approximation to a Poisson process, ruling out the possibility of having arbitrary demand distributions and dependence between demands over different time intervals. Considering a revenue management problem with a single resource, \cite{GaRy94} show that if we scale the expected demand and the capacity of the resource with the same rate $\theta$, then a policy from a fluid approximation has a performance guarantee of~\mbox{$1-\Ommega(\frac{1}{\sqrt{\theta}})$}.~\cite{GaRy97} generalize this result to a network of resources, where the sale of different products consumes capacities of different combinations of resources. The policies in these papers use the primal solution to a fluid approximation.
\cite{TaRy98} use the dual solution to construct an asymptotically optimal policy in the same regime. \cite{LiVa2008} and \cite{gallego2004managing} construct similar asymptotically optimal policies under customer choice behavior, where the customers choose among the sets of products offered to them.~Considering~the~case where the customers with bookings do not necessarily show up at the time of service, \cite{KuTa11} give an asymptotically optimal policy that allows overbooking. 
The papers discussed so far solve the fluid approximation once at the beginning of the selling horizon. \cite{JaKu12} show that solving the fluid approximation periodically over the selling horizon provides policies with substantially better performance guarantees. Letting $c_{\min}$ be the smallest resource capacity, \cite{RuSuTo17} give a policy with a performance guarantee of \mbox{$1 -\Ommega\Big( \frac{1}{\sqrt[3]{c_{\min}}}\Big)$}. The asymptotic regime in their paper is different in the sense that the authors allow the expected demand to be scaled in an arbitrary fashion.~In the same asymptotic regime, \cite{MaSi21}, \cite{AoMa22}, \cite{BaHo22} and \cite{FeRa20} give policies with a performance guarantee of $1 - \Ommega\Big(\frac{1}{\sqrt{c_{\min}}}\Big)$. \cite{AoMa22} allow random number of customer arrivals. \cite{BaBe21} give a unified analysis for fluid approximations, allowing the possibility of solving the fluid approximation periodically over the selling horizon.

Under a random number of customer arrivals, \cite{BaHo23} give a policy with a performance guarantee of \mbox{1 - $\Ommega\Big(\frac{1}{\sqrt{c_{\min}}}\Big)$}. In their paper, the customer arrivals occur in one stage, so there is no possibility of introducing dependence between the numbers of customer arrivals in successive stages. Thus, the authors do not model the dependence between the numbers of customer arrivals over different time intervals, whereas our focus  is to precisely deal with such dependence. Also, even if the demands in successive stages are independent, when the uncertainty in the demand resolves sequentially over multiple stages, we demonstrate that intuitive modifications of the model in \cite{BaHo23} do not yield asymptotically optimal policies or upper bounds. Lastly, the performance guarantee in \cite{BaHo23} directly follows from one-sided Bernstein inequality, whereas we derive our concentration inequalities from scratch. Letting $L$ be the maximum number of resources used by a product, \cite{Ji23} gives a policy with a performance guarantee of $\frac{1}{1+L}$ when the distribution governing the products requested by the customers evolves from one time period to the next according to an exogenous Markov chain. The performance guarantee for the policy does not improve when we deal with large systems commanding large product demands and involving large resource capacities. Our performance guarantee converges to one when we deal with large systems. To establish the performance guarantee of $\frac{1}{1+L}$, \cite{Ji23} uses a linear program obtained by using linear value function approximations. We show that our fluid approximation is at least as tight as this linear program. Also, \cite{Ji23} does not use the linear program to construct a policy, whereas our policy directly uses the optimal solution to our fluid approximation.

{\bf \underline{Organization}:} In Section \ref{sec:form}, we formulate our revenue management model with  arbitrary demand distributions in each stage and dependence between the demands in successive stages. In Section \ref{sec:lp}, we give the fluid approximation corresponding to our model and show that its optimal objective value is an upper bound on the optimal total expected revenue. In Section \ref{sec:main_result}, we describe the approximate policy from the fluid approximation and give a performance guarantee for the approximate policy. In Section \ref{sec:perf_asymp}, we prove our performance guarantee. In Section \ref{sec:conc}, we conclude. In \cite{LiRuTo24}, we give extended results for our paper, where we focus our results to the case with independent demands in different stages and give computational experiments.

\newpage

\newpage

\section{Problem Formulation}
\label{sec:form}

The set of resources is $\Lcal$. The capacity of resource $i$ is $c_i$. The set of products is $\Jcal$. The revenue of product $j$ is $f_j$. The resources used by product $j$ are given by the vector \mbox{$\avec_j = (a_{ij} : i \in \Lcal) \in \mathbb \{0,1\}^{|\Lcal|}$}, where $a_{ij} = 1$ if and only if product $j$ uses resource $i$. We divide the selling horizon into $K$ stages indexed by $\Kcal = \{1,\ldots,K\}$. We use the random variable $D^k$ to capture the number of customer arrivals in stage $k$. There are at most $T$ customer arrivals in each stage. We divide each stage into $T$ time periods indexed by $\Tcal = \{1,\ldots,T\}$. We use $\lambda_{jt}^k$ to denote the probability that the customer arriving at time period $t$ in stage $k$ requests product $j$, so we have $\sum_{j \in \Jcal} \lambda_{jt}^k=1$. We refer $D^k$ as the demand in stage $k$. The demands in successive stages follow a Markov process. Thus, conditional on $D^k$, $D^{k+1}$ is independent of $D^1,\ldots,D^{k-1}$. We characterize the evolution of the demands by the survival rate function $\theta_t^k(q) = \mathbb P \{ D^k \geq {t+1} \ts | \ts  D^k \geq t,~D^{k-1} =q \}$, capturing the probability that the demand in stage $k$ is at least $t+1$, given that the demand in stage $k$ is at least $t$ and the demand in the previous stage was $q$. We assume that $\mathbb P \{ D^{k+1} = p \ts | \ts D^k = q\} \geq \epsilon$ for all $p,q \in \Tcal$ and $k \in \Kcal$ for some $\epsilon > 0$, so the demand in any stage takes values over its full support. %irrespective of the demand in the previous stage.

Each stage is a canonical interval of time on the calendar, such as, a day, a week or a month. We are aware of the calendar in the sense that we know when the current stage starts, but we only have probabilistic information about the number of customer arrivals in each stage. Customers in the current stage arrive one by one. Each arriving customer makes a request for a product. We decide whether to accept the product request. Our goal is to find a policy to decide which customer requests to accept so that we maximize the total expected revenue over the selling horizon. We give a dynamic program to compute the optimal policy. We use $\yvec = (y_i : i \in \Lcal) \in \mathbb Z_+^{|\Lcal|}$ to capture the state of the system, where $y_i$ is the remaining capacity of resource $i$. We use $\uvec = (u_j : j \in \Jcal) \in \{0,1\}^{|\Jcal|}$ to capture the decisions, where $u_j=1$ if and only if we accept a request for product $j$. The set of feasible decisions is given by $\Fcal(\yvec) = \{ \uvec \in \{0,1\}^{|\Jcal|} : a_{ij} \ts u_j \leq y_i ~~\forall \ts i \in \Lcal,~j \in \Jcal\}$, ensuring that if we want to accept a request for product $j$ and the product uses resource $i$, then we need to have at least one unit of remaining capacity for resource $i$. We can find the optimal policy by computing the value functions $\{J_t^k : t \in \Tcal,~k \in \Kcal\}$ through the dynamic program
\begin{align}
J_t^k(\yvec,q)
= \max_{\uvec \in \Fcal(\yvec)} \Bigg\{ \sum_{j \in \Jcal} \lambda_{jt}^k \ts \bigg\{ f_j \ts u_j + \theta_t^k(q) \ts J_{t+1}^k(\yvec -\avec_j \ts u_j ,q) + ( 1 - \theta_t^k(q)) \ts J_1^{k+1}(\yvec - \avec_j   \ts u_j,t) \bigg\} \Bigg\},
\label{eqn:dp}
\end{align}
with the boundary condition that $J_1^{K+1} = 0$. Note that the state variable above keeps both the remaining capacities of the resources and the demand in the previous stage.

In (\ref{eqn:dp}), we have a request for product $j$ at time period $t$ in stage $k$ with probability $\lambda_{jt}^k$. If we accept this request, then we generate a revenue of $f_j$ and consume the capacities of the resources used by the product. Given that the demand in the current stage is already $t$ and the demand in the previous stage was $q$, we have one more demand in the current stage with probability $\theta_t^k(q)$.~The demand in the stage right before the beginning of the selling horizon is  $D^0$ and it is fixed at $\Dhat^0$ as problem data. Using $\cvec = (c_i : i \in \Lcal)$ to denote the initial resource capacities, the optimal total expected revenue is $\opt = J_1^1(\cvec,\Dhat^0)$. %We allow dependent demands in different stages. 
The form of the survival rate function is general, so the demand in each stage can have an arbitrary distribution. In our model, we specify the distribution of the demand in the next stage conditional on the demand in the current stage, which implies a certain distribution for the total demand. In \cite{LiRuTo24}, we give extended results for our paper. In Extended Results \ref{sec:total_demand}, we show that we can calibrate our model to match a given distribution for the total demand. In our model, we decide whether to accept each product request, but our results hold when we make pricing or assortment offering decisions. We will consider such extensions.

%Our model is unique in the sense that it allows dependent demands in different stages. The form of the survival rate function is general, so the demand in each stage can have an arbitrary distribution. Therefore, our model can capture general demand distributions, while allowing dependent demands in successive stages. 

%We focus on developing fluid approximations corresponding to the dynamic program in (\ref{eqn:dp}) with the goal of obtaining tractable upper bounds and policies with performance guarantees. 

\section{Fluid Approximation}
\label{sec:lp}

The dynamic program in (\ref{eqn:dp}) involves a high-dimensional state variable, so it is computationally difficult to compute the optimal policy by solving this dynamic program. We give a fluid approximation that will serve two purposes. First, we will use the fluid approximation to obtain an upper bound on the optimal total expected revenue. In this case, we can compare the total expected revenue obtained by any policy with the upper bound to assess the optimality gap of the policy. Second, we will use the fluid approximation to construct a policy that is asymptotically optimal as the capacities of the resources and demands get large. In our fluid approximation, we use the decision variable $x_{jt}^k(q)$ to capture the probability of accepting a request for product $j$ at time period $t$ in stage $k$ given that the demand in the previous stage was~$q$.~Using the vector of decision variables $\xvec = (x_{jt}^k(q) : j \in \Jcal,~t,q \in \Tcal,~k \in \Kcal)$, to approximate the optimal total expected revenue over the selling horizon, consider linear program
\begin{align}
\Zbar_\lp ~=~ \max_{\xvec \in \mathbb R_+^{|\Jcal| T^2K}} ~~&  \sum_{k \in \Kcal} \sum_{t \in \Tcal} \sum_{q \in \Tcal} \sum_{j \in \Jcal} f_j \ts \mathbb P \{  D^k \geq t, \ts D^{k-1} = q \} \ts x_{jt}^k(q) 
\label{eqn:lp}
\\
\mbox{st}~~&\sum_{\ell =1}^{k-1} \sum_{s \in \Tcal} \sum_{p \in \Tcal} \sum_{j \in \Jcal} a_{ij} \ts \mathbb P \{ D^\ell \geq s, \ts D^{\ell-1} =p \ts|\ts D^k \geq t, \ts D^{k-1} = q \} \ts x_{js}^\ell(p) 
\nonumber
\\
& \qquad \qquad \qquad \qquad \qquad  + \sum_{s=1}^t \sum_{j \in \Jcal} a_{ij} \ts x_{js}^k(q) \leq c_i \qquad \forall \ts i \in \Lcal,~t,q \in \Tcal,~k \in \Kcal
\nonumber
\\
& x_{jt}^k(q) \leq \lambda_{jt}^k \qquad \forall \ts j \in \Jcal,~t,q \in \Tcal,~k \in \Kcal. \phantom{\Bigg\}}
\nonumber
\end{align}

In the linear program above, the objective function accounts for the total expected revenue over the selling horizon. In particular, we can make a sale for product $j$ at time period $t$ in stage~$k$ only if the demand in stage $k$ is at least $t$. Furthermore, if the demand in stage $k-1$ is $q$, then we make a sale for product $j$ at time period $t$ in stage $k$ with probability $x_{jt}^k(q)$. Therefore, the expression $\sum_{q \in \Tcal} \mathbb P \{  D^k \geq t, \ts D^{k-1} = q \} \ts x_{jt}^k(q)$ corresponds to the expected sales for product $j$ at time period $t$ in stage $k$. The left side of the first constraint corresponds to the total expected capacity consumption of resource $i$ up to and including time period $t$ in stage~$k$ \underline{\it conditional\phantom{p}\!\!\!} on the fact that the demand in stage $k$ is at least $t$ and the demand in stage $k-1$ is $q$. In the first sum, similar to the objective function, the expression $\sum_{p \in \Tcal} \mathbb P \{ D^\ell \geq s, \ts D^{\ell-1} =p \ts|\ts D^k \geq t, \ts D^{k-1} = q \} \ts x_{js}^\ell(p)$ corresponds to the expected sales for product $j$ at time period $s$ in stage $\ell$ conditional on the fact that the demand in stage $k$ is at least $t$ and the demand in stage $k-1$ is $q$. If product $j$ uses resource $i$, then these sales consume the capacity of resource $i$. In the second sum, conditional on the fact that the demand in stage $k$ is at least $t$ and the demand in stage $k-1$ is $q$, we can make a sale for product $j$ at all time periods in stage $k$ up to and including time period $t$. Furthermore, we accept a request for product $j$ at time~period $s$ in stage $k$ with probability $x_{js}^k(q)$. Thus, the first constraint is the capacity constraint. The second constraint is the demand constraint, ensuring that the probability of accepting a request for a product at any time period in any stage does not exceed the probability of getting the request.  We emphasize two novel aspects of the linear program given above.

First, because the demand in  stage~$k$ depends on the demand in stage $k-1$, the probability of accepting a request for a product at any time period in stage $k$ depends on the demand in stage $k-1$ as well. Second, perhaps more surprisingly, the first constraint keeps track of the total expected capacity consumption of a resource up to and including time period $t$ in stage $k$, conditional on the fact that the demand in stage $k$ is at least $t$ and the demand in stage $k-1$ is $q$. The form of this conditioning is unexpected. We can use the Markovian structure of the demands to slightly simplify the first constraint. If $\ell \leq k-1$, then conditional on $D^{k-1}$, $D^\ell$ and $D^{\ell-1}$ are independent of $D^k$. Thus, we can replace the probability $\mathbb P \{ D^\ell \geq s, \ts D^{\ell-1} =p \ts|\ts D^k \geq t, \ts D^{k-1} = q \}$ in the first sum with \mbox{$\mathbb P \{ D^\ell \geq s, \ts D^{\ell-1} =p \ts|\ts D^{k-1} = q \}$}. Furthermore, the second sum is increasing in $t$, so we can replace the sum $\sum_{s=1}^t \sum_{j \in \Jcal} a_{ij} \ts x_{js}^k(q)$ with $\sum_{s \in \Tcal} \sum_{j \in \Jcal} a_{ij} \ts x_{js}^k(q)$. Therefore, we can express the first constraint equivalently as $\mbox{$\sum_{\ell =1}^{k-1} \sum_{s \in \Tcal} \sum_{p \in \Tcal} \sum_{j \in \Jcal} a_{ij} \ts \mathbb P \{ D^\ell \geq s, \ts D^{\ell-1} =p \ts|\ts  D^{k-1} = q \} \ts x_{js}^\ell(p)$} + \sum_{s \in \Tcal} \sum_{j \in \Jcal} a_{ij} \ts x_{js}^k(q) \leq c_i$ for all $i \in \Lcal$, $q \in \Tcal$ and $k \in \Kcal$. In this way, we can reduce the number of constraints in the first constraint by a factor of $T$. Nevertheless, we believe that our fluid approximation, as stated in (\ref{eqn:lp}), is more instructive, so we keep it in its full form. Lastly, we can write \mbox{$\mathbb P \{ D^\ell \geq s, \ts D^{\ell-1} =p \ts|\ts  D^{k-1} = q \} \ts = \ts \mathbb P \{  D^{\ell-1} =p \ts|\ts  D^\ell \geq s, \ts D^{k-1} = q \} \ts \ts \ts \mathbb P\{ D^\ell \geq s \ts | \ts D^{k-1} = q \}$}, but it is not true that $\mathbb P \{  D^{\ell-1} =p \ts|\ts  D^\ell \geq s, \ts D^{k-1} = q \}  \ts = \ts \mathbb P \{  D^{\ell-1} =p \ts|\ts  D^\ell \geq s \} $, so there is no  more simplification.

Our demand model is a generalization of the standard demand model based on a discrete-time approximation to a Poisson process. If there is one stage in the selling horizon and the number of time periods in the stage takes a deterministic value, then our demand model reduces to the standard demand model based on a discrete-time approximation to a Poisson process. Accordingly, if there is one stage in the selling horizon and the number of time periods in the stage takes a deterministic value, then problem (\ref{eqn:lp}) reduces to the standard fluid approximation. In particular, if there is one stage and the number of time periods in the stage takes the deterministic value of $T$, then we have $\mathbb P \{ D^1 \geq t,~D^0 = q\} = 1$ for all $t = 1,\ldots,T$ and $q = \Dhat^0$, whereas \mbox{$\mathbb P \{ D^1 \geq t,~D^0 = q\} = 0$} for all $t = 1,\ldots,T$ and $q \neq \Dhat^0$.~Thus, noting that there is one stage, only the decision variables $(x_{jt}^1(\Dhat^0) : j \in \Jcal,~t \in \Tcal)$ have a non-zero objective function coefficient in (\ref{eqn:lp}), so we can drop all other decision variables. In this case, using $y_{jt}$ to succinctly denote the decision variable $x_{jt}^1(\Dhat^0)$, the objective function in (\ref{eqn:lp}) becomes $\sum_{t \in \Tcal} \sum_{j \in \Jcal} f_j \ts y_{jt}$. Because there is one stage, the first sum in the first constraint with $k =1$ is zero, so noting that we drop the decision variables other than \mbox{$(x_{jt}^1(\Dhat^0) : j \in \Jcal,~t \in \Tcal)$}, the first constraint in (\ref{eqn:lp}) becomes $\sum_{s=1}^t \sum_{j \in \Jcal} a_{ij} \ts y_{js} \leq c_i$ for all $i \in \Lcal$ and $t \in \Tcal$.  Lastly, the second constraint in (\ref{eqn:lp}) becomes $y_{jt} \leq \lambda_{jt}^1$ for all $j \in \Jcal$ and $t \in \Tcal$. Therefore, the linear program in (\ref{eqn:lp}) reduces to the standard fluid approximation.

{\bf \underline{Upper Bound on the Optimal Total Expected Revenue}:}
\\
\indent 
Turning back to our demand model with multiple stages in the selling horizon with the demand in each stage having an arbitrary distribution, possibly with dependence between the demands in different stages, we show that the optimal objective value of the linear program in (\ref{eqn:lp}) is an upper bound on the optimal total expected revenue. There are two uses of this result. From the practical side, it is difficult to compute the optimal policy, but we can compare the performance of any heuristic policy with the upper bound on the optimal total expected revenue to assess the optimality gap of the heuristic policy. If the total expected revenue of the heuristic policy is close to the upper bound, then we can conclude that the heuristic policy is satisfactory. From the theoretical side, we will give performance guarantees for a policy that is obtained by using the linear program in (\ref{eqn:lp}). Because it is difficult to compute the optimal policy, we will establish these performance guarantees by comparing the total expected revenue of the policy with the upper bound on the optimal total expected revenue. Thus, the linear program in (\ref{eqn:lp}) also becomes useful to establish performance guarantees. In the next theorem, we show that the optimal objective value of problem (\ref{eqn:lp}) is indeed an upper bound on the optimal total expected revenue. 

\begin{thm}[Upper Bound] 
\label{thm:bound}
Using $\opt$ to denote the optimal total expected revenue and $\Zbar_\lp$ to denote the optimal objective value of problem $(\ref{eqn:lp})$, we have $\Zbar_\lp \geq \opt$. 
\end{thm}

The proof of the theorem, which is based on relaxing the resource availability constraints in (\ref{eqn:dp}) through Lagrange multipliers, is in Extended Results \ref{sec:bound}. There are other possible approaches to show that problem (\ref{eqn:lp}) provides an upper bound on the optimal total expected revenue. 
In Extended Results \ref{sec:lin_app}, we use linear approximations to the value functions. To compute the slopes and intercepts of the  approximations, we plug the approximations into the linear programming formulation of the dynamic program in (\ref{eqn:dp}). A suitable relaxation of this linear program is the dual of problem (\ref{eqn:lp}). \cite{Ji23} uses a similar outline to develop a linear program that provides an upper bound. In the development in our extended results, we also show that our relaxation is at least as tight as the one in \cite{Ji23}. In Extended Results \ref{sec:path_bound}, we use the decisions made by the optimal policy to construct a feasible solution to problem (\ref{eqn:lp}). This approach is not constructive in the sense that it does not allow us to derive the form of problem (\ref{eqn:lp}), but it simply verifies that the optimal objective value of problem (\ref{eqn:lp}) is an upper bound.
Lastly, another approach for obtaining an upper bound on the optimal total expected revenue is to use an offline bound, where we observe the realizations of all product requests and solve an optimization problem to choose the product requests to accept. Under the standard demand model with no demand dependence, the offline bound is at least as tight as the one provided by the fluid approximation; see \cite{TaRy99}. The analogue of this result does not hold under our demand model with demand dependence.~In Extended Results \ref{sec:offline_bound}, we give problem instances to demonstrate that the offline bound can be tighter or looser than the optimal objective value of problem (\ref{eqn:lp}).

\section{Approximate Policy and Asymptotic Optimality}
\label{sec:main_result}

We use an optimal solution to problem (\ref{eqn:lp}) to construct an approximate policy. We show that this policy has a constant-factor performance guarantee, but if both the number of stages in the selling horizon and capacities of the resources increase with the same rate, then the policy is asymptotically optimal. Thus, we expect the approximate policy to perform particularly well under large product demands and large resource capacities, but because of its constant-factor guarantee, the approximate policy can never perform arbitrarily badly. In our approximate policy, we solve the linear program in (\ref{eqn:lp}) once at the beginning of the selling horizon.~Letting \mbox{$\xvecbar = (\xbar_{jt}^k(q) : j \in \Jcal,~t,q \in \Tcal,~k \in \Kcal)$} be an optimal solution, we make the decisions as follows.

{\bf \underline{Approximate Policy from the Fluid Approximation}:}
\\
\indent Using $\gamma \in [0,1]$ to denote a tuning parameter, if we have a request for product $j$ at time period~$t$ in stage $k$ and the demand in stage $k-1$ was $q$, then we are willing to accept the request with probability $\gamma \ts \frac{\xbar_{jt}^k(q)}{\lambda_{jt}^k}$. If we are willing to accept the request and there are enough resource capacities to accept the request, then we accept the request. Otherwise, we reject. 

The description of the approximate policy as given above will be adequate to establish its performance guarantee. In Extended Results \ref{sec:policy_steps}, we give a detailed description of our approximate policy. Letting $c_{\min} = \min_{i \in \Lcal} c_i$ be the smallest resource capacity and $L = \max_{j \in \Jcal} \sum_{i \in \Lcal} a_{ij}$ be the maximum number of resources used by a product, we have the next performance guarantee. 

\begin{thm}[Performance Guarantee] 
\label{thm:perf}
Using $\apx$ to denote the total expected revenue of the approximate policy, there exists a choice of the tuning parameter $\gamma$ such that we have 
\begin{align*}
\frac{\apx}{\opt} 
~\geq~ 
\frac{\apx}{\Zbar_\lp} 
~\geq~ 
\max\Bigg\{ \frac{1}{4L} , \Bigg( 1 - 4 \ts \frac{\sqrt{(c_{\min} + \frac{1}{\epsilon^6} \ts (K-1) ) \log c_{\min}}}{c_{\min}} - \frac{L}{c_{\min}} \Bigg) \Bigg\}.
\end{align*}
\end{thm}

We devote the next section to the proof of the theorem. To our knowledge, this theorem gives the first policy with an asymptotic performance guarantee under dependent demands. The proof involves novel ideas. We use the moment generating function of the resource capacity consumptions to lower bound the probability that we have enough resource capacities to accept a product request. Because of the dependence between the demands in different stages, it is difficult to characterize the moment generation function of the resource capacity consumptions.~We bound the moment generating functions by using martingales and the method of bounded differences; see Chapter 5 in \cite{DuPa09}. Dependence between the demands requires us to derive our own moment generating function bounds, which ultimately yield the tail probability bounds needed for our performance guarantee. If there is one stage in the selling horizon and the number of time periods in the stage takes a deterministic value, then our demand model reduces to the standard demand model based on a discrete-time approximation to a Poisson process. Other papers, such as \cite{MaRuSuTo18} and \cite{BaMa22}, focus on performance guarantees under the standard demand model. Our demand model is a generalization of the standard demand model, so Theorem \ref{thm:perf} holds under the standard demand model. Setting $K=1$ in Theorem \ref{thm:perf}, our approximate policy has a performance guarantee of $\max \Big\{ \frac{1}{4L}, 1 - 4 \sqrt{\frac{\log c_{\min}}{c_{\min}}} - \frac{L}{c_{\min}} \Big\}$ under the standard demand model. \cite{FeRa20} also study fluid approximations under the standard demand model.

We proceed to interpreting the two parts in the performance guarantee in Theorem \ref{thm:perf}. In many network revenue management settings, the number of resources and number of products can be large, but the number of resources used by a particular product remains bounded. In the airline setting, for example, we may have hundreds of flight legs and thousands of itineraries, but the number of flight legs in an itinerary rarely exceeds two, corresponding to $L=2$. Thus, the first part in the performance guarantee provides a constant-factor performance guarantee for the approximate policy when $L$ is uniformly bounded. On the other hand, consider a regime where we scale both the number of stages and resource capacities with the same rate $\theta$, so that $K = \theta \ts \Kbar$ and $c_{\min} = \theta \ts \cbar$ for some fixed $\Kbar,\cbar \in \mathbb Z_+$. If $\theta$ gets large, then the expected demands for the products and  capacities for the resources both get large. Letting $\apx^\theta$ be the total expected revenue from the approximate policy, $\opt^\theta$ be the optimal total expected revenue and $\Zbar_\lp^\theta$ be the optimal objective value of problem (\ref{eqn:lp}) when we scale the number of stages and resource capacities with $\theta$, by Theorem \ref{thm:perf}, we have $1 \geq \frac{\apx^\theta}{\opt^\theta} \geq \frac{\apx^\theta}{\Zbar_\lp^\theta} \geq 1 - \frac{4}{\epsilon^3} \frac{\sqrt{(\cbar + \Kbar) \ts \log(\theta \ts \cbar)}}{ \cbar \sqrt{\theta} } - \frac{L}{\theta \ts \cbar}$. In this case, ignoring the logarithmic terms, as $\theta$ gets large, the relative gap between the total expected revenue of the approximate policy and the optimal total expected revenue converges to one with rate $1 - \frac{1}{\sqrt \theta}$. Therefore, as the number of stages in the selling horizon and capacities of the resources increase with the same rate, the approximate policy is asymptotically optimal. Similarly, as $\theta$ gets large, the relative gap between the total expected revenue of the approximate policy and the optimal objective value of the fluid approximation converges to one with rate $1 - \frac{1}{\sqrt \theta}$ as well. In Extended Results \ref{sec:asymp_ub}, we also give a problem instance such that if we scale the number of stages and resource capacities with the same rate $\theta$, then the optimal total expected revenue and the optimal objective value of problem (\ref{eqn:lp}) satisfy $\frac{\opt^\theta}{\Zbar_\lp^\theta} \leq 1 - \frac{1}{\sqrt{26 \ts \theta}}$.  In other words, the optimal total expected revenue is at most within a factor of $O\Big(1 - \frac{1}{\sqrt \theta} \Big)$ of the optimal objective value of the fluid approximation. Therefore, we cannot improve the asymptotic performance guarantee for the approximate policy as long as we establish this performance guarantee by comparing the total expected revenue of the approximate policy with the optimal objective value of the fluid approximation. In this sense, the performance guarantee that we give for the approximate policy is asymptotically tight. 

\vspace{-1mm}

The scaling regime in the previous paragraph increases the number of stages and resource capacities with the same rate. Because the demands in different stages are dependent, increasing the number of stages in the selling horizon is perhaps the most natural approach to increase the expected demands for the products. In this way, we can increase the expected demands for the products without distorting the correlation structure for the demands in different stages. Another approach to increase the expected demands for the products could be to increase the support of the demand in each stage, while keeping the number of stages in the selling horizon fixed. Scaling the expected demand in this fashion can potentially distort the correlation structure for the demands in different stages. In Extended Results \ref{sec:large_capacity}, we also give a counterexample to demonstrate that the relative gap between the total expected revenue of the approximate policy and the optimal objective value of the fluid approximation in (\ref{eqn:lp}) does not necessarily converge to one as we increase the support of the demand in each stage and capacities of the resources with the same rate, while keeping the number of stages constant. In our counterexample, we give a problem instance with  three stages. There is a single resource with a capacity of $C+1$. The largest value of the demand in a week is $C$. There are two products. We show that the optimal total expected revenue is $C$, whereas the optimal objective value of problem (\ref{eqn:lp}) is $\frac{5}{4} \ts C$. Thus, we have $\frac{\apx^\theta}{\Zbar_\lp^\theta} \leq \frac{\opt^\theta}{\Zbar_\lp^\theta} = \frac 45$. For this problem instance, no matter how large $C$ is, the ratio $\frac{\apx^\theta}{\Zbar_\lp^\theta}$ always stays away from one. 

\section{Performance Guarantee}
\label{sec:perf_asymp}

In this section, we give a proof for Theorem \ref{thm:perf}. We focus on showing the performance guarantee $\frac{\apx}{\opt} \geq \frac{\apx}{\Zbar_\lp}  \geq 1 - 4 \ts \frac{\sqrt{(c_{\min} + \frac{1}{\epsilon^6} \ts (K-1) ) \log c_{\min}}}{c_{\min}} - \frac{L}{c_{\min}}$. To establish this performance guarantee, we use ideas that have not been used in the revenue management literature to analyze fluid approximations. In Extended Results \ref{sec:constant_factor}, we turn our attention to showing the performance guarantee \mbox{$\frac{\apx}{\opt} \geq \frac{\apx}{\Zbar_\lp}  \geq \frac{1}{4L}$}, which is more straightforward. A common approach for analyzing approximate policies from fluid approximations involves using an auxiliary random variable to upper bound the capacity consumption of a resource under the approximate policy. Thus, we can use a concentration inequality to upper bound the tail probabilities of the auxiliary random variable, in which case, we can lower bound the probability that there is enough capacity to accept different product requests at different time periods in the selling horizon; see, for example, \cite{FeRa20}. This approach usually exploits the fact that the auxiliary random variable can be expressed as a sum of independent random variables, which facilitates using concentration inequalities for sums~of independent random variables. Because the demands in different stages are dependent in our setting, we cannot construct similar auxiliary random variables that can be expressed as sums of independent random variables.~Thus, we resort to new ideas.

{\bf \underline{Preliminary Random Variables and Availability Probabilities}:}
\\
\indent We define four classes of Bernoulli random variables for each $k \in \Kcal$ and $t \in \Tcal$. Analogues of these random variables have been used in the analysis of other fluid approximations.

\indent {$\bullet$ \it {Demand in Each Stage}.} For each $q \in \Tcal$, the random variable $\Psi_t^k(q)$ takes value one if we reach time period $t$ in stage $k$ before this stage is over and the demand in stage $k-1$ is $q$. In other words, letting ${\bf 1}(\cdot)$ be the indicator function, $\Psi_t^k(q) = {\bf 1}(D^k \geq t,~D^{k-1} = q)$. 

\indent {$\bullet$ \it {Product Request}.} For each $j \in \Jcal$, the random variable $A_{jt}^k$ takes value one if the customer arriving at time period $t$ in stage $k$ requests product $j$. We have $\mathbb P \{ A_{jt}^k  =1 \} = \lambda_{jt}^k$. The random variables $\{ A_{jt}^k : t \in \Tcal,~k \in \Kcal\}$ are independent of each other. 

\indent {$\bullet$ \it {Policy Decision}.} For each $j \in \Jcal$ and $q \in \Tcal$, the random variable $X_{jt}^k(q)$ takes value one if the approximate policy is willing to accept a request for product $j$ at time period $t$ in stage $k$ when the demand in stage $k-1$ was $q$. By our approximate policy, $\mathbb P \{ X_{jt}^k(q) = 1\} = \gamma \ts \frac{\xbar_{jt}^k(q)}{\lambda_{jt}^k}$.

\indent {$\bullet$ \it {Availability}.} For each $j \in \Jcal$, the random variable $G_{jt}^k$ takes value one if we have enough capacity to accept a request for product $j$ at time period $t$ in stage $k$ under the approximate policy. Instead of calculating the probability $\mathbb P \{ G_{jt}^k =1\}$, we will lower bound $\mathbb P \{ G_{jt}^k =1 \ts | \ts \Psi_t^k(q) = 1\}$.

The random variables $A_{jt}^k$ and $X_{jt}^k(q)$ are both simple Bernoulli draws independent of the decisions of the approximate policy, remaining capacities of the resources or realizations of the demands in different stages. Under the approximate policy, the sales for product $j$ at time period~$t$ in stage $k$ is given by $\sum_{q \in \Tcal} \Psi_t^k(q) \ts G_{jt}^k \ts A_{jt}^k \ts X_{jt}^k(q)$, where we use the fact that we sell product $j$ at time period~$t$~in stage $k$ if we reach time period $t$ in stage $k$, there is enough capacity to accept a request for product $j$, we have a request for the product and the approximate policy is willing to accept the product request. The remaining capacities of the resources at time period $t$ in stage $k$ depend on the requests for the products and decisions of the approximate policy at the earlier time periods, but not at time period $t$ in stage $k$. Thus, taking expectations, the expected sales for product $j$ at time period $t$ in stage~$k$ is $\sum_{q \in \Tcal} \mathbb P\{ \Psi_t^k(q) = 1\} \ts \mathbb P \{ G_{jt}^k = 1 \ts | \ts \Psi_t^k(q) = 1\} \ts \mathbb P\{ A_{jt}^k = 1\} \ts \mathbb P \{ X_{jt}^k(q) = 1 \}$. In this case, noting that $\mathbb P \{ A_{jt}^k = 1\} = \lambda_{jt}^k$ and $\mathbb P \{ X_{jt}^k(q) = 1\} = \gamma \ts \frac{\xbar_{jt}^k(q)}{\lambda_{jt}^k}$, we can write the last expectation equivalently as $\sum_{q \in \Tcal} \mathbb P \{ D^k \geq t,~D^{k-1} = q \} \ts \mathbb P \{ G_{jt}^k = 1 \ts | \ts \Psi_t^k(q) = 1\} \ts \gamma \ts \xbar_{jt}^k(q)$, capturing the expected sales for product $j$ at time period $t$ in stage $k$ under the approximate policy. Thus, the total expected revenue of the approximate policy is given by
\begin{align}
\apx
~=~ 
\sum_{k \in \Kcal} \sum_{t \in \Tcal} \sum_{j \in \Jcal} \sum_{q \in \Tcal} f_j \ts \mathbb P \{ D^k \geq t,~D^{k-1} = q \} \ts\ts \mathbb P \{ G_{jt}^k = 1 \ts | \ts \Psi_t^k(q) = 1\} \ts\ts \gamma \ts \xbar_{jt}^k(q).
\label{eqn:apx_rev}
\end{align}
By the definition of $\xvecbar$, we have \mbox{$\Zbar_\lp = \sum_{k \in \Kcal} \sum_{t \in \Tcal} \sum_{j \in \Jcal} \sum_{q \in \Tcal} f_j \ts \mathbb P \{  D^k \geq t, \ts D^{k-1} = q \} \ts \xbar_{jt}^k(q)$}, so if we can show that $\mathbb P \{ G_{jt}^k = 1 \ts | \ts \Psi_t^k(q) = 1\}  \geq \alpha$, then we get $\apx \geq \gamma \ts \alpha \ts \Zbar_\lp$.

Motivated by the discussion in the previous paragraph, we focus on lower bounding the availability probability $\mathbb P \{ G_{jt}^k = 1 \ts | \ts \Psi_t^k(q) = 1\}$. Under the approximate policy, the sales for product~$j$ at time period $t$ in stage $k$ is $\sum_{q \in \Tcal} \Psi_t^k(q) \ts G_{jt}^k \ts A_{jt}^k \ts X_{jt}^k(q)$, so $\sum_{q \in \Tcal} \Psi_t^k(q) \ts A_{jt}^k \ts X_{jt}^k(q)$ is an upper bound on these sales for product $j$ at time period~$t$ in stage~$k$.~In this case, $\sum_{j \in \Jcal} \sum_{q \in \Tcal} a_{ij} \ts \Psi_t^k(q) \ts A_{jt}^k \ts X_{jt}^k(q)$ is an upper bound on the capacity consumption of resource $i$ at time period $t$ in stage $k$. Letting $N_{it}^k(q) = \sum_{j \in \Jcal} a_{ij} \ts A_{jt}^k \ts  X_{jt}^k(q)$, we express our upper bound on the capacity consumption of resource~$i$ at time period $t$ in stage $k$ as $\sum_{q \in \Tcal} \Psi_t^k(q) \ts N_{it}^k(q)$. Note that  $\{ N_{it}^k(q) : t \in \Tcal,~k \in \Kcal\}$ are Bernoulli random variables and they are independent of each other. Having $\sum_{\ell=1}^{k-1} \sum_{s \in \Tcal} \sum_{q \in \Tcal} \Psi_s^\ell(q) \ts N_{is}^\ell(q) + \sum_{s=1}^t \sum_{q \in \Tcal} \Psi_s^k(q) \ts N_{is}^k(q) < c_i$ implies that the total capacity consumption of resource $i$ up to and including time period $t$ in stage $k$ does not exceed the capacity of the resource, in which case, we have capacity available for resource $i$ at time period $t$ in stage $k$. Therefore, letting \mbox{$\Lcal_j = \{ i \in \Lcal: a_{ij} = 1\}$} to capture the set of resources used by product $j$, we obtain 
\begin{align*}
\mathbb P \{ G_{jt}^k = 1 \ts | \ts \Psi_t^k(q) = 1\} \geq 
\mathbb P \bigg\{\! \sum_{\ell=1}^{k-1} \sum_{s \in \Tcal} \sum_{p \in \Tcal} \Psi_s^\ell(p) \ts N_{is}^\ell(p) + \sum_{s=1}^t \sum_{p \in \Tcal} \Psi_s^k(p) \ts N_{is}^k(p) < c_i ~\forall \ts i \in \Lcal_j \ts \Big| \ts \Psi_t^k(q) =1 \! \bigg\},
\end{align*}
where we use the fact if the upper bounds on the consumption of the resources used by product~$j$ do not exceed their capacities, then we have capacities to accept a request for product $j$.

By the definition of $\Psi_t^k(q)$, having $\Psi_t^k(q) = 1$ is equivalent to having $D^k \geq t$ and $D^{k-1} = q$. Thus, if $\Psi_t^k(q) = 1$, then we have $D^k \geq s$ for all $s=1,\ldots,t$, $D^{k-1} = q$ and $D^{k-1} \neq p$ for all $p \in \Tcal \setminus \{q\}$. Therefore, if $\Psi_t^k(q) = 1$, then $\Psi_s^k(q) = 1$ for all $s=1,\ldots,t$ and $\Psi_s^k(p) =0$ for all $p \in \Tcal \setminus \{q\}$. In this case, the inequality above is equivalent to
\begin{align}
& \mathbb P \{ G_{jt}^k = 1 \ts | \ts \Psi_t^k(q) = 1\} 
~\geq~ 
\mathbb P \bigg\{ \sum_{\ell=1}^{k-1} \sum_{s \in \Tcal} \sum_{p \in \Tcal} \Psi_s^\ell(p) \ts N_{is}^\ell(p) + \sum_{s=1}^t  N_{is}^k(q) < c_i ~\forall \ts i \in \Lcal_j \ts \Big| \ts \Psi_t^k(q) =1  \bigg\}
\nonumber
\\
& \qquad \qquad \qquad =~
1 - \mathbb P \bigg\{\sum_{\ell=1}^{k-1} \sum_{s \in \Tcal} \sum_{p \in \Tcal} \Psi_s^\ell(p) \ts N_{is}^\ell(p) + \sum_{s=1}^t  N_{is}^k(q) \geq c_i ~\mbox{for some } i \in \Lcal_j \ts \Big| \ts \Psi_t^k(q) =1  \bigg\}
\nonumber
\\
& \qquad \qquad \qquad \geq~
1 - \sum_{i \in \Lcal_j} \mathbb P \bigg\{\sum_{\ell=1}^{k-1} \sum_{s \in \Tcal} \sum_{p \in \Tcal} \Psi_s^\ell(p) \ts N_{is}^\ell(p) + \sum_{s=1}^t N_{is}^k(q) \geq c_i \ts \Big| \ts \Psi_t^k(q) =1  \bigg\}
\nonumber
\\
& \qquad \qquad \qquad =~
1 - \sum_{i \in \Lcal_j} \mathbb P \bigg\{\sum_{\ell=1}^{k-1} \sum_{s \in \Tcal} \sum_{p \in \Tcal} \Psi_s^\ell(p) \ts N_{is}^\ell(p) + \sum_{s=1}^t N_{is}^k(q) \geq c_i \ts \Big| \ts D^{k-1} = q  \bigg\},
\label{eqn:union}
\end{align}
where the second inequality is the union bound and the second equality holds because given $D^{k-1}$, $D^1,\ldots,D^{k-1}$ are independent of $D^k$. Thus, it is enough to upper bound the last probability in (\ref{eqn:union}).

{\bf \underline{Moment Generating Function Bounds}:}
\\
\indent The discussion so far in this section has been following standard arguments, but we proceed to introducing new ideas. To upper bound the last probability on the right side of (\ref{eqn:union}), letting \mbox{$n_{it}^k(q) = \mathbb E \{ N_{it}^k(q)\}$}, for all $i \in \Lcal$ and $k \in \Kcal$, we define $U_i^k = \sum_{\ell=1}^k \sum_{s \in \Tcal} \sum_{p \in \Tcal} \Psi_s^\ell(p) \ts N_{is}^\ell(p)$ and $V_i^k = \sum_{\ell=1}^k \sum_{s \in \Tcal} \sum_{p \in \Tcal} \Psi_s^\ell(p) \ts n_{is}^\ell(p)$. Using the vector of random variables $\Dvec^{[\ell,k]} = (D^\ell,\ldots,D^k)$ for notational brevity, note that the random variable $V_i^k$ is a deterministic function of $\Dvec^{[1,k]}$. Because $e^x$ is convex in~$x$, using the Jensen inequality, it is simple to show that $\mathbb E\{ e^{\lambda  V_i^k} \ts | \ts D^k \} \leq \mathbb E\{ e^{\lambda  U_i^k} \ts | \ts D^k\}$ for all $\lambda \geq 0$, so the moment generating function of the random variable $U_i^k$ conditional on $D^k$ upper bounds its counterpart for the random variable $V_i^k$. In the next lemma, we characterize the gap between the two moment generating functions. 

\begin{lem}[Moment Generating Function Gap]
\label{lem:mgf_gap}
For all $k \in \Kcal$, $i \in \Lcal$ and \mbox{$\lambda \geq 0$}, we have $\mathbb E\{ e^{\lambda U_i^k} \ts | \ts D^k \} \leq e^{\frac{1}{2\epsilon} k  \lambda^2} \ts \mathbb E\{ e^{\lambda  V_i^k} \ts | \ts D^k\}$.
\end{lem}

\noindent{\it Proof:} The random variables $\{ N_{is}^\ell(p) : s \in \Tcal,~\ell =1,\ldots,k\}$ are independent of each other. 
\
Also, for $\ell=1,\dots,k$ and $s,p \in \Tcal$, $\Psi_s^\ell(p)$ is a deterministic function of $\Dvec^{[1,k]}$. Thus, we have
\begin{multline}
\mathbb E\{ e^{\lambda(U_i^k - V_i^k)} \ts | \ts \Dvec^{[1,k]} \}
~=~
\mathbb E\{ e^{\lambda \sum_{\ell=1}^k \sum_{s \in \Tcal} \sum_{p \in \Tcal} \Psi_s^\ell(p) (N_{is}^\ell(p) - n_{is}^\ell(p))} \ts | \ts \Dvec^{[1,k]} \}
\\
~=~
\prod_{\ell =1}^k \prod_{s \in \Tcal}
\mathbb E\{ e^{\lambda \sum_{p \in \Tcal} \Psi_s^\ell(p) (N_{is}^\ell(p) - n_{is}^\ell(p))} \ts | \ts \Dvec^{[1,k]} \}.
\label{eqn:diff_caps}
\end{multline}
Because $\sum_{p \in \Tcal} \Psi_s^\ell(p) \leq 1$ and  the random variable $N_{is}^\ell(p)$ is Bernoulli with expectation $n_{is}^\ell(p)$, we have $\sum_{p \in \Tcal} \Psi_s^\ell(p) (N_{is}^\ell(p) - n_{is}^\ell(p)) \in [-1,1]$. If the mean-zero random variable $Z$ is bounded by $[a,b]$, then we have $\mathbb E\{ e^{\lambda Z}\} \leq e^{\frac 18 (b-a)^2 \lambda^2}$ for any $\lambda \geq 0$; see Lemma 5.1 in \cite{DuPa09}. Thus, the last conditional expectation on the right side of (\ref{eqn:diff_caps}) is upper bounded by $e^{\frac12 \lambda^2}$, so by  (\ref{eqn:diff_caps}), we obtain \mbox{$\mathbb E\{ e^{\lambda(U_i^k - V_i^k)} \ts | \ts \Dvec^{[1,k]} \} \leq e^{\frac 12 k T \lambda^2}$}. The random variable $V_i^k$ is a deterministic function of $\Dvec^{[1,k]}$. In this case, using the tower property of conditional expectations, we get 
$\mathbb E\{ e^{\lambda U_i^k} \ts | \ts D^k \} = \mathbb E \{ \ts \mathbb E\{ e^{\lambda \ts V_i^k} \ts e^{\lambda \ts (U_i^k - V_i^k)} \ts | \ts \Dvec^{[1,k]} \} \ts | \ts D^k \} = \mathbb E \{ e^{\lambda V_i^k} \ts \mathbb E\{ e^{\lambda \ts (U_i^k - V_i^k)} \ts | \ts \Dvec^{[1,k]}  \} \ts | \ts D^k \}$. Using the fact that $\mathbb E\{ e^{\lambda(U_i^k - V_i^k)} \ts | \ts \Dvec^{[1,k]} \} \leq e^{\frac 12 k T \lambda^2}$,  we obtain \mbox{$\mathbb E\{ e^{\lambda U_i^k} \ts | \ts D^k \} \leq e^{\frac 12 k T \lambda^2} \mathbb E\{ e^{\lambda V_i^k} \ts | \ts D^k\}$} by the last chain of equalities. The result follows by noting the assumption that $\mathbb P \{ D^{k+1} = p \ts | \ts D^k = q\} \geq \epsilon$ for all $p,q \in \Tcal$, so  $1 = \sum_{p \in \Tcal} \mathbb P \{ D^{k+1} = p \ts | \ts D^k = q\} \geq T \epsilon$, which implies that $T \leq 1/\epsilon$. \qed

If the random variable $Z$ satisfies $\mathbb E\{ e^{\lambda Z} \} \leq f(\lambda)$ for all $\lambda \geq 0$, then we can upper bound its tail probabilities as $\mathbb P\{ Z \geq c\} = \mathbb P\{ e^{\lambda Z} \geq e^{\lambda c}\} \leq \frac 1{e^{\lambda c}} \mathbb E\{ e^{\lambda Z}\}\leq \frac{f(\lambda)}{e^{\lambda c}}$, where the first inequality uses the Markov inequality.~In the last probability in (\ref{eqn:union}), we have $U_i^{k-1} = \sum_{\ell=1}^{k-1} \sum_{s \in \Tcal} \sum_{p \in \Tcal} \Psi_s^\ell(p) \ts N_{is}^\ell(p)$, so we may use the moment generation function of $U_i^k$ to lower bound the availability probabilities. By Lemma~\ref{lem:mgf_gap}, the moment generating function of $V_i^k$ can be a proxy for the moment generating function of~$U_i^k$. For all $i \in \Lcal$, $k \in \Kcal$ and $\ell = 1,\ldots,k$, we define $M_i^k(\ell)= \mathbb E\{ V_i^k \ts | \ts \Dvec^{[\ell,k]}\}$. Therefore, the random variable $M_i^k(\ell)$ is a deterministic function of $\Dvec^{[\ell,k]}$. Noting that $V_i^k$ is a deterministic function of $\Dvec^{[1,k]}$, we have $M_i^k(1)= \mathbb E\{ V_i^k \ts | \ts \Dvec^{[1,k]}\} = V_i^k$ with probability one.

%To upper bound the moment generating function of $V_i^k$, we can upper bound the same for $M_i^k(1)$. 
In the next lemma, we upper bound the moment generating function of $M_i^k(\ell)$ for~all~$\ell = 1,\ldots,k$, which, noting that $M_i^k(1) = V_i^k$, will yield an upper bound on the same for $V_i^k$.

\begin{lem}[Moment Generating Function Bound]
\label{lem:mgf_bound}
Letting $M_i^0(0) = 0$, for all $k \in \Kcal$, $i \in \Lcal$, $\ell = 0,\ldots,k-1$ and $\lambda \geq 0$, we have 
\begin{align*}
\mathbb E \{ e^{\lambda ( M_i^{k-1}(\ell) + \sum_{s=1}^t n_{is}^k(q))} \ts | \ts D^{k-1} = q\} 
~\leq~
e^{\frac{2}{\epsilon^6} \ts (k-1-\ell) \lambda^2 + \gamma \ts c_i \lambda}.
\end{align*}
\end{lem}

\noindent{\it Proof:} We show the result by using induction over $\ell = 1,\ldots,k-1$. Consider the case $\ell = k-1$. We have $\mathbb E\{ M_i^{k-1}(k-1) \ts | \ts D^{k-1} \} = \mathbb E\{\ts \mathbb E\{ V_i^{k-1} \ts | \ts D^{k-1} \} \ts | \ts D^{k-1}\} = \mathbb E\{ V_i^{k-1} \ts | \ts D^{k-1}\}$, so we get 
\begin{align*}
& \mathbb E\{ M_i^{k-1}(k-1) \ts | \ts D^{k-1} = q \} + \sum_{s=1}^t n_{is}^k(q)
~~=~~
\mathbb E\{ V_i^{k-1} \ts | \ts D^{k-1} = q\} + \sum_{s=1}^t n_{is}^k(q)
\\
& \qquad \quad \stackrel{(a)}=~ 
\sum_{\ell=1}^{k-1} \sum_{s \in \Tcal} \sum_{p \in \Tcal} \mathbb P \{D^\ell \geq s,~D^{\ell-1} = p \ts | \ts D^{k-1} = q \} \ts n_{is}^\ell(p) 
+
\sum_{s=1}^t n_{is}^k(q)
\\
& \qquad \quad \stackrel{(b)}=~ \sum_{\ell=1}^{k-1} \sum_{s \in \Tcal} \sum_{p \in \Tcal} \sum_{j \in \Jcal} a_{ij} \ts  \mathbb P \{D^\ell \geq s,~D^{\ell-1} = p \ts | \ts D^{k-1} = q \} \ts \gamma \ts \xbar_{js}^\ell(p) 
+
\sum_{s=1}^t \sum_{j \in \Jcal} a_{ij} \ts \gamma \ts \xbar_{js}^k(q)  
~\ts \stackrel{(c)}\leq ~\ts 
\gamma\ts  c_i,
\end{align*}
where $(a)$ uses the definition of $V_i^{k-1}$, $(b)$ follows because we have $n_{is}^\ell(p) = \mathbb E\{ N_{is}^\ell(p)\}$, in which case,  by the definition of $N_{is}^\ell(p)$, we get $n_{is}^\ell(p) = \mathbb E\{ N_{is}^\ell(p) \} = \sum_{j \in \Jcal} a_{ij} \ts \gamma \ts \xbar_{js}^\ell(p)$ and $(c)$ holds by the first constraint in problem (\ref{eqn:lp}), as well as noting that conditional on $D^{k-1}$, $D^\ell$ is independent of $D^k$ for $\ell \leq k-1$, so  \mbox{$\mathbb P \{ D^\ell \geq s, \ts D^{\ell-1} =p \ts|\ts D^k \geq t, \ts D^{k-1} = q \} = \mathbb P \{ D^\ell \geq s, \ts D^{\ell-1} =p \ts|\ts \ts D^{k-1} = q \}$}.~By its definition, $M_i^{k-1}(k-1)$ is a deterministic function of $D^{k-1}$, which implies that given $D^{k-1} = q$, $M_i^{k-1}(k-1)$ is a deterministic quantity. Thus, using the chain of inequalities above, we~obtain $\mathbb E \{ e^{\lambda ( M_i^{k-1}(k-1) + \sum_{s=1}^t n_{is}^k(q))} \ts | \ts D^{k-1} = q\} = e^{\mathbb E\{ \lambda ( M_i^{k-1}(k-1) + \sum_{s=1}^t n_{is}^k(q)) \ts | \ts D^{k-1} = q\}} \leq e^{\gamma \ts c_i \lambda}$, so the result holds for $\ell = k-1$. Assuming that the result holds for $\ell+1 \leq k-1$, we show that the result holds for \mbox{$\ell \leq k-1$}. By Lemma \ref{lem:martingale_diff} in Extended Results \ref{sec:auxiliary}, $|M_i^{k-1}(\ell) - M_i^{k-1}(\ell+1)| \leq \frac{2}{\epsilon^3}$ with probability one. Also, by the tower property of conditional expectations, using the definition of $M_i^{k-1}(\ell)$, we have $\mathbb E\{ M_i^{k-1}(\ell) \ts | \ts \Dvec^{[\ell+1,k-1]} \} = \mathbb E\{\ts \mathbb E\{ V_i^{k-1} \ts | \ts \Dvec^{[\ell,k-1]}\} \ts | \ts \Dvec^{[\ell+1,k-1]} \} = \mathbb E\{ V_i^{k-1} \ts | \ts  \Dvec^{[\ell+1,k-1]} \}$. Using precisely the same argument, we can verify that $\mathbb E\{ M_i^{k-1}(\ell+1) \ts | \ts \Dvec^{[\ell+1,k-1]} \} = 
%\mathbb E\{\ts \mathbb E\{ V_i^{k-1} \ts | \ts \Dvec^{[\ell+1,k-1]}\} \ts | \ts \Dvec^{[\ell+1,k-1]} \} = 
\mathbb E\{ V_i^{k-1} \ts | \ts  \Dvec^{[\ell+1,k-1]} \}$ as well. Thus, conditional on $\Dvec^{[\ell+1,k-1]}$, the random variable $M_i^{k-1}(\ell) - M_i^{k-1}(\ell+1)$ is mean-zero and bounded by $[-\frac{2}{\epsilon^3},\frac{2}{\epsilon^3}]$. Recall that if the mean-zero random variable $Z$ is bounded by $[a,b]$, then~we have $\mathbb E\{ e^{\lambda Z}\} \leq e^{\frac 18 (b-a)^2 \lambda^2}$. In this case, we get $\mathbb E\{ e^{\lambda( M_i^{k-1}(\ell) - M_i^{k-1}(\ell+1) )} \ts | \ts \Dvec^{[\ell+1,k-1]} \} \leq e^{\frac{2}{\epsilon^6}  \lambda^2}$. Thus, using the fact that $M_i^{k-1}(\ell+1)$ is a deterministic function of $\Dvec^{[\ell+1,k-1]}$, we have
\begin{align*}
&\mathbb E \{ e^{\lambda M_i^{k-1}(\ell)} \ts | \ts D^{k-1} \} 
~\ts = ~\ts
\mathbb E \{\ts \mathbb E\{  e^{\lambda M_i^{k-1}(\ell+1)} \ts e^{\lambda( M_i^{k-1}(\ell) - M_i^{k-1}(\ell+1))} \ts | \ts \Dvec^{[\ell+1,k-1]} \} \ts | \ts D^{k-1} \} 
\\
~& \quad = ~\ts
\mathbb E \{e^{\lambda M_i^{k-1}(\ell+1)} \ts \mathbb E\{  e^{\lambda( M_i^{k-1}(\ell) - M_i^{k-1}(\ell+1))} \ts | \ts \Dvec^{[\ell+1,k-1]} \} \ts | \ts D^{k-1} \} 
~\ts \leq ~\ts
e^{\frac{2}{\epsilon^6}  \lambda^2} \ts \mathbb E\{ e^{\lambda M_i^{k-1}(\ell+1)} \ts | \ts D^{k-1}\}. 
\end{align*}
Thus, we get $\mathbb E \{ e^{\lambda ( M_i^{k-1}(\ell) + \sum_{s=1}^t n_{is}^k(q))} \ts | \ts D^{k-1} = q\} \leq e^{\frac{2}{\epsilon^6}  \lambda^2} \ts \mathbb E\{ e^{\lambda (M_i^{k-1}(\ell+1)+\sum_{s=1}^t n_{is}^k(q))} \ts | \ts D^{k-1} = q\}\leq e^{\frac{2}{\epsilon^6}\lambda^2} \ts  e^{\frac{2}{\epsilon^6} \ts (k-2-\ell) \lambda^2 + \gamma \ts c_i \lambda}$, where the last inequality is by the induction assumption. \qed

Note that $\{ M_i^k(\ell) : \ell = 1,\ldots,k\}$ is a martingale adapted to $\{ \Dvec^{[\ell,k]} : \ell = 1,\ldots,k\}$ in the sense that $\mathbb E \{ M_i^k(\ell) \ts | \ts \Dvec^{[\ell+1,k]} \} = \mathbb E \{ \ts \mathbb E\{ V_i^k \ts | \ts \Dvec^{[\ell,k]} \} \ts | \ts \Dvec^{[\ell+1,k]} \} = \mathbb E\{ V_i^k \ts | \ts \Dvec^{[\ell+1,k]} \}  = M_i^k(\ell+1)$.

{\bf \underline{Performance Guarantee for the Approximate Policy}:}
\\
\indent In the next lemma, we use the moment generating function bounds given in Lemmas \ref{lem:mgf_gap} and \ref{lem:mgf_bound} to lower bound  the availability probability on the right side of (\ref{eqn:union}). 

\begin{lem}[Availability Probability Bound]
\label{lem:avail}
Letting $U_i^0 = 0$, for all $k \in \Kcal$, $i \in \Lcal$, $t,q \in \Tcal$ and $\lambda \in [0,1]$, we have 
\begin{align*}
\mathbb P \bigg\{ U_i^{k-1} + \sum_{s=1}^t N_{is}^k(q) \ts \geq \ts c_i \ts \Big| \ts D^{k-1} = q  \bigg\}
~\leq~ 
e^{(c_i + \frac{3}{\epsilon^6} (k-1) )\lambda^2 - (1-\gamma) c_i \lambda}.
\end{align*}
\end{lem}

\noindent{\it Proof:} By the discussion just after the definition of $M_i^k(\ell)$, we have $V_i^{k-1} = M_i^{k-1}(1)$ with probability one, so $\mathbb E \{ e^{\lambda ( V_i^{k-1} + \sum_{s \in \Tcal} n_{is}^k(q))} \ts | \ts D^{k-1} = q \} = \mathbb E \{ e^{\lambda ( M_i^{k-1}(1) + \sum_{s \in \Tcal} n_{is}^k(q))} \ts | \ts D^{k-1} = q \} \leq e^{\frac{2}{\epsilon^6} \ts (k-2) \lambda^2 + \gamma \ts c_i \lambda}$, where the last inequality uses Lemma \ref{lem:mgf_bound}. On the other hand, by the first constraint in (\ref{eqn:lp}), we have $\sum_{s=1}^t n_{is}^k(q) = \sum_{s=1}^t \sum_{j \in \Jcal} a_{ij} \ts \gamma \ts \xbar_{js}^k(q) \leq \gamma \ts c_i$. By a simple lemma, given as Lemma \ref{lem:mgf_bern} in Extended Results \ref{sec:auxiliary}, if the Bernoulli random variable $Z$ has mean $\mu$, then $\mathbb E\{ e^{\lambda \ts (Z - \mu)} \} \leq e^{\mu \lambda^2}$ for all $\lambda \in [0,1]$, so because $n_{is}^k(q) = \mathbb E\{ N_{is}^k(q)\}$, we get $\mathbb E\{ e^{\lambda \sum_{s=1}^t(N_{is}^k(q) - n_{is}^k(q))} \} \leq e^{\lambda^2 \sum_{s=1}^t n_{is}^k(q)} \leq e^{\gamma \ts c_i \ts \lambda^2}$. We get
%
%for any random variable $Z$, we have $\mathbb P \{ Z \geq c \} = \mathbb P \{ e^{\lambda Z} \geq e^{\lambda c} \} \leq \frac{1}{e^{\lambda c}} \mathbb E\{ e^{\lambda Z} \}$  and $\lambda \geq 0$, we obtain the chain of inequalities 
%
%
\begin{align*}
& \mathbb P \bigg\{ U_i^{k-1} + \sum_{s=1}^t N_{is}^k(q) \ts \geq \ts c_i \ts \Big| \ts D^{k-1} = q  \bigg\}
~=~
\mathbb P\{ e^{\lambda \ts ( U_i^{k-1} + \sum_{s=1}^t N_{is}^k(q) ) } \geq e^{\lambda \ts c_i} \ts | \ts D^{k-1} =q \} 
\\
~& \qquad \qquad \stackrel{(a)}\leq~
\frac{1}{e^{\lambda c_i}} \ts \mathbb E\{ e^{\lambda \ts ( U_i^{k-1} + \sum_{s=1}^t N_{is}^k(q) ) } \ts | \ts D^{k-1} =q \} 
\\
~&\qquad \qquad \stackrel{(b)}=~
\frac{1}{e^{\lambda c_i}} \ts \mathbb E\{ e^{\lambda \sum_{s=1}^t (N_{is}^k(q) - n_{is}^k(q) ) } \} \ts \ts \mathbb E\{ e^{\lambda \ts ( U_i^{k-1} + \sum_{s=1}^t n_{is}^k(q) ) } \ts | \ts D^{k-1} =q \} 
\\
~&\qquad \qquad  \stackrel{(c)}\leq~
\frac{1}{e^{\lambda c_i}} \ts \mathbb E\{ e^{\lambda \sum_{s=1}^t (N_{is}^k(q) - n_{is}^k(q) ) } \} \ts \ts e^{\frac{1}{2\epsilon} (k-1)  \lambda^2} \ts \mathbb E\{ e^{\lambda \ts ( V_i^{k-1} + \sum_{s=1}^t n_{is}^k(q) ) } \ts | \ts D^{k-1} =q \} 
\\
~& \qquad \qquad  \stackrel{(d)}\leq~
\frac{1}{e^{\lambda c_i}} \ts e^{\gamma \ts c_i \ts \lambda^2} \ts \ts e^{\frac{1}{2\epsilon} (k-1)  \lambda^2} \ts 
e^{\frac{2}{\epsilon^6} \ts (k-2) \lambda^2 + \gamma \ts c_i \lambda}
\\
~& \qquad \qquad  =~
e^{(\gamma c_i + \frac{1}{2\epsilon} \ts (k-1) + \frac{2}{\epsilon^6} (k-2))\lambda^2 - (1-\gamma) c_i \lambda }
~\stackrel{(e)}\leq~
e^{(c_i + \frac{3}{\epsilon^6} \ts (k-1))\lambda^2 - (1-\gamma) c_i \lambda },
\end{align*}
where $(a)$ uses the Markov inequality, $(b)$ holds because $N_{is}^k(q)$ is independent of $D^{k-1}$, $(c)$ is by Lemma \ref{lem:mgf_gap}, $(d)$ uses the two inequalities at the beginning of the proof and $(e)$ uses $\epsilon \leq 1$. \qed

Using specific values for $\gamma$ and $\lambda$ in Lemma \ref{lem:avail}, we will bound the availability probabilities. Using this bound in (\ref{eqn:apx_rev}) will yield the performance guarantee in Theorem \ref{thm:perf}.

{\it \underline{Proof of Theorem \ref{thm:perf}}:}
\\
\indent
We use Lemma \ref{lem:avail} with specific values of $\gamma$ and $\lambda$. Letting $\delta = \frac{1}{\epsilon^6}$ for notational brevity, fix $\gammabar = 1 - \frac{\sqrt{4 \ts (c_{\min} + 3 \ts \delta \ts (K-1)) \log c_{\min}}}{c_{\min}}$ and $\lambdabar = \frac{(1-\gammabar) \ts c_i}{2 (c_i + 3 \ts \delta \ts (K-1))}$. For these  values of $\gammabar$ and $\lambdabar$, we have 
\begin{align}
\!
(c_i + 3 \ts \delta  (K-1) )\lambdabar^2 - (1-\gammabar) \ts c_i \lambdabar 
\ts \stackrel{(a)}= 
- \frac{[(1-\gammabar) \ts c_i]^2}{4 (c_i + 3 \ts \delta \ts (K-1))} 
\ts \stackrel{(b)}\leq 
- \frac{[(1-\gammabar) \ts c_{\min}]^2}{4 (c_{\min} + 3 \ts \delta \ts (K-1))}
\ts \stackrel{(c)}= 
- \log c_{\min},
\label{eqn:avail_exp}
\end{align}
where $(a)$ follows by direct computation with the specific value of $\lambdabar$, $(b)$ holds because we can check the first derivative to verify that $\frac{[(1-\gamma) \ts x]^2}{4 (x + 3 \ts \delta \ts (K-1))}$ is increasing in $x$ for $x \geq 0$ and $(c)$ follows by noting that the value of $\gammabar$ satisfies $(1-\gammabar)^2 \ts c_{\min}^2 = 4 \ts (c_{\min} + 3 \ts \delta \ts (K-1)) \log c_{\min}$. Without loss of generality, we can assume that $c_{\min} > 4 \sqrt{(c_{\min} + \delta \ts (K-1)) \log c_{\min}}$. Otherwise, the second term in the max operator in the theorem becomes a negative number and $\frac{\apx}{\Zbar_\lp}$ is trivially lower bounded by a negative number, so the result immediately holds. Therefore, we have \mbox{$c_{\min} > 4 \sqrt{(c_{\min} + \delta \ts (K-1)) \log c_{\min}} \geq \sqrt{4 \ts (c_{\min} + 3 \ts \delta \ts (K-1)) \log c_{\min}}$}. In this case, our choice of $\gammabar$ satisfies $\gammabar \in [0,1]$. If $\gammabar \in [0,1]$, then our choice of $\lambdabar$ satisfies $\lambdabar \in [0,1]$ as well. Therefore, we can use Lemma \ref{lem:avail} with $\gamma = \gammabar$ and $\lambda = \lambdabar$, so noting (\ref{eqn:avail_exp}), we obtain 
\begin{align*}
\mathbb P \bigg\{ U_i^{k-1} + \sum_{s=1}^t N_{is}^k(q) \ts \geq \ts c_i \ts \Big| \ts D^{k-1} = q  \bigg\}
~\leq~
e^{(c_i + \frac{3}{\epsilon^6} (K-1) )\lambdabar^2 - (1-\gammabar) c_i \lambdabar} 
~\leq~
e^{- \log c_{\min}} 
~=~
\frac{1}{c_{\min}}.
\end{align*}
By the definition of $U_i^{k-1}$, the probability on the left side above is the same as the probability on the right side of (\ref{eqn:union}). Using the inequality above on the right side of (\ref{eqn:union}), because $|\Lcal_j| \leq L$, we get $\mathbb P \{ G_{jt}^k = 1 \ts | \ts \Psi_t^k(q) = 1\}  \geq 1 - \frac{L}{c_{\min}}$. By the discussion just after (\ref{eqn:apx_rev}), if $\mathbb P \{ G_{jt}^k = 1 \ts | \ts \Psi_t^k(q) = 1\}  \geq \alpha$, then $\frac{\apx}{\Zbar_\lp} \geq \gamma \ts \alpha$. Thus, using the specific value of $\gammabar$, we get 
\begin{align*}
& \frac{\apx}{\Zbar_\lp} \geq  \Bigg( 1 - \frac{\sqrt{4 \ts (c_{\min} + 3 \ts \delta \ts (K-1)) \log c_{\min}}}{c_{\min}} \Bigg)  \Bigg( 1 - \frac{L}{c_{\min}} \Bigg).
\end{align*}
The result follows by noting that the right side of the chain of inequalities above is lower bounded by  $1 - 4 \frac{\sqrt{(c_{\min} + \delta \ts (K-1)) \log c_{\min}}}{c_{\min}} - \frac{L}{c_{\min}}$, as well as using the fact that $\frac{\apx}{\opt} \geq \frac{\apx}{\Zbar_\lp}$ by Theorem \ref{thm:bound}.  \qed 

In Extended Results \ref{sec:tuning}, we unpack the proof of Theorem~\ref{thm:perf} further by elaborating on the tradeoffs that lead to the specific choices for the values of $\gamma$ and $\lambda$ in the proof.
\color{black}

\vspace{-2mm}

\section{Extensions, Numerical Performance and Research Directions}
\label{sec:conc}

\vspace{-2mm}

We can study performance guarantees when the demand random variables $D^1,\ldots,D^K$ in different stages are independent. In Extended Results \ref{sec:indep}, assuming that the demand random variable in each stage is \mbox{sub-Gaussian} with variance proxy $\sigma^2$ and independent of the demand random variables in other stages, we give an approximate policy with a performance guarantee of \mbox{$1 - \Ommega\Big( \frac{\sqrt{(c_{\min} +  K \ts \sigma^2 ) \log c_{\min}}}{c_{\min}} \Big)$}. In Extended Results \ref{sec:assort}, we show that our fluid approximation naturally extends to the case where  we make assortment offering or pricing decisions. In Extended Results~\ref{sec:exp}, we give computational experiments to demonstrate that policies driven by the correct fluid approximation can make a dramatic impact in practice. 
\color{black}
Our work opens up several research directions. It would be useful to test our model in practical revenue management systems.~The additional input requirement of our model is the distribution of the customer arrivals in the next stage conditional on the number of customer arrivals in the current stage.~Assuming that each stage corresponds to a week, one approach to estimate this input requirement is that we start by unconstraining the customer arrivals so that we obtain the raw data that corresponds to the volume of the customer arrivals in each week; see  \cite{WePo02}, \cite{QuFe07} and \cite{KoLi19}. Obtaining such unconstrained customer arrival volumes is often done in revenue management systems. 
%In this step, one may need the volume of the customer arrivals in each week for different markets separately.~In the airline setting, for example, one may need the volume of the customer arrivals in each week separately for  each origin and destination pair.
Next, we can fit an autoregressive process to the customer arrival data; see \cite{BoJe15}. Lastly, there is work on capturing the evolution of an autoregressive process through a Markov chain, allowing us to capture the distribution of the demand at the next stage as a function of the demand in the current stage; see \cite{Ta86}, \cite{Bu01} and \cite{FeGa19}. This approach is one possibility for estimating the key input requirement of our model and it is useful to test such approaches in practice. Another research direction is that the demand in a stage depends on the demand in the previous stage in our model. One can focus on more complicated dependence structures. 

{\bf Acknowledgements and Notes:} The authors thank the department editor, associate editor and two referees whose comments significantly improved the depth and exposition in the paper. Extended results for our paper are available in \cite{LiRuTo24}.

\vspace{-4mm}

\bibliographystyle{ormsv080}
\renewcommand*{\bibfont}{\footnotesize}
\renewcommand*{\bibfont}{\normalfont\footnotesize\linespread{1}\selectfont}
\setlength{\bibsep}{-2.1pt}
\bibliography{references}

\clearpage

\thispagestyle{empty}

~\\
~\\
~\\
\begin{center}

\large
Details and Extended Results: \\  Technical Note -- Revenue Management with Calendar\hspace{0.5mm}-Aware and Dependent Demands: Asymptotically Tight Fluid Approximations

~\\
~\\

\normalfont \normalsize Weiyuan Li, Paat Rusmevichientong, Huseyin Topaloglu 
\\
\today

\end{center}

\newpage

\setcounter{page}{1}

\begin{APPENDICES}

%\ECSwitch
%\ECDisclaimer
%%% Main head for the e-companion
%\ECHead{
%\vspace{-1.5cm}
%\begin{center}
%Appendix: Fluid Approximations for Calendar\hspace{0.5mm}-Aware Network Revenue Management with Dependent Demands
%\\
%{\normalfont \normalsize Weiyuan Li, Paat Rusmevichientong, Huseyin Topaloglu \\ \vspace{-2mm} \today }
%\end{center}
%}

\vspace{-1.0cm}
\begin{center}
{\large
\underline{Details and Extended Results}: \\  
Technical Note -- Revenue Management with Calendar\hspace{0.5mm}-Aware and Dependent Demands: Asymptotically Tight Fluid Approximations}
\\
{\normalfont \normalsize Weiyuan Li, Paat Rusmevichientong, Huseyin Topaloglu \\ %\vspace{-2mm} \today 
}
\end{center}

\vspace{-0.2cm}

\section{Matching a Given Total Demand Distribution}
\label{sec:total_demand}

Consider a target distribution for the total demand over the selling horizon with a support bounded by $KT$. We show that we can calibrate our demand model such that the distribution of the total demand implied by our demand model matches the target distribution. In our demand model, there are $K$ stages and the demand in each stage takes values between one and $T$, so the total demand takes values in $\{K,\ldots,KT\}$. The assumption that the demand in each stage is at least one is without loss of generality. We can add a dummy time period to each stage such that the probability that there is a request for any product at the dummy time period is zero, in which case, there is at least one time period in each stage, but we may not have any demand in a particular stage.~Nevertheless, for consistency with our demand model, we consider an arbitrary total demand random variable $D^\total$ that takes values in $\{K,\ldots,KT\}$. We show that we can calibrate our demand model such that the distribution of the total demand implied by our demand model matches the distribution of $D^\total$. Assuming that we are given the random variable $D^\total$ for the total demand, for stage $k$, we define the probability mass function over  $\{1,\ldots,T\}$ as  
\begin{align}
\!\!\!\!\!\!\!
f^k(\ell) = 
\begin{cases}
\mathbb P \bigl\{ D^\total = (k-1) \ts (T-1) + K + \ell - 1 \ts | \ts D^\total \geq (k-1)(T-1) + K \bigr\}
& \mbox{if $\ell = 1,\ldots,T-1$}
\\
\mathbb P \bigl\{ D^\total \geq (k-1) \ts (T-1) + K + T - 1 \ts | \ts D^\total \geq (k-1)(T-1) + K \bigr\}
& \mbox{if $\ell = T$.}
\end{cases}
\!\!\!
\label{eqn:pmf}
\end{align}

Note that we have $\sum_{\ell =1}^T f^k(\ell) = 1$, so $f^k$ is indeed a probability mass function. The way we use the probability mass function shortly becomes clear, but to provide intuition into this probability mass function, note that the random variable $D^\total - K$ takes values in $\{0,\ldots,K(T-1)\}$. In our model, on the other hand, the demand in each stage takes at least a value of one. In addition to this value of one, the demand in each stage can contribute at most $T-1$ units to the total demand. In this case, we observe that if the condition $D^\total \geq (k-1)(T-1) + K$ in the conditional probability in (\ref{eqn:pmf}) holds, then we have \mbox{$D^\total - K \geq (k-1) \ts (T-1)$}, which implies that the total demand surely exceeds the demands in the first $k-1$ stages of our model. In the rest of this section, we will construct the random variables $D^1,\ldots,D^K$ that satisfy the following three properties. First, each of the random variables $D^1,\ldots,D^K$ takes values in $\{1,\ldots,T\}$. Second, conditional on $D^k$, the random variable $D^{k+1}$ is independent of the random variables $D^1,\ldots,D^{k-1}$. Third, the random variable $D^1+\ldots+D^K$ has the same distribution as the random variable $D^\total$. We accomplish this goal as follows. We use $Z^1,\ldots,Z^K$ to denote a sequence of independent random variables, where $Z^k$ has the probability mass function $f^k$. Using the random variables $Z^1,\ldots,Z^K$, we define the random variables $D^1,\ldots,D^K$ recursively as 
\begin{align}
D^k =
\begin{cases}
1 & \mbox{if $D^{k-1} \neq T$}
\\
Z^k & \mbox{if $D^{k-1} = T$,}
\end{cases}
\label{eqn:rec_vars}
\end{align}
where we fix $D^0 =T$. Note that $D^k$ takes values in $\{1,\ldots,T\}$. Also, $D^k$ depends only on $D^{k-1}$ and $Z^k$. Thus, since the random variables $Z^1,\ldots,Z^K$ are independent of each other, given $D^{k-1}$,~the random variable $D^k$ is independent of $D^1,\ldots,D^{k-2}$. In the next theorem, we show that the random variable \mbox{$D^1+ \ldots+ D^K$} has the same distribution as $D^\total$. Before we state this theorem, by  (\ref{eqn:rec_vars}), note that only one of the random variables $D^1,\ldots,D^K$ can take a value other than one or $T$. Also, if the random variable $D^k$ takes a value other than one or $T$, then the random variables $D^1,\ldots,D^{k-1}$ all take value $T$ and the random variables $D^{k+1}\ldots,D^K$ all take value one.

\begin{thm}
Letting the random variables $D^1,\ldots,D^K$ be recursively defined as in $(\ref{eqn:rec_vars})$, the random \mbox{$D^1+\ldots+D^K$} has the same distribution as the random variable $D^\total$. 
\end{thm}

\noindent{\it Proof:} Fixing some $\ell \in \{K,\ldots,KT\}$, we show that $\mathbb P \{ D^1+ \ldots + D^K = \ell\} = \mathbb P \{ D^\total = \ell\}$.  Using $\ru{\cdot}$ to denote the round up function, set $r = \ruBig{\frac{\ell - K}{T-1}}$, which takes values in $\{0,\ldots,K\}$. By the discussion just before the theorem, only one of the random variables in the sequence $D^1,\ldots,D^K$ can take a value other than one or $T$. If $D^k$ is such a random variable, then all of the random variables $D^1,\ldots,D^{k-1}$ take value $T$ and all of the random variables $D^{k+1},\ldots,D^K$ take value one. In this case, if $D^1+\ldots + D^K = \ell$, then we must have $D^1 = T,\ldots, D^{r-1} = T$, $D^{r+1} = 1,\ldots, D^K = 1$ and $D^r = \ell - (r-1) \ts T - (K-r)$. Considering the case $\ell - (r-1) \ts T - (K-r) \neq T$, we have
\vspace{-1mm}
\begin{align*}
& \mathbb P \bigl\{ D^1+\ldots + D^K = \ell \bigr\} 
\\
& \quad =~
\mathbb P \bigl\{ D^1 = T,\ldots, D^{r-1} = T ,~ D^r = \ell - (r-1) \ts T - (K-r) ,~ D^{r+1} = 1,\ldots, D^K = 1 \bigr\}
\\
& \quad \stackrel{(a)} =~ 
\mathbb P \bigl\{ D^1 = T \bigr\} \times \mathbb P \bigl\{ D^2 = T \ts | \ts D^1 = T \bigr\} \times \ldots \times \mathbb P \bigl\{ D^r = \ell - (r-1) \ts T - (K-r) \ts | \ts D^{r-1} = T \bigr\}
\\
& \quad \quad \times \mathbb P \bigl\{ D^{r+1} = 1 \ts | \ts D^r =  \ell - (r-1) \ts T - (K-r) \bigr\} 
\\
& \quad  \quad \times
\mathbb P \bigl\{ D^{r+2} = 1 \ts | \ts D^{r+1} = 1 \bigr\} \times \ldots \times \mathbb P \{ D^K  =1 \ts | \ts D^{K-1} = 1 \bigr\}
\\
& \quad \stackrel{(b)} =~
\mathbb P \bigl\{ Z^1 = T \bigr\} \times \mathbb P \bigl\{ Z^2 = T \bigr\} \times \ldots \times \mathbb P \bigl\{ Z^r = \ell - (r-1) \ts T - (K-r) \bigr\}
\\
& \quad \stackrel{(c)}=~
\mathbb P \bigl\{ D^\total \geq T - 1 + K \ts | \ts D^\total \geq K  \bigr\}
\times
\mathbb P \bigl\{ D^\total \geq  2(T - 1) + K  \ts | \ts D^\total \geq T-1 + K  \bigr\}
\\
& \quad \quad   \times \ldots \times 
\mathbb P \bigl\{ D^\total = (r-1) \ts (T-1) + K + \ell - (r-1) \ts T - (K-r) - 1 \ts | \ts D^\total \geq (r-1) \ts (T-1) + K  \bigr\}
\\
& \quad \stackrel{(d)} = ~ \mathbb P \bigl\{ D^\total = \ell,~ D^\total \geq (r-1) \ts (T-1) + K \ts | \ts D^\total \geq K \bigr\}
~\stackrel{(e)} =~
\mathbb P \bigl\{ D^\total = \ell \bigr\},
\end{align*}
where $(a)$ holds because conditional on $D^k$, the random variable $D^{k+1}$ is independent of the random variables $D^1,\ldots,D^{k-1}$, $(b)$ uses (\ref{eqn:rec_vars}) along with the fact that $\ell - (r-1) \ts T - (K-r) \neq T$, $(c)$ follows from (\ref{eqn:pmf}), $(d)$ uses the Bayes rule and $(e)$ holds because the support of $D^\total$ is lower bounded by $K$, as well as the fact that $r \leq \frac{\ell - K}{T-1}$, which implies that $\ell \geq r \ts (T-1) + K$. Considering the case $\ell - (r-1) \ts T - (K-r) = T$, we can follow an argument similar to the one in this paragraph to show that we have $\mathbb P \{ D^1 + \ldots + D^K = \ell \} = \mathbb P \{ D^\total = \ell \}$ in this case as well. \qed

\section{Proof of Theorem \ref{thm:bound}}
\label{sec:bound}

We relax the constraint $a_{ij} \ts u_j \leq y_i$ at time period $t$ in stage $k$ in (\ref{eqn:dp}) using the Lagrange multiplier $\mu_{it}^k(q)$. Letting $\muvec  = (\mu_{it}^k(q) : i \in \Lcal,~t,q \in \Tcal,~k \in \Kcal)$, we obtain the relaxed dynamic program
\begin{align}
\Jtilde_t^k(\yvec,q; \muvec)
& =
\max_{\uvec \in \{0,1\}^{|\Jcal|}} \Bigg\{ \sum_{j \in \Tcal} \lambda_{jt}^k \ts \Big\{ f_j \ts u_j + \theta_t^k(q) \ts \Jtilde_{t+1}^k(\yvec -\avec_j \ts u_j ,q; \muvec) 
\nonumber
\\
& \qquad \qquad \qquad \qquad ~~~ + ( 1 - \theta_t^k(q)) \ts \Jtilde_1^{k+1}(\yvec - \avec_j   \ts u_j,t;\muvec) \Big\}
+ 
\sum_{i \in \Lcal} \sum_{j \in \Jcal} \lambda_{jt}^k \ts \mu_{it}^k(q) \ts[y_i - a_{ij} \ts u_j ] \Bigg\}
\nonumber
\\
& =
\max_{\uvec \in \{0,1\}^{|\Jcal|}} \Bigg\{ \sum_{j \in \Tcal} \lambda_{jt}^k \ts \Big\{ \Big[ f_j - \sum_{i \in \Lcal} a_{ij} \ts \mu_{it}^k(q) \Big] \ts u_j + \theta_t^k(q) \ts \Jtilde_{t+1}^k(\yvec -\avec_j \ts u_j ,q; \muvec) 
\nonumber
\\
& \qquad \qquad \qquad \qquad ~~~ + ( 1 - \theta_t^k(q)) \ts \Jtilde_1^{k+1}(\yvec - \avec_j   \ts u_j,t;\muvec) \Big\} \Bigg\}
+ 
\sum_{i \in \Lcal} \mu_{it}^k(q) \ts y_i,
\label{eqn:relax_dp}
\end{align}
with the boundary condition that $\Jtilde_1^{K+1} = 0$. Note that the value functions of the relaxed dynamic program depend on the choice of the Lagrange multipliers. In the first equality above, we scale the Lagrange multiplier $\mu_{it}^k(q)$ with $\lambda_{jt}^k$ for notational uniformity. The second equality follows by arranging the terms and using the fact that $\sum_{j \in \Jcal} \lambda_{jt}^k=1$. If the Lagrange multipliers are \mbox{non-negative}, then the value functions from the relaxed dynamic program in (\ref{eqn:relax_dp}) are upper bounds on the value functions from the dynamic program in (\ref{eqn:dp}). We do not show this result. This result is considered standard and analogues of this result have been shown in other settings; see Proposition~2 in \cite{AdMe04}.  Therefore, we have $\Jtilde_t^k(\yvec,q; \muvec) \geq J_t^k(\yvec,q)$ for all $\yvec \in \mathbb Z_+^{|\Lcal|}$, $q \in \Tcal$ as long as $\muvec \in \mathbb R_+^{|\Lcal|T^2K}$. We can solve the problem $\min_{\muvec \in \mathbb R_+^{|\Lcal|T^2K}} \Jtilde_1^1(\cvec,\Dhat^0;\muvec)$ to obtain an upper bound on the optimal total expected revenue. One of the useful features of the relaxed dynamic program is that the value functions computed through this dynamic program are linear in the remaining capacities. In the next lemma, we show that $\Jtilde_t^k(\yvec,q;\muvec) = \sum_{i \in \Lcal} \alpha_{it}^k(q;\muvec) \ts y_i + \beta_t^k(q;\muvec)$, where the slope $\alpha_{it}^k(q;\muvec)$ and the intercept $\beta_t^k(q;\muvec)$ are recursively computed as 
\begin{gather}
\alpha_{it}^k(q;\muvec) ~=~ \mu_{it}^k(q) + \theta_t^k(q) \ts \alpha_{i,t+1}^k(q;\muvec) +  (1-\theta_t^k(q)) \ts \alpha_{i1}^{k+1}(t;\muvec)
\label{eqn:slope}
\\
\beta_t^k(q;\muvec) = \sum_{j \in \Jcal} \lambda_{jt}^k \ts \Big[ f_j - \sum_{i \in \Lcal} a_{ij} \ts \alpha_{it}^k(q;\muvec) \Big]^+
+
\theta_t^k(q) \ts \beta_{t+1}^k(q;\muvec) +  (1-\theta_t^k(q)) \ts \beta_1^{k+1}(t;\muvec),
\nonumber
\end{gather}
with the boundary condition that $\alpha_{i1}^{K+1} = 0$ and $\beta_1^{K+1} = 0$. The linear form of the value functions from the relaxed dynamic program will be useful to show Theorem \ref{thm:bound}.

\begin{lem}
Letting $\alpha_{it}^k(q;\muvec)$ and $\beta_t^k(q;\muvec)$ be as in $(\ref{eqn:slope})$, the value functions computed through the dynamic program in $(\ref{eqn:relax_dp})$ satisfy $\Jtilde_t^k(\yvec,q;\muvec) = \sum_{i \in \Lcal} \alpha_{it}^k(q;\muvec) \ts y_i + \beta_t^k(q;\muvec)$ for all $t \in \Tcal$ and $k \in \Kcal$. 
\end{lem}

\noindent{\it Proof:} We show the result by using induction over the time periods. At the last time period in the last stage, by (\ref{eqn:relax_dp}), we have $J_T^K(\yvec,q;\muvec) = \sum_{j \in \Tcal} \lambda_{jT}^K \ts [ f_j - \sum_{i \in \Lcal} a_{ij} \ts \mu_{iT}^K(q) ]^+ + \sum_{i \in \Lcal} \mu_{iT}^K(q) \ts y_i = \beta_T^K(q;\muvec) + \sum_{i \in \Lcal} \alpha_{iT}^K(q;\muvec) \ts y_i$, where the last equality uses (\ref{eqn:slope}). Therefore, the result holds at the last time period in the last stage. Assuming that the result holds at all time periods after time period~$t$ in stage $k$, we show that the result holds at time period $t$ in stage $k$ as well.~Using the induction assumption on the right side of (\ref{eqn:relax_dp}), we have 
\begin{align*}
\Jtilde_t^k(\yvec,q; \muvec)
& =  \!\!\!\!\!
\max_{\uvec \in \{0,1\}^{|\Jcal|}}\!\! \Bigg\{ \sum_{j \in \Tcal}\! \lambda_{jt}^k \bigg\{\! \Big[ f_j - \sum_{i \in \Lcal} a_{ij} \ts \mu_{it}^k(q) \Big] u_j 
+ 
\theta_t^k(q) \Big[ \beta_{t+1}^k(q;\muvec) \!+\!\! \sum_{i \in \Lcal} \alpha_{i,t+1}^k(q;\muvec)\ts (y_i -a_{ij} \ts u_j ) \Big]
\\
& \qquad \qquad  ~~~ + ( 1 - \theta_t^k(q)) \Big[ \beta_1^{k+1}(t;\muvec) \!+\!\! \sum_{i \in \Lcal} \alpha_{i1}^{k+1}(t;\muvec) \ts (y_i - a_{ij}   \ts u_j)   \Big]\! \bigg\} \! \Bigg\}
+ 
\sum_{i \in \Lcal} \mu_{it}^k(q) \ts y_i
\\
& \stackrel{(a)}= \!\!\!\!\!
\max_{\uvec \in \{0,1\}^{|\Jcal|}}\!\! \Bigg\{ \sum_{j \in \Tcal}\! \lambda_{jt}^k  \Big[ f_j - \sum_{i \in \Lcal} a_{ij} \ts \alpha_{it}^k(q;\muvec) \Big] u_j \Bigg\}
\\
& \qquad \qquad  ~~~  +
\theta_t^k(q) \ts \beta_{t+1}^k(q;\muvec) +  (1-\theta_t^k(q)) \ts \beta_1^{k+1}(t;\muvec) +
\sum_{i \in \Lcal} \alpha_{it}^k(q;\muvec) \ts y_i \phantom{\Bigg\}}
\\
&= 
\sum_{j \in \Tcal}\! \lambda_{jt}^k  \Big[ f_j - \sum_{i \in \Lcal} a_{ij} \ts \alpha_{it}^k(q;\muvec) \Big]^+
\!\! \!\!+
\theta_t^k(q) \ts \beta_{t+1}^k(q;\muvec) +  (1-\theta_t^k(q)) \ts \beta_1^{k+1}(t;\muvec) \!+\!\!
\sum_{i \in \Lcal} \alpha_{it}^k(q;\muvec) \ts y_i \phantom{\Bigg\}}
\\
&\stackrel{(b)}= \beta_t^k(q;\muvec) + \sum_{i \in \Lcal} \alpha_{it}^k(q;\muvec) \ts y_i,  \phantom{\Bigg\}}
\end{align*}
where $(a)$ follows by arranging the terms and using the definition of $\alpha_{it}^k(q;\muvec)$, as well as noting the fact that $\sum_{j \in \Jcal} \lambda_{jt}^k = 1$, whereas $(b)$ uses the definition of $\beta_t^k(q;\muvec)$. \qed

By the lemma above, we have $\Jtilde_1^1(\cvec,\Dhat^0;\muvec) = \sum_{i \in \Lcal} \alpha_{i1}^1(\Dhat^0;\muvec) \ts c_i + \beta_1^1(\Dhat^0,\muvec)$. In this case, the problem $\min_{\muvec \in \mathbb R_+^{|\Lcal|T^2K}} \Jtilde_1^1(\cvec,\Dhat^0;\muvec)$ is equivalent to the linear program
\begin{align}
& \min_{(\alphavec,\betavec,\muvec,\etavec) \in \mathbb R^{T^2K (|\Lcal|+1)} \times \mathbb R_+^{T^2K (|\Lcal| + |\Jcal|)}}~~~  \sum_{i \in \Lcal} \alpha_{i1}^1(\Dhat^0) \ts c_i + \beta_1^1(\Dhat^0)
\label{eqn:tighten_bound}
\\
& \qquad \qquad  \qquad \mbox{st}~~~ 
\alpha_{it}^k(q) ~=~ \mu_{it}^k(q) + \theta_t^k(q) \ts \alpha_{i,t+1}^k(q) +  (1-\theta_t^k(q)) \ts \alpha_{i1}^{k+1}(t) \qquad \forall \ts i \in \Lcal,~t,q \in \Tcal,~k \in \Kcal \phantom{\sum_{i \in \Lcal}}
\nonumber
\\
& \qquad \qquad  \qquad ~~~~\ts\ts
\beta_t^k(q) ~=~ \sum_{j \in \Jcal} \lambda_{jt}^k \ts \eta_{jt}^k(q)
+
\theta_t^k(q) \ts \beta_{t+1}^k(q) +  (1-\theta_t^k(q)) \ts \beta_1^{k+1}(t)
\qquad \forall \ts t,q \in \Tcal,~k \in \Kcal
\nonumber
\\
& \qquad \qquad  \qquad ~~~~\ts\ts
\eta_{jt}^k(q) ~\geq~ f_j - \sum_{i \in \Lcal} a_{ij} \ts \alpha_{it}^k(q) \qquad \forall \ts j \in \Jcal,~t,q \in \Tcal,~k \in \Kcal,
\nonumber
\end{align}
where we use the decision variables $\alphavec = (\alpha_{it}^k(q) : i \in \Lcal,~t,q \in \Tcal,~k \in \Kcal)$, \mbox{$\betavec = (\beta_t^k(q) : q,t \in \Tcal,~k \in \Kcal)$}, $\muvec  = (\mu_{it}^k(q) : i \in \Lcal,~t,q \in \Tcal,~k \in \Kcal)$ and $\etavec = (\eta_{jt}^k(q) : j \in \Jcal,~t,q \in \Tcal,~k \in \Kcal)$. We follow the convention that $\alpha_{i1}^{K+1}(q) = 0$ and $\beta_1^{K+1}(q) = 0$ for all $i \in \Lcal$ and $q \in \Tcal$. In the linear program above, the first constraint computes the slopes of the value functions of the relaxed dynamic program. By the third constraint, noting  the non-negativity constraints, we have \mbox{$\eta_{jt}^k(q) = [f_j - \sum_{i \in \Lcal} a_{ij} \ts \alpha_{it}^k(q)]^+$} at an optimal solution to the linear program, in which case, the second constraint computes the intercepts of the value functions of the relaxed dynamic program. We write the objective function of problem (\ref{eqn:tighten_bound}) as $\sum_{q \in \Tcal} \sum_{i \in \Lcal} {\bf 1}(\Dhat^0 = q) \ts c_i \ts \alpha_{i1}^1(q) + \sum_{q \in \Tcal} {\bf 1}(\Dhat^0 = q) \ts \beta_1^1(q)$. We work with the dual of problem~(\ref{eqn:tighten_bound}).~We associate the dual variables \mbox{$\wvec = (w_t^k(q) : t,q \in \Tcal,~k \in \Kcal)$} with the second constraint in (\ref{eqn:tighten_bound}). The decision variables \mbox{$(\beta_t^k(q) : t,q \in \Tcal,~k \in \Kcal)$} appear only in the second constraint in problem (\ref{eqn:tighten_bound}), so in the dual of problem (\ref{eqn:tighten_bound}), the constraints associated with the decision variables  \mbox{$(\beta_t^k(q) : t,q \in \Tcal,~k \in \Kcal)$} are given by $w_1^1(q) \ts = \ts {\bf 1}(\Dhat^0 = q)$ for all~\mbox{$q \in \Tcal$}, \mbox{$w_1^k(q) = \sum_{p \in \Tcal} (1 - \theta_q^{k-1}(p)) \ts w_q^{k-1}(p)$} for all $q \in \Tcal$, $k \in \Kcal \setminus \{1\}$ and $w_t^k(q) = \theta_{t-1}^k(q) \ts w_{t-1}^k(q)$ for all~\mbox{$t \in \Tcal \setminus \{1\}$}, $q \in \Tcal$, $k \in \Kcal$. We capture these constraints by defining the set
\begin{multline*}
\Wcal ~=~ \Bigg\{ \wvec \in \mathbb R^{T^2K} ~:~ 
w_t^k(q) = \theta_{t-1}^k(q) \ts w_{t-1}^k(q) ~~ \forall \ts t \in \Tcal \setminus \{1\},~q \in \Tcal,~k \in \Kcal,
\\
w_1^k(q) = \sum_{p \in \Tcal} (1 - \theta_q^{k-1}(p)) \ts w_q^{k-1}(p)~~ \forall \ts q \in \Tcal,~k \in \Kcal \setminus \{1\},~~ w_1^1(q) = {\bf 1} (\Dhat^0 =q) ~~\forall \ts q \in \Tcal \Bigg\}.
\end{multline*}
In the next lemma, we show that if we have $\wvec \in \Wcal$, then the vector $\wvec$ is closely related to the joint distribution of demands in a pair of successive stages.

\begin{lem}
\label{lem:lp_prob}
If $\wvec \in \Wcal$, then we have $w_t^k(q) = \mathbb P \{ D^k \geq t,~D^{k-1} = q\}$ for all $t,q \in \Tcal$ and $k \in \Kcal$.
\end{lem}

\noindent{\it Proof:} We show the result by using induction over the time periods. At the first time period in the first stage, we have $w_1^1(q) \ts = \ts {\bf 1}(\Dhat^0 = q) \ts= \ts \mathbb P\{ D^0 = q\} \ts = \ts \mathbb P \{ D^1 \geq 1,~D^0 = q\}$, where the first equality holds by noting the third constraint in the definition of $\Wcal$, the second equality follows by noting that $D^0$ is deterministically fixed at $\Dhat^0$ and the third equality holds because the support of $D^1$ is $\{1,\ldots,T\}$. Assuming that the result holds at all time periods before time period $t$ in stage $k$, we show that the result holds at time period $t$ in stage $k$ as well. If $t \neq 1$, then using the first constraint in the definition of $\Wcal$, we have the chain of equalities $w_t^k(q) \ts = \ts  \theta_{t-1}^k(q) \ts w_{t-1}^k(q) \ts = \ts  \mbox{$\mathbb P \{ D^k \geq t \ts | \ts D^k \geq t-1,~D^{k-1} = q\}$} \ts\ts  \mathbb P \{ D^k \geq t-1 ,~D^{k-1} = q\} \ts= \ts \mbox{$\mathbb P \{ D^k \geq t ,~D^{k-1} = q\}$}$,  where the second equality is by the definition of $\theta_{t-1}^k(q)$ and the induction assumption. Similarly, if $t = 1$, then using the second constraint in the definition of $\Wcal$, we have $w_1^k(q) = \sum_{p \in \Tcal} (1 - \theta_q^{k-1}(p)) \ts w_q^{k-1}(p) = \sum_{p \in \Tcal} \mathbb P \{ D^{k-1} = q \ts | \ts D^{k-1} \geq q,~D^{k-2} = p\} \ts \mathbb P \{ D^{k-1} \geq q,~D^{k-2} = p\}$, but the last sum expression is equal to~$\mathbb P \{ D^{k-1} =q\}$, so $w_1^k(q) = \mathbb P \{ D^{k-1} =q\} = \mathbb P \{ D^k \geq 1,~D^{k-1} = q\}$. \qed

By the lemma above, there exists a single element in $\Wcal$. To write the dual of problem (\ref{eqn:tighten_bound}), we associate the dual variables $\yvec = ( y_{it}^k(q) : i \in \Lcal,~t,q \in \Tcal,~k \in \Kcal)$, $\wvec = (w_t^k(q) : t,q \in \Tcal,~k \in \Kcal)$ and $\uvec = (u_{jt}^k(q) : j \in \Jcal,~t,q \in \Tcal,~k \in \Kcal\}$ with the first, second and third constraints, respectively, in problem (\ref{eqn:tighten_bound}). In this case, the dual of problem (\ref{eqn:tighten_bound}) is given by
\begin{align}
& \max_{(\yvec,\uvec,\wvec) \in \mathbb R_+^{T^2K (|\Lcal| + |\Jcal|)} \times \Wcal} ~~~ \sum_{k \in \Kcal} \sum_{t \in \Tcal} \sum_{q \in \Tcal} \sum_{j \in \Jcal} f_j \ts u_{jt}^k(q)
\label{eqn:lp_long}
\\
& \qquad \qquad   \mbox{st}~~~  
y_{it}^k(q) + \sum_{j \in \Jcal} a_{ij} \ts u_{jt}^k(q) ~=~ \theta_{t-1}^k(q) \ts y_{i,t-1}^k(q) \qquad \forall \ts i \in \Lcal,~ t \in \Tcal \setminus \{1\},~q \in \Tcal,~k \in \Kcal
\nonumber
\\
& \qquad \qquad ~~~~\ts\ts
y_{i1}^k(q) + \sum_{j \in \Jcal} a_{ij} \ts u_{j1}^k(q) ~=~ \sum_{p \in \Tcal} (1 - \theta_q^{k-1}(p)) \ts\ts y_{iq}^{k-1}(p) \qquad \forall \ts i \in \Lcal,~q \in \Tcal,~k \in \Kcal \setminus \{1\}
\nonumber
\\
& \qquad \qquad ~~~~\ts\ts
y_{i1}^1(q) + \sum_{j \in \Jcal} a_{ij} \ts u_{j1}^1(q) ~ =~ {\bf 1}(\Dhat^0 = q) \ts c_i  \qquad \forall \ts i \in \Lcal,~ q \in \Tcal
\nonumber
\\
& \qquad \qquad ~~~~\ts\ts
u_{jt}^k(q) \leq \lambda_{jt}^k \ts w_t^k(q) \qquad \forall \ts j \in \Jcal,~t,q \in \Tcal,~k \in \Kcal, \phantom{\sum_{j \in \Jcal}}
\nonumber
\end{align}
where the constraints above are associated with the decision variables $\alphavec$ and $\etavec$ in (\ref{eqn:tighten_bound}). The constraint for the decision variables $\muvec$ translates into the non-negativity constraint for $\yvec$. 

In (\ref{eqn:lp_long}), we capture the constraint associated with the decision variables $\betavec$ as $\wvec \in \Wcal$. In the next lemma, we give an equality that is satisfied by all feasible solutions to problem (\ref{eqn:lp_long}).

\vspace{-0.5mm}

\begin{lem}
\label{lem:lp_capacity}
Letting $( \yvec,\uvec,\wvec)$ be a feasible solution to the linear program in $(\ref{eqn:lp_long})$, for all $i \in \Lcal$, $t,q \in \Tcal$ and $k \in \Kcal$, we have
\begin{multline*}
\mathbb P \{ D^k \geq t,~D^{k-1}=q\} \ts c_i - \sum_{\ell=1}^{k-1} \sum_{s \in \Tcal} \sum_{p \in \Tcal} \sum_{j \in \Jcal} a_{ij} \ts \mathbb P \{ D^k \geq t ,~D^{k-1} = q \ts | \ts D^\ell \geq s,~D^{\ell-1} = p \} \ts u_{js}^\ell(p)
\\
-
\sum_{s=1}^t \sum_{j \in \Jcal} a_{ij} \ts \mathbb P \{ D^k \geq t \ts | \ts D^k \geq s,~D^{k-1} = q\} \ts u_{js}^k(q) ~=~ y_{it}^k(q).
\end{multline*}

\end{lem}

\noindent{\it Proof:} We show the result by using induction over the time periods. At the first time period in the first stage, by the third constraint in (\ref{eqn:lp_long}), we have $y_{i1}^1(q) = {\bf 1}(\Dhat^0 = q) \ts c_i - \sum_{j \in \Jcal} a_{ij} \ts u_{j1}^1(q) = \mbox{$\mathbb P \{ D^1 \geq 1,~D^0 = q\}$} \ts c_i - \sum_{j \in \Jcal} a_{ij} \ts \mathbb P\{ D^1\geq 1 \ts | \ts D^1 \geq 1,~D^0 = q\} \ts u_{j1}^1(q)$, where the last equality holds because $D^0$ is a deterministic quantity and the support of $D^1$ is $\{1,\ldots,T\}$. Assuming that the result holds at all time periods up to and including time period $t$ in stage $k$, we show that the~result holds at the subsequent time period as well. Consider the case $t\neq T$. We will use three identities.~First, for $\ell \leq k-1$, given $D^{k-1}$, $D^k$ is independent of $D^1,\ldots,D^\ell$, in which case, we obtain 
\begin{align*}
& \theta_t^k(q) \ts\ts 
\mathbb P \{ D^k \geq t,~D^{k-1} = q \ts | \ts D^\ell \geq s,~ D^{\ell-1} = p\}
\\
& \qquad \quad  ~=~ \mathbb P \{ D^k \geq t+1 \ts | \ts D^k \geq t,~ D^{k-1} = q\} \ts \ts \ts
\mathbb P \{ D^k \geq t,~D^{k-1} = q \ts | \ts D^\ell \geq s,~ D^{\ell-1} = p\}
\\
& \qquad \qquad \qquad \qquad \qquad \qquad \qquad \qquad \qquad \quad =~
\mathbb P \{ D^k \geq {t+1},~D^{k-1} = q \ts | \ts D^\ell \geq s,~ D^{\ell-1} = p\}.
\end{align*}
Second, by the Bayes rule and definition of $\theta_t^k(q)$, we can show that $\theta_t^k(q) \ts\ts \mbox{$\mathbb P \{ D^k \geq t, \ts D^{k-1}\! = q\}$} =  \mbox{$\mathbb P \{ D^k \geq t+1, \ts D^{k-1} \!= q\}$}$. 
Third, for $s \leq t$, we can, once more, use the Bayes rule and definition of $\theta_t^k(q)$ to show that $\theta_t^k(q) \ts\ts \mathbb P \{ D^k \geq t \ts | \ts D^k \geq s,~D^{k-1} = q\} \ts = \ts \mathbb P \{ D^k \geq t+1 \ts | \ts D^k \geq s,~D^{k-1} = q\}$. Noting that the solution $(\yvec,\uvec,\wvec)$ is feasible to (\ref{eqn:lp_long}), it satisfies the first constraint. Thus, we obtain 
\begin{align*}
y_{i,t+1}^k(q) ~& =~ \theta_t^k(q) \ts y_{it}^k(q) -  \sum_{j \in \Jcal} a_{ij} \ts u_{j,t+1}^k(q)
\\
& \stackrel{(a)}=~ 
\theta_t^k(q) \ts \ts \mathbb P \{ D^k \geq t,~D^{k-1}=q\} \ts c_i 
\\
& \qquad - 
\sum_{\ell=1}^{k-1} \sum_{s \in \Tcal} \sum_{p \in \Tcal} \sum_{j \in \Jcal} a_{ij} \ts \theta_t^k(q) \ts  \ts \mathbb P \{ D^k \geq t ,~D^{k-1} = q \ts | \ts D^\ell \geq s,~D^{\ell-1} = p \} \ts u_{js}^\ell(p)
\\
& \qquad -
\sum_{s=1}^t \sum_{j \in \Jcal} a_{ij} \ts \theta_t^k(q) \ts\ts \mathbb P \{ D^k \geq t \ts | \ts D^k \geq s,~D^{k-1} = q\} \ts  u_{js}^k(q) - \sum_{j \in \Jcal} a_{ij} \ts u_{j,t+1}^k(q)
\\
& \stackrel{(b)}=~
\mathbb P \{ D^k \geq t+1,~ D^{k-1} = q\} \ts c_i 
\\
& \qquad - 
\sum_{\ell=1}^{k-1} \sum_{s \in \Tcal} \sum_{p \in \Tcal} \sum_{j \in \Jcal} a_{ij} \ts \mathbb P \{ D^k \geq t+1 ,~D^{k-1} = q \ts | \ts D^\ell \geq s,~D^{\ell-1} = p \} \ts u_{js}^\ell(p)
\\
& \qquad -
\sum_{s=1}^t \sum_{j \in \Jcal} a_{ij} \ts  \mathbb P \{ D^k \geq t+1 \ts | \ts D^k \geq s,~D^{k-1} = q\} \ts\ts u_{js}^k(q) - \sum_{j \in \Jcal} a_{ij} \ts u_{j,t+1}^k(q)
\\
& \stackrel{(c)}=~
\mathbb P \{ D^k \geq t+1,~ D^{k-1} = q\} \ts c_i 
\\
& \qquad - 
\sum_{\ell=1}^{k-1} \sum_{s \in \Tcal} \sum_{p \in \Tcal} \sum_{j \in \Jcal} a_{ij} \ts \mathbb P \{ D^k \geq t+1 ,~D^{k-1} = q \ts | \ts D^\ell \geq s,~D^{\ell-1} = p \} \ts u_{js}^\ell(p)
\\
& \qquad -
\sum_{s=1}^{t+1} \sum_{j \in \Jcal} a_{ij} \ts  \mathbb P \{ D^k \geq t+1 \ts | \ts D^k \geq s,~D^{k-1} = q\} \ts \ts u_{js}^k(q),
\end{align*}
where $(a)$ is by the induction assumption, $(b)$ uses the three identities given earlier in the proof and $(c)$ holds by noting that $\mathbb P \{ D^k \geq t+1 \ts | \ts D^k \geq t+1,~D^{k-1} = q\} = 1$ and collecting the terms.

The chain of equalities above shows that if $t \neq T$, then the result holds at the subsequent time period. We can use a similar argument to show that the result holds when $t=T$ as well. \qed

Using Lemmas \ref{lem:lp_prob} and \ref{lem:lp_capacity}, we give a proof for Theorem \ref{thm:bound}.

{\it \underline{Proof of Theorem \ref{thm:bound}}:}
\\
\indent Any feasible solution $(\yvec,\uvec,\wvec)$ to problem (\ref{eqn:lp_long}) satisfies $y_{it}^k(q) \geq 0$, in which case, dividing both sides of the equality in Lemma \ref{lem:lp_capacity} by $\mathbb P \{ D^k \geq t,~D^{k-1}=q\}$, we obtain 
\begin{align*}
& \sum_{\ell=1}^{k-1} \sum_{s \in \Tcal} \sum_{p \in \Tcal} \sum_{j \in \Jcal} a_{ij} \ts \frac{\mathbb P \{ D^k \geq t ,~D^{k-1} = q \ts | \ts D^\ell \geq s,~D^{\ell-1} = p \}}{\mathbb P \{ D^k \geq t,~D^{k-1}=q\} } \ts u_{js}^\ell(p)
\\
& \qquad \qquad \qquad \qquad \qquad \qquad -
\sum_{s=1}^t \sum_{j \in \Jcal} a_{ij} \ts \frac{\mathbb P \{ D^k \geq t \ts | \ts D^k \geq s,~D^{k-1} = q\}}{\mathbb P \{ D^k \geq t,~D^{k-1}=q\} } \ts u_{js}^k(q)
~\leq~ c_i.
\end{align*}
By the Bayes rule, the two fractions on the left side of the inequality above are, respectively, given by $\frac{\mathbb P \{D^\ell \geq s,~D^{\ell-1} = p \ts | \ts  D^k \geq t ,~D^{k-1} = q\}}{\mathbb P \{ D^\ell \geq s,~D^{\ell-1} = p \} } $ and $\frac{\mathbb P \{ \ts D^k \geq s,~D^{k-1} = q \ts| \ts D^k \geq t, ~D^{k-1} = q \}}{\mathbb P \{ D^k \geq s,~D^{k-1}=q\} } $, but for $s \leq t$, the last probability is equal to  $\frac{1}{\mathbb P \{ D^k \geq s,~D^{k-1}=q\} } $.
In this case, any feasible solution to the linear program in (\ref{eqn:lp_long}) is also a feasible solution to the linear program
\begin{align}
\max_{(\uvec,\wvec) \in \mathbb R_+^{|\Jcal| T^2K} \times \Wcal} ~~&  \sum_{k \in \Kcal} \sum_{t \in \Tcal} \sum_{q \in \Tcal} \sum_{j \in \Jcal} f_j \ts u_{jt}^k(q) 
\label{eqn:lp_scaled}
\\
\mbox{st}~~&\sum_{\ell =1}^{k-1} \sum_{s \in \Tcal} \sum_{p \in \Tcal} \sum_{j \in \Jcal} a_{ij} \ts \frac{\mathbb P \{ D^\ell \geq s, \ts D^{\ell-1} =p \ts|\ts D^k \geq t, \ts D^{k-1} = q \}}{\mathbb P\{ D^\ell \geq s ,~ D^{\ell-1} = p\}} \ts u_{js}^\ell(p) 
\nonumber
\\
& \qquad \quad   + \sum_{s=1}^t \sum_{j \in \Jcal} a_{ij} \ts \frac{1}{\mathbb P\{ D^k \geq s,~D^{k-1} = q\}} \ts u_{js}^k(q) ~\leq~ c_i \qquad \forall \ts i \in \Lcal,~t,q \in \Tcal,~k \in \Kcal
\nonumber
\\
& u_{jt}^k(q) \leq \lambda_{jt}^k \ts w_t^k(q) \qquad \forall \ts j \in \Jcal,~t,q \in \Tcal,~k \in \Kcal. \phantom{\Bigg\}}
\nonumber
\end{align}
Thus, the optimal objective value of problem (\ref{eqn:lp_scaled}) is an upper bound on that of problem (\ref{eqn:lp_long}), which is, in turn, an upper bound on the optimal total expected revenue.

By Lemma \ref{lem:lp_prob}, for any $\wvec \in \Wcal$, we have $w_t^k(q) = \mathbb P \{ D^k \geq t,~D^{k-1} = q\}$. In this case, making the change of variables $x_{jt}^k(q) = \frac{1}{\mathbb P\{ D^k \geq t,~D^{k-1} = q\}} \ts u_{jt}^k(q)$, problem (\ref{eqn:lp_scaled}) is equivalent to problem (\ref{eqn:lp}). \qed

\section{Fluid Approximation Through Linear Value Function Approximations}
\label{sec:lin_app}

The proof for Theorem \ref{thm:bound} that we give in Extended Results \ref{sec:bound} is based on relaxing the capacity constraints in the dynamic program in (\ref{eqn:dp}) through Lagrangian relaxation. In this section, we give an alternative proof for Theorem \ref{thm:bound} using linear approximations to the value functions in the dynamic program.~Using this approach, we also compare our fluid approximation with a linear programming approximation used by \cite{Ji23}. Because the remaining capacity of resource $i$  takes values in $\{0,\ldots,c_i\}$, we use \mbox{$\Ccal = \prod_{i \in \Lcal} \{ 0,\ldots, c_i\}$} to denote the set of all possible remaining capacity vectors. We can formulate the dynamic program in (\ref{eqn:dp}) as a linear program. In this linear program, we have one decision variable for each possible value of the state variable at each time period in each stage, whereas we have one constraint for each possible action that we can take at each possible value of the state variable at each time period in each state. Using the decision variables $\Jvec = (J_t^k(\yvec,q) : \yvec \in \Ccal,~t,q \in \Tcal,~k \in \Kcal)$, as well as recalling that we use $\opt$ to denote the optimal total expected revenue, the linear program corresponding to (\ref{eqn:dp}) is given by 
\begin{align}
& \opt = \!\!\! 
\min_{\Jvec \in \mathbb R^{T^2 K |\Ccal|}} ~~  J_1^1(\cvec, \Dhat^0)
\label{eqn:dp_lp}
\\
& ~~~~~~~~~~~\ts \mbox{st}~~  J_t^k(\yvec,q) ~\geq~ \sum_{j \in \Jcal} \lambda_{jt}^k \ts \bigg\{ f_j \ts u_j + \theta_t^k(q) \ts J_{t+1}^k(\yvec -\avec_j \ts u_j ,q) + ( 1 - \theta_t^k(q)) \ts J_1^{k+1}(\yvec - \avec_j   \ts u_j,t) \bigg\} 
\nonumber
\\
& \qquad \qquad \qquad \qquad \qquad \qquad \qquad \qquad \qquad \qquad \qquad \qquad \qquad \qquad 
\forall \ts \yvec \in \Ccal,~\uvec \in \Fcal(\yvec),~t,q \in \Tcal,~k \in \Kcal,
\nonumber 
\end{align}
where we follow the convention that $J_1^{K+1}(\yvec,q)$ is fixed at zero for all $\yvec \in \Ccal$ and $q \in \Tcal$. We use a linear value function approximation of the form \mbox{$\Jtilde_t^k(\yvec,q) = \beta_t^k(q) + \sum_{i \in \Lcal} \alpha_{it}^k(q) \ts y_i$}. To choose the slope and intercept parameters $\alphavec = (\alpha_{it}^k(q) : i \in \Lcal,~t,q \in \Tcal,~k \in \Kcal)$ and $\betavec = (\beta_t^k(q) : t,q \in \Tcal,~k \in \Kcal)$, we plug the linear value function approximation into (\ref{eqn:dp_lp}) to obtain the linear program
\begin{align}
& \Zbar_\linear =\!\!\!\!\!\!\! \min_{(\alphavec,\betavec) \in \mathbb R^{T^2 K (|\Lcal|+ 1 )}}  \beta_1^1(\Dhat^0) + \sum_{i \in \Lcal} \alpha_{i1}^1(\Dhat^0) \ts c_i
\label{eqn:alp}
\\
& ~~\ts \qquad \qquad \mbox{st}~~ \beta_t^k(q) + \sum_{i \in \Lcal} \alpha_{it}^k(q) \ts y_i 
\nonumber
\\
& ~~~~ \qquad \qquad \qquad  \geq~ \sum_{j \in \Jcal} \lambda_{jt}^k \ts u_j \ts \bigg\{ f_j - \theta_t^k(q) \sum_{i \in \Lcal} \alpha_{i,t+1}^k(q) \ts a_{ij} - (1-\theta_t^k(q)) \ts \sum_{i \in \Lcal} \alpha_{i1}^{k+1}(t) \ts a_{ij} \bigg\}  
\nonumber
\\
& ~~~\ts \qquad \qquad \qquad \qquad + \theta_t^k(q) \bigg\{ \beta_{t+1}^k(q) +  \sum_{i \in \Lcal} \alpha_{i,t+1}^k(q) \ts y_i \bigg\} + (1-\theta_t^k(q)) \bigg\{ \beta_1^{k+1}(t) +  \sum_{i \in \Lcal} \alpha_{i1}^{k+1}(t) \ts y_i \bigg\}
\nonumber
\\
& ~~ \qquad \qquad \qquad \qquad \qquad \qquad \qquad \qquad \qquad \qquad \qquad \qquad \qquad \quad 
\forall \ts \yvec \in \Ccal,~\uvec \in \Fcal(\yvec),~t,q \in \Tcal,~k \in \Kcal,
\nonumber 
\end{align}
where we follow the convention that $\beta_1^{K+1}(q)$ and $\alpha_{i1}^{K+1}(q)$ are fixed at zero for all $q \in \Tcal$ and $i \in \Lcal$. In the constraint above, we use the fact that $\Jtilde_t^k(\yvec,q)$ is linear in $\yvec$, in which case, we have $\mbox{$\Jtilde_{t+1}^k(\yvec - \avec_j \ts u_j , q)$} = \beta_{t+1}^k(q) + \sum_{i \in \Lcal} \alpha_{i,t+1}^k(q) \ts y_i - \sum_{i \in \Lcal} \alpha_{i,t+1}^k(q) \ts a_{ij} \ts u_j$ and $\Jtilde_1^{k+1}(\yvec - \avec_j \ts u_j,t) = \beta_1^{k+1}(t) + \sum_{i \in \Lcal} \alpha_{i1}^{k+1}(t) \ts y_i - \sum_{i \in \Lcal} \alpha_{i1}^{k+1}(t) \ts a_{ij} \ts u_j$, as well as the fact that $\sum_{j \in \Jcal} \lambda_{jt}^k = 1$. We obtain the linear program in (\ref{eqn:alp}) by constraining the value functions in the linear program in (\ref{eqn:dp_lp}) to be linear in the remaining capacities. Therefore, noting that we minimize the objective function in the linear program in (\ref{eqn:dp_lp}), the optimal objective value of problem (\ref{eqn:alp}) cannot be smaller than the optimal objective value of problem (\ref{eqn:dp_lp}), so we have $\Zbar_\linear \geq \opt$. The number of constraints in (\ref{eqn:alp}) is exponential in the numbers of resources and products. To obtain our fluid approximation, we use a relaxation of this linear program with the number of constraints polynomial in the numbers of resources and products. In particular, arranging the terms in the constraint in (\ref{eqn:alp}), we write the linear program in (\ref{eqn:alp}) equivalently as 
\begin{align}
& \Zbar_\linear =\!\!\!\!\!\!\! \min_{(\alphavec,\betavec) \in \mathbb R^{T^2 K (|\Lcal|+ 1 )}}   \beta_1^1(\Dhat^0) + \sum_{i \in \Lcal} \alpha_{i1}^1(\Dhat^0) \ts c_i
\label{eqn:alp_eq}
\\
& ~~\ts \qquad \qquad \mbox{st}~~ \beta_t^k(q) - \theta_t^k(q) \ts \beta_{t+1}^k(q) - (1-\theta_t^k(q)) \ts \beta_1^{k+1}(t) \phantom{\Bigg\}}
\nonumber
\\
& ~~~~ \qquad \qquad \qquad  \geq\!\!\! \max_{\yvec \in \Ccal, \uvec \in \Fcal(\yvec)} \Bigg\{ \sum_{j \in \Jcal} \lambda_{jt}^k \ts u_j \ts \bigg\{ f_j - \theta_t^k(q) \sum_{i \in \Lcal} \alpha_{i,t+1}^k(q) \ts a_{ij} - (1-\theta_t^k(q)) \ts \sum_{i \in \Lcal} \alpha_{i1}^{k+1}(t) \ts a_{ij} \bigg\}  
\nonumber
\\
& \qquad \qquad \qquad \qquad \qquad + \sum_{i \in \Lcal} y_i \ts \bigg\{ \theta_t^k(q)  \ts  \alpha_{i,t+1}^k(q) +  (1-\theta_t^k(q)) \ts  \alpha_{i1}^{k+1}(t) - \alpha_{it}^k(q) \bigg\} \Bigg\}
~~~~
\forall \ts t,q \in \Tcal,~k \in \Kcal,
\nonumber 
\end{align}
We relax the sets $\Ccal$ and $\Fcal(\yvec)$, in which case, the optimal objective value of the maximization problem in the constraint becomes larger, so the constraint becomes tighter.

Letting $\Ccaltilde = \prod_{i \in \Lcal} [0,c_i]$, we have $\Ccaltilde \supseteq \Ccal$. Noting that \mbox{$\Fcal(\yvec) = \{ \uvec \in \{0,1\}^{|\Jcal|} \!:\! a_{ij} \ts u_j \leq y_i ~\forall \ts i \in \Lcal,~j \in \Jcal\}$}, multiplying the constraint corresponding to product $j$ with $\lambda_{jt}^k$ and adding over all products, we define \mbox{$\Fcaltilde_t^k(\yvec) = \{ \uvec \in \{0,1\}^{|\Jcal|} : \sum_{j \in \Jcal} \lambda_{jt}^k \ts a_{ij} \ts u_j \leq y_i ~~\forall \ts i \in \Lcal\}$}, in which case,  $\Fcaltilde_t^k(\yvec) \supseteq \Fcal(\yvec)$. Thus, we can upper bound on the optimal objective value of (\ref{eqn:alp_eq}) by using the linear program
\begin{align}
& \Ztilde_\linear =\!\!\!\!\!\!\! \min_{(\alphavec,\betavec) \in \mathbb R^{T^2 K (|\Lcal|+ 1 )}}   \beta_1^1(\Dhat^0) + \sum_{i \in \Lcal} \alpha_{i1}^1(\Dhat^0) \ts c_i
\label{eqn:alp_relaxed}
\\
& ~~\ts \qquad \qquad \mbox{st}~~ \beta_t^k(q) - \theta_t^k(q) \ts \beta_{t+1}^k(q) - (1-\theta_t^k(q)) \ts \beta_1^{k+1}(t) \phantom{\Bigg\}}
\nonumber
\\
& ~~~~~~~ \qquad \qquad   \geq\!\!\!  \max_{\yvec \in \Ccaltilde, \uvec \in \Fcaltilde_t^k(\yvec)} \Bigg\{ \sum_{j \in \Jcal} \lambda_{jt}^k \ts u_j \ts \bigg\{ f_j - \theta_t^k(q) \sum_{i \in \Lcal} \alpha_{i,t+1}^k(q) \ts a_{ij} - (1-\theta_t^k(q)) \ts \sum_{i \in \Lcal} \alpha_{i1}^{k+1}(t) \ts a_{ij} \bigg\}  
\nonumber
\\
& \qquad \qquad \qquad \qquad \quad + \sum_{i \in \Lcal} y_i \ts \bigg\{ \theta_t^k(q)  \ts  \alpha_{i,t+1}^k(q) +  (1-\theta_t^k(q)) \ts  \alpha_{i1}^{k+1}(t) - \alpha_{it}^k(q) \bigg\} \Bigg\}
~~~~
\forall \ts t,q \in \Tcal,~k \in \Kcal.
\nonumber 
\end{align}
Because $\Ccaltilde \supseteq \Ccal$ and $\Fcaltilde_t^k(\yvec) \supseteq \Fcal(\yvec)$, the constraint in (\ref{eqn:alp_relaxed}) is at least as tight as the one in (\ref{eqn:alp_eq}), so $\Ztilde_\linear \geq \Zbar_\linear$. We can show that there exists an optimal solution to problem (\ref{eqn:alp_relaxed}), where the decision variables $(\alpha_{it}^k(q) : i \in \Lcal,~t,q \in \Tcal,~k \in \Kcal)$ satisfy \mbox{$\alpha_{it}^k(q) \geq \theta_t^k(q) \ts \alpha_{i,t+1}^k(q) + (1 - \theta_t^k(q)) \ts \alpha_{i1}^{k+1}(t)$} for all $i \in \Lcal$, $t,q \in \Tcal$ and $k \in \Kcal$. Theorem 2 in \cite{Ad05} shows this result when there is a single stage with a fixed number of customer arrivals, but we can use the same proof technique in our problem setting. In particular, letting if $(\alphavecbar,\betavecbar)$ is an optimal solution to (\ref{eqn:alp_relaxed}) with  \mbox{$\alphabar_{it}^k(q) < \theta_t^k(q) \ts \alphabar_{i,t+1}^k(q) + (1 - \theta_t^k(q)) \ts \alphabar_{i1}^{k+1}(t)$} for some $i \in \Lcal$, $t,q \in \Tcal$ and $k \in \Kcal$, then we can increase $\alphabar_{it}^k(q)$ by $\epsilon$ and decrease $\betabar_t^k(q)$ by $c_i \ts \epsilon$ for small enough $\epsilon > 0$ to obtain a feasible solution to (\ref{eqn:alp_relaxed}) that provides an objective value that is at least as small as the optimal objective value. Thus, we can add the constraint \mbox{$\alpha_{it}^k(q) \geq \theta_t^k(q) \ts \alpha_{i,t+1}^k(q) + (1 - \theta_t^k(q)) \ts \alpha_{i1}^{k+1}(t)$} for all $i \in \Lcal$, $t,q \in \Tcal$ and $k \in \Kcal$ to problem (\ref{eqn:alp_relaxed}) without changing its optimal objective value, in which case, the coefficient of $y_i$ in the maximization problem in the constraint is non-positive. Thus, we need to choose the value of $y_i$ as small as possible. By the definition of $\Fcaltilde_t^k(\yvec)$, the smallest value of the decision variable $y_i$ is $\sum_{j \in \Jcal} \lambda_{jt}^k \ts a_{ij} \ts u_j$. Thus, the  maximization problem in the constraint in (\ref{eqn:alp_relaxed}) is equivalent to
\begin{align*}
& \max_{\uvec \in \{0,1\}^{|\Jcal|}} \Bigg\{ \sum_{j \in \Jcal} \lambda_{jt}^k \ts u_j \ts \bigg\{ f_j - \theta_t^k(q) \sum_{i \in \Lcal} \alpha_{i,t+1}^k(q) \ts a_{ij} - (1-\theta_t^k(q)) \ts \sum_{i \in \Lcal} \alpha_{i1}^{k+1}(t) \ts a_{ij} \bigg\}  
\nonumber
\\
&  \qquad \qquad \qquad \qquad \qquad + \sum_{i \in \Lcal} \sum_{j \in \Jcal} \lambda_{jt}^k \ts a_{ij} \ts u_j \ts \bigg\{ \theta_t^k(q)  \ts  \alpha_{i,t+1}^k(q) +  (1-\theta_t^k(q)) \ts  \alpha_{i1}^{k+1}(t) - \alpha_{it}^k(q) \bigg\} \Bigg\}
\\
& \qquad = \max_{\uvec \in \{0,1\}^{|\Jcal|}} \Bigg\{ \sum_{j \in \Jcal} \lambda_{jt}^k \ts u_j \ts \bigg\{ f_j - \sum_{i \in \Lcal} a_{ij} \ts \alpha_{it}^k(q) \bigg\}  \Bigg\}
~=~
\sum_{j \in \Jcal} \lambda_{jt}^k \ts \Big[ f_j - \sum_{i \in \Lcal} a_{ij} \ts \alpha_{it}^k(q) \Big]^+,
\end{align*}
where the first equality follows by arranging the terms and the second equality holds by noting the set of feasible solutions in the maximization problem on the left side of the equality.

We replace the optimal objective value of the maximization problem in the constraint in (\ref{eqn:alp_relaxed}) with $\sum_{j \in \Jcal} \lambda_{jt}^k \ts [ f_j - \sum_{i \in \Lcal} a_{ij} \ts \alpha_{it}^k(q) ]^+$. To linearize the constraint, we  use the decision variable $\eta_{jt}^k(q)$ to capture the expression  $[f_j - \sum_{i \in \Lcal} a_{ij} \ts \alpha_{it}^k(q)]^+$. Recalling that we can add the constraint \mbox{$\alpha_{it}^k(q) \geq \theta_t^k(q) \ts \alpha_{i,t+1}^k(q) + (1 - \theta_t^k(q)) \ts \alpha_{i1}^{k+1}(t)$} for all $i \in \Lcal$, $t,q \in \Tcal$ and $k \in \Kcal$, (\ref{eqn:alp_relaxed}) is equivalent to 
\begin{align}
& \Ztilde_\linear ~~=\!\!\!\!\!\!\! \min_{(\alphavec,\betavec,\etavec) \in \mathbb R^{T^2 K (|\Lcal| + 1)} \times \mathbb R_+^{T^2K |\Jcal|}}  \beta_1^1(\Dhat^0) + \sum_{i \in \Lcal} \alpha_{i1}^1(\Dhat^0) \ts c_i
\label{eqn:alp_final}
\\
& ~~~~\ts \qquad \qquad \qquad  \mbox{st}~~ 
\alpha_{it}^k(q) ~\geq~ \theta_t^k(q) \ts \alpha_{i,t+1}^k(q) + (1 - \theta_t^k(q)) \ts \alpha_{i1}^{k+1}(t)
\qquad \forall \ts i \in \Lcal,~t,q \in \Tcal,~k \in \Kcal \phantom{\sum_{i \in \Lcal}}
\nonumber
\\
& ~~\ts \qquad \qquad \qquad \qquad 
\beta_t^k(q) - \theta_t^k(q) \ts \beta_{t+1}^k(q) - (1-\theta_t^k(q)) \ts \beta_1^{k+1}(t) ~\geq~ \sum_{j \in \Jcal} \lambda_{jt}^k \ts \eta_{jt}^k(q) \quad
\forall \ts t,q \in \Tcal,~k \in \Kcal \phantom{\sum_{i \in \Lcal}}
\nonumber 
\\
& ~~\ts \qquad \qquad \qquad \qquad \eta_{jt}^k(q) ~\geq~ f_j - \sum_{i \in \Lcal} a_{ij} \ts \alpha_{it}^k(q) \qquad \forall \ts j \in \Jcal,~t,q \in \Tcal,~k \in \Kcal,
\nonumber
\end{align}
where we use the decision variables $\etavec = (\eta_{jt}^k(q) : j \in \Jcal,~t,q \in \Tcal,~k \in \Kcal)$. The objective function in (\ref{eqn:alp_final}) matches that in (\ref{eqn:tighten_bound}). Viewing $\mu_{it}^k(q)$ as the slack variable for the first constraint in problem (\ref{eqn:tighten_bound}) and expressing this constraint as an inequality constraint by dropping the slack variable, the first constraint in  (\ref{eqn:alp_final}) matches that in (\ref{eqn:tighten_bound}). We can show that there exists an optimal solution to problem (\ref{eqn:alp_final}) such that the second constraint is satisfied as equality. In particular, if $(\alphavecbar,\betavecbar,\etavecbar)$ is an optimal solution to (\ref{eqn:alp_final}) with $\betabar_t^k(q) - \theta_t^k(q) \ts \betabar_{t+1}^k(q) - (1-\theta_t^k(q)) \ts \betabar_1^{k+1}(t) >\sum_{j \in \Jcal} \lambda_{jt}^k \ts \etabar_{jt}^k(q)$ for some $t,q \in \Tcal$ and $k \in \Kcal$, then we can decrease $\betabar_t^k(q)$ by $\epsilon$ for some small enough $\epsilon > 0$ to obtain a feasible solution to (\ref{eqn:alp_final}) that provides an objective value that is at least as small as the optimal objective value. Thus, the second constraint in (\ref{eqn:alp_final}) matches that in (\ref{eqn:tighten_bound}). Lastly, the third constraint in (\ref{eqn:alp_final}) matches that in (\ref{eqn:tighten_bound}). Therefore, problems (\ref{eqn:tighten_bound}) and (\ref{eqn:alp_final}) are equivalent to each other. Because we obtain our fluid approximation through the dual of problem (\ref{eqn:tighten_bound}), the upper bound from our fluid approximation is $\Ztilde_\linear$ and we can obtain our fluid approximation by using linear value function approximations.~It turns out the linear program in \cite{Ji23} is based on a further relaxation of the maximization problem in the constraint in (\ref{eqn:alp_relaxed}). We have $\Fcaltilde_t^k(\cvec) \supseteq \Fcaltilde_t^k(\yvec)$. Also, we have 
\begin{align}
& \max_{\yvec \in \Ccaltilde,~\uvec \in \Fcaltilde_t^k(\cvec)} \Bigg\{ \sum_{j \in \Jcal} \lambda_{jt}^k \ts u_j \ts \bigg\{ f_j - \theta_t^k(q) \sum_{i \in \Lcal} \alpha_{i,t+1}^k(q) \ts a_{ij} - (1-\theta_t^k(q)) \ts \sum_{i \in \Lcal} \alpha_{i1}^{k+1}(t) \ts a_{ij} \bigg\}  
\nonumber
\\
& \qquad \qquad \qquad \qquad \qquad \qquad \qquad + \sum_{i \in \Lcal} y_i \ts \bigg\{ \theta_t^k(q)  \ts  \alpha_{i,t+1}^k(q) +  (1-\theta_t^k(q)) \ts  \alpha_{i1}^{k+1}(t) - \alpha_{it}^k(q) \bigg\} \Bigg\}
\nonumber
\\
&  \qquad = \sum_{j \in \Jcal} \lambda_{jt}^k \ts \Big[ f_j - \theta_t^k(q) \sum_{i \in \Lcal} \alpha_{i,t+1}^k(q) \ts a_{ij} - (1-\theta_t^k(q)) \ts \sum_{i \in \Lcal} \alpha_{i1}^{k+1}(t) \ts a_{ij} \Big]^+  
\nonumber
\\
&  \qquad \qquad \qquad \qquad \qquad \qquad \qquad + \sum_{i \in \Lcal} c_i \Big[ \theta_t^k(q)  \ts  \alpha_{i,t+1}^k(q) +  (1-\theta_t^k(q)) \ts  \alpha_{i1}^{k+1}(t) - \alpha_{it}^k(q) \Big]^+\!\!\!,
\label{eqn:full_relaxed}
\end{align}
where the equality holds by $\Fcaltilde_t^k(\cvec) = \{0,1\}^{|\Jcal|}$, as well as by the fact that if the coefficient of $y_i$ in the maximization problem above is non-negative, then it is optimal to set this variable to $c_i$. 

Because the maximization problem in (\ref{eqn:full_relaxed}) is a relaxation of the maximization problem in the constraint in (\ref{eqn:alp_relaxed}), replacing the maximization problem in the constraint in (\ref{eqn:alp_relaxed}) with the~optimal objective value of the maximization problem in (\ref{eqn:full_relaxed}), we can obtain an upper bound on the optimal~objective value of the linear program in (\ref{eqn:alp_relaxed}) through the linear program 
\begin{align}
& \Zhat_\linear =\!\!\!\!\!\!\! \min_{(\alphavec,\betavec) \in \mathbb R^{T^2 K (|\Lcal| + 1)}}  \beta_1^1(\Dhat^0) + \sum_{i \in \Lcal} \alpha_{i1}^1(\Dhat^0) \ts c_i
\label{eqn:alp_full_relaxed}
\\
& ~~\ts \qquad \qquad \mbox{st}~~ \beta_t^k(q) - \theta_t^k(q) \ts \beta_{t+1}^k(q) - (1-\theta_t^k(q)) \ts \beta_1^{k+1}(t) \phantom{\sum_{j \in \Jcal} \Big]^+}
\nonumber
\\
& ~~~~~~~ \qquad \qquad   ~\geq~  \sum_{j \in \Jcal} \lambda_{jt}^k \ts \Big[ f_j - \theta_t^k(q) \sum_{i \in \Lcal} \alpha_{i,t+1}^k(q) \ts a_{ij} - (1-\theta_t^k(q)) \ts \sum_{i \in \Lcal} \alpha_{i1}^{k+1}(t) \ts a_{ij} \Big]^+
\nonumber
\\
& \qquad \qquad \qquad \qquad \quad + \sum_{i \in \Lcal} c_i \Big[ \theta_t^k(q)  \ts  \alpha_{i,t+1}^k(q) +  (1-\theta_t^k(q)) \ts  \alpha_{i1}^{k+1}(t) - \alpha_{it}^k(q) \Big]^+
~~~~
\forall \ts t,q \in \Tcal,~k \in \Kcal.
\nonumber 
\end{align}
Because the constraint in (\ref{eqn:alp_full_relaxed}) is at least as tight as the  constraint in (\ref{eqn:alp_relaxed}), the optimal objective values of these two linear programs satisfy $\Zhat_\linear \geq \Ztilde_\linear$. Using the fact that $\Ztilde_\linear \geq \Zbar_\linear \geq \opt$, we get $\Zhat_\linear \geq \Ztilde_\linear \geq \Zbar_\linear \geq \opt$. Thus, noting that the optimal objective value of problem (\ref{eqn:alp_final}) is $\Ztilde_\linear$, the linear programs in (\ref{eqn:alp_final}) and (\ref{eqn:alp_full_relaxed}) both provide upper bounds on the optimal total expected revenue, but the upper bound provided by (\ref{eqn:alp_final}) is at least as tight as the one provided by (\ref{eqn:alp_full_relaxed}). The linear program in (\ref{eqn:alp_full_relaxed}) is precisely the analogue of the linear program given in Section 3 of \cite{Ji23} when the demands occur in multiple stages, there is a random number of customer arrivals in each stage and the numbers of customer arrivals in successive stages are dependent on each other. The author uses this linear program to compute an upper bound on the optimal total expected revenue, but does not construct an approximate policy by using an optimal solution to this linear program. By the preceding discussion, the upper bound provided by our fluid approximation is at least as tight as the one provided by the linear program in \cite{Ji23}.

\vspace{-1mm}

\section{Verifying the Fluid Approximation Through Decisions of the Optimal Policy}
\label{sec:path_bound}

\vspace{-2mm}

We use the decisions made by the optimal policy to construct a feasible solution to the problem (\ref{eqn:lp}) in such a way that the objective value provided by this solution for problem (\ref{eqn:lp}) matches the optimal total expected revenue. %In this case, it follows that the optimal objective value of problem (\ref{eqn:lp}) is an upper bound on the optimal total expected revenue. 
We use the Bernoulli random variable $X_{jt}^k$ to capture the decision of the optimal policy for product $j$ at time period $t$ in stage $k$, where $X_{jt}^k = 1$ if and only if the optimal policy accepts a request for product $j$ at time period $t$ in stage $k$. By the dynamic program in (\ref{eqn:dp}), given that we reach time period $t$ in stage~$k$ before this stage is over, the random variable $X_{jt}^k$ depends on the remaining capacities of the resources at time period $t$ in stage $k$ and the demand in stage $k-1$. Given that we reach time period $t$ in stage $k$, the remaining capacities of the resources at time period $t$ in stage $k$ depend on the product requests at all time periods up to time period~$t$ in stage $k$, as well as the demand random variables in stages up to stage~$k$.~Thus, given that we reach time period $t$ in stage $k$ before this stage is over, the random variable $X_{jt}^k$ depends on the demand random variables $D^1,\ldots, D^{k-1}$, but not on the other demand random variables.~Given that we do not reach time period $t$ in stage~$k$, we have $X_{jt}^k=0$. By the preceding discussion, the random variable $X_{jt}^k$ depends on the demand random variables $D^1,\ldots,D^{k-1}$, as well as \mbox{${\bf 1}(D^k \geq t)$}, but not on the other demand random variables.  Letting \mbox{$\xbar_{jt}^k(q) = \mathbb E\{ X_{jt}^k \ts | \ts D^k \geq t,~D^{k-1} = q\}$}, we show that the solution \mbox{$\xvecbar = (\xbar_{jt}^k(q) : j \in \Jcal,~t,q \in \Tcal,~k \in \Kcal)$} provides an objective value for problem (\ref{eqn:lp}) that is equal to the optimal total expected revenue. 
The optimal total expected revenue is $\opt = \sum_{k \in \Kcal} \sum_{t \in \Tcal} \sum_{j \in \Jcal} f_j \ts \mathbb E\{ X_{jt}^k \}$.~Given $D^k < t$ so that we do not reach time period~$t$ in stage $k$, we have $X_{jt}^k = 0$, so by the tower property of conditional expectations, we obtain \mbox{$\mathbb E\{ X_{jt}^k\} = \sum_{q \in \Tcal} \mathbb E\{ X_{jt}^k \ts | \ts D^k \geq t, \ts D^{k-1} = q\} \ts \mathbb P \{ D^k \geq t, \ts D^{k-1} = q\}$}. By the definition of $\xbar_{jt}^k(q)$, the last equality is equivalent to $\mathbb E\{ X_{jt}^k\} = \sum_{q \in \Tcal} \xbar_{jt}^k(q) \ts \mathbb P \{ D^k \geq t,~D^{k-1} = q\}$. Therefore, we can express the optimal total expected revenue as $\opt = \sum_{k \in \Kcal} \sum_{t \in \Tcal} \sum_{j \in \Jcal} f_j \ts \sum_{q \in \Tcal} \xbar_{jt}^k(q) \ts \mathbb P\{ D^k \geq t, \ts D^{k-1} = q\}$, in which case, noting the objective function of problem (\ref{eqn:lp}),  the objective value provided by the solution $\xvecbar$ for problem (\ref{eqn:lp}) is equal to the optimal total expected revenue. 

We show that the solution $\xvecbar$ is feasible to problem (\ref{eqn:lp}). We use the Bernoulli random variable $\Lambda_{jt}^k$ to capture whether there is a request for product $j$ at time period $t$ in stage $k$, where $\Lambda_{jt}^k = 1$ if and only if there is a request for product $j$ at time period $t$ in stage $k$. The Bernoulli random variable $\Lambda_{jt}^k$ has parameter $\lambda_{jt}^k$ and it is independent of the demand random variables. If $D^k < t$, then we have $X_{jt}^k = 0$, which implies that $X_{jt}^k \leq {\bf 1}(D^k \geq t)$. On the other hand, we can accept a request for a product only if there is request for it, which implies that $X_{jt}^k \leq \Lambda_{jt}^k$. In this case, we obtain $X_{jt}^k \leq {\bf 1}(D^k \geq t,~\Lambda_{jt}^k=1)$. If we take the expectation of both sides of the last inequality conditional on $D^k \geq t$ and $D^{k-1} = q$, then we obtain $\xbar_{jt}^k(q) \leq \mathbb E\{ \Lambda_{jt}^k \ts | \ts D^k \geq t,~D^{k-1} = q\} = \lambda_{jt}^k$, where the equality holds because the random variable $\Lambda_{jt}^k$ is independent of $D^{k-1}$ and $D^k$.  Thus, the solution $\xvecbar$ satisfies the second constraint in (\ref{eqn:lp}). We give a useful identity. Given that $D^\ell < s$, so that we do not reach time period $s$ in stage $\ell$, we have \mbox{$X_{js}^\ell = 0$}.~Therefore, considering some stage $\ell=1,\ldots,k-1$, by the tower property of conditional expectations, we get
\begin{align}
& \mathbb E\{ X_{js}^\ell  \ts | \ts D^k \geq t, \ts D^{k-1} = q\}
\nonumber
\\
& \qquad \quad =~ \sum_{p \in \Tcal} \mathbb E\{ X_{js}^\ell  \ts | \ts D^\ell \geq s,\ts D^{\ell-1} = p, \ts D^k \geq t, \ts D^{k-1} = q\} \ts\ts \mathbb P \{  D^\ell \geq s,\ts D^{\ell-1} = p \ts | \ts D^k \geq t, \ts D^{k-1} = q\}
\nonumber
\\
& \qquad \quad \stackrel{(a)}=~ \sum_{p \in \Tcal} \mathbb E\{ X_{js}^\ell  \ts | \ts D^\ell \geq s,\ts D^{\ell-1} = p\} \ts\ts \mathbb P \{  D^\ell \geq s,\ts D^{\ell-1} = p \ts | \ts D^k \geq t, \ts D^{k-1} = q\},
\nonumber
\\
& \qquad \quad =~ \sum_{p \in \Tcal} \xbar_{js}^\ell(p) \ts \mathbb P \{  D^\ell \geq s,\ts D^{\ell-1} = p \ts | \ts D^k \geq t, \ts D^{k-1} = q\},
\label{eqn:cond_accept}
\end{align}
where $(a)$ holds because $X_{js}^\ell$ depends on the demand random variables $D^1,\ldots,D^{\ell-1}$, as well as ${\bf 1}(D^\ell \geq s)$, but given that $D^{\ell -1}=p$, $D^1,\ldots,D^{\ell-1}$ are independent of $D^{k-1}$ and $D^k$.

The total capacity consumption of resource $i$ up to time period $t$ in stage $k$ cannot exceed the capacity of the resource, so we have $\sum_{\ell = 1}^{k-1} \sum_{s \in \Tcal} \sum_{j \in \Jcal} a_{ij} \ts X_{js}^\ell + \sum_{s=1}^t \sum_{j \in \Jcal} a_{ij} \ts X_{js}^k \leq c_i$ with probability one.~Taking the expectation of both sides of the last inequality conditional on $D^k \geq t$ and $D^{k-1} = q$, we obtain the chain of inequalities 
\begin{align*}
c_i~& \geq~
\sum_{\ell = 1}^{k-1} \sum_{s \in \Tcal} \sum_{j \in \Jcal} a_{ij} \ts \mathbb E\{ X_{js}^\ell  \ts | \ts D^k \geq t, \ts D^{k-1} = q\} + \sum_{s=1}^t \sum_{j \in \Jcal} a_{ij} \ts \mathbb E\{ X_{js}^k \ts | \ts D^k \geq t, \ts D^{k-1} = q\}
\\
& \stackrel{(b)}=~
\sum_{\ell = 1}^{k-1} \sum_{s \in \Tcal} \sum_{j \in \Jcal} a_{ij} \sum_{p \in \Tcal} \xbar_{js}^\ell(p) \ts\ts  \mathbb P \{  D^\ell \geq s,\ts D^{\ell-1} = p \ts | \ts D^k \geq t, \ts D^{k-1} = q\} + \sum_{s=1}^t \sum_{j \in \Jcal} a_{ij} \ts \xbar_{js}^k(q),
\end{align*}
where $(b)$ uses (\ref{eqn:cond_accept}) and the fact that $\mathbb E\{ X_{js}^k \ts | \ts D^k \geq t, \ts D^{k-1} = q\} = \mathbb E\{ X_{js}^k \ts | \ts D^k \geq s, \ts D^{k-1} = q\} = \xbar_{js}^k(q)$ because $X_{js}^k$ depends on  ${\bf 1}(D^k \geq s), D^{k-1},\ldots,D^1$. Thus, the solution $\xvecbar$ satisfies the first constraint in (\ref{eqn:lp}) as well. By the discussion in this section, the solution $\xvecbar$ is feasible to problem (\ref{eqn:lp}) and provides an objective value for this problem that is equal to the optimal total expected revenue.~In this case, the optimal objective value of problem (\ref{eqn:lp}) is at least as large as the optimal total expected revenue. This discussion gives an alternative proof of the fact that the optimal objective value of our fluid approximation is an upper bound on the optimal total expected revenue.~This discussion is not constructive in the sense that it does not allow us to derive the form of problem~(\ref{eqn:lp}).~Instead, it verifies that the optimal objective value of problem (\ref{eqn:lp}) is an upper bound on the optimal total expected revenue once we are given the form of problem (\ref{eqn:lp}). On the other hand, using the Lagrangian relaxation strategy as in Extended Results~\ref{sec:bound} or using the linear value function approximations as in Extended Results \ref{sec:lin_app} allows us to derive the form of problem (\ref{eqn:lp}), while also establishing that problem (\ref{eqn:lp}) provides an upper bound on the optimal total expected revenue.

\section{Comparing the Fluid Approximation and the Offline Bound}
\label{sec:offline_bound}

We give problem instances to demonstrate that the upper bound on the optimal total expected revenue provided by our fluid approximation can be tighter or looser than the offline bound.~First, we give a problem instance for which the upper bound provided by our fluid approximation is tighter than the offline bound. There are two stages. The support of the demand in each stage is one. There is one resource and two products. Thus, we have $K = 1$, $T = 1$, $\Lcal = \{1\}$ and $\Jcal = \{1,2\}$. The capacity of the resource is $c_1 = 1$. The revenues of the products are $f_1 = 1$ and $f_2 = 2$. Fixing the demand in the stage right before the beginning of the selling horizon at $\Dhat^0$, the distribution of the demand random variables in each stage is characterized by the conditional probabilities $\mathbb P \{ D^1 =0 \ts | \ts D^0 = \Dhat^0\} = \mathbb P \{ D^1 = 1 \ts | \ts D^0 = \Dhat^0 \} = \frac 12$, $\mathbb P \{ D^2 = 1 \ts | \ts D^1 = 0 \} = 1$  and $\mathbb P \{ D^2 = 0 \ts | \ts D^1  = 1 \} = \mathbb P \{ D^2 = 1 \ts | \ts D^1 = 1 \} = \frac 12$. Lastly, the probability of getting a request for each product at each time period in each stage is given by $\lambda_{11}^1 = 1$, $\lambda_{11}^2 = \frac 12$ and $\lambda_{21}^2 = \frac 12$. Considering the offline bound, because the capacity of the resource is one and the revenue of the second product is larger than that of the first product, if there is a request for the second product at any time period in any stage, then we accept the request for the second product. If there is a request for the first product but not the second product, then we accept the request for the first product. The probability that we have a request for the second product in the first stage is zero. To have a request for the second product in the second stage, we must have a demand of one in the second stage and this demand must be for the second product. Thus, the probability that we have a request for the second product is $(\frac 12 + \frac 12 \times \frac 12) \times \frac 12 = \frac 38$. On the other hand, the probability that we have a request for the first product but not the second product is given by $ \frac 12 \times \frac 12 + \frac 12 \times \frac 12 + \frac 12 \times \frac 12 \times \frac 12 = \frac 58$, where the three terms in the sum, respectively, correspond to having one demand in the first stage and no demand in the second stage, no demand in the first stage and one demand in the second stage and one demand in both stages, along with a demand for the first product and no demand for the second product. Thus, the offline bound is $\Zbar_\text{\sf offline} = \frac 38 \times 2 + \frac 58 \times 1 = \frac{11}{8}$. Considering problem (\ref{eqn:lp}), given that the demand in the previous stage is $q$, if we never have a request for product $j$ at time period $t$ in stage $k$, then we can drop the decision variable $x_{jt}^k(q)$. For our problem instance, problem (\ref{eqn:lp}) is given by
\vspace{1mm}
\begin{align*}
\Zbar_\lp ~=~ \max ~~ & \frac 12 \ts x_{11}^1(\Dhat^0) + \frac 12 \ts x_{11}^2(0) +  \frac 12 \ts 2 \ts x_{21}^2(0) + \frac 14 \ts x_{11}^2(1) +   \frac 14  \ts 2 \ts x_{21}^2(1)
\\
\mbox{st} ~~& x_{11}^1(\Dhat^0) + x_{11}^2(1) + x_{21}^2(1) \leq 1
\\
& x_{11}^2(0) + x_{21}^2(0) \leq 1
\\
& x_{11}^1(\Dhat^0) \in [0,1],~ x_{11}^2(0), x_{21}^2(0), x_{11}^2(1), x_{21}^2(1) \in [0 , \textstyle \frac 12].
\end{align*}
The optimal objective value of the problem above is $\frac 54$. Therefore, we have $\Zbar_\lp = \frac 54 \leq \frac{11}{8} = \Zbar_\text{\sf offline}$, so the upper bound from our fluid approximation can be tighter than the offline bound.

\vspace{1mm}

Second, we give a problem instance for which the upper bound provided by our fluid approximation is looser than the offline bound.  This problem instance has the same parameters as the problem instance in the previous paragraph other than the distribution of the demand random variables in different stages and the request probabilities for different products at different time periods in different stages. In particular, we have $\mathbb P \{ D^1 = 1 \ts | \ts D^0 = \Dhat^0 \} = 1$, $\mathbb P\{ D^2 = 1 \ts | \ts D^1 = 1 \} = 1$, so the demand in each stage is always one. Furthermore, we have $\lambda_{11}^1 = \lambda_{21}^1 = \frac 12$ and $\lambda_{21}^2 = \frac 12$, so we have a request for each of the two products at the first time period in the first stage with equal probabilities, whereas we have a request for the second product at the first time period in the second stage with probability $\frac 12$ and we do not have a request at the first time period in the second stage with probability $\frac 12$. To ensure that we have a product request at each time period in each stage with probability one, we can introduce a third product with a revenue of zero such that we have a request for the product at the first time period in the second stage with probability $\frac 12$. To keep our problem instance succinct, we do not explicitly introduce such a product. We have  a request for the second product at some period in some stage with probability $\frac 12 \times \frac 12 + \frac 12 \times \frac 12 + \frac 12 \times \frac 12 = \frac 34$, where the three terms in the sum, respectively, correspond to having a request for the second product in the first stage but not in the second stage, not in the first stage but in the second stage and in both stages. On the other hand, the probability that we have a request for the first product but not the second product is $\frac 12 \times \frac 12 = \frac 14$. Thus, the offline bound is $\Zbar_\text{\sf offline} = \frac 34 \times 2 + \frac 14 = \frac 74$. Considering problem (\ref{eqn:lp}), for our problem instance, this problem is given by 
\vspace{1mm}
\begin{align*}
\Zbar_\lp ~=~ \max ~~ & x_{11}^1(\Dhat^0) +  2 \ts x_{21}^1(\Dhat^0) + 2 \ts x_{21}^2(1) 
\\
\mbox{st} ~~&  x_{11}^1(\Dhat^0) +  x_{21}^1(\Dhat^0) + x_{21}^2(1)  \leq 1
\\
& x_{11}^1(\Dhat^0),  x_{21}^1(\Dhat^0), x_{21}^2(1) \in [0 , \textstyle \frac 12].
\end{align*}
The optimal objective value of the problem above is $2$. Therefore, we have $\Zbar_\lp = 2 \geq \frac 74 = \Zbar_\text{\sf offline}$, so the upper bound from our fluid approximation can be looser than the offline bound.

\section{Detailed Description of the Approximate Policy}
\label{sec:policy_steps}

In Table \ref{tab:policy_steps}, we give a detailed description of our approximate policy. We use the variable $y_{it}^k$ to keep track of the remaining capacity of resource $i$ at the beginning of time period $t$ in stage $k$.~In Line~1, we initialize the remaining capacities of the resources, as well as set the demand in the stage right before the beginning of the selling horizon to $\Dhat^0$, which is a part of the problem data. In Line 2, the for loop enumerates over the stages in the selling horizon. We use the variable $k$ to keep track of the current stage. In Line 3, we use the variable $t$ to keep track of the current time period and $S^k$ to keep track of whether the current stage should end. We initialize the current time period to one and ensure that the current stage should not end at least until after the first time period. In Line 4, the while loop enumerates over the time periods in the current stage.

In Line 5, we sample the product requested at the current time period. We use the variable $J_t^k$ to keep track of the product requested at the current time period. In Line 6, we sample whether the approximate policy would accept the product request at the current time period. The binary variable $A_t^k$ takes value one if and only if the approximate policy would accept the product request.~In~Line 7, we check whether we have enough capacity to serve the product requested at the current time period. The binary variable $C_t^k$ takes value one if and only if we have enough capacity. In Lines 8 to 10, we check whether the approximate policy would accept the  product requested at the current time period and we have enough capacity to serve the requested product.~If so, then we decrease the capacities of the resources used by the requested product.

In Line 11, recalling that $\theta_t^k(q) = \mathbb P \{ D^k \geq t+1 \ts | D^k \geq t,~D^k = q\}$, we sample whether there is one more time period in the current stage. The binary variable $E^k$ takes value one if and only if there is one more time period in the current stage. In Lines 12 to 16, we check whether there is one more time period. If so, then we increment the current time period. Otherwise, we set the variable $S^k$ to indicate that the current stage is to end and set the value of the demand in the current stage, as well as the remaining capacities of the resources at the beginning of the next stage.  

\begin{table}[t]
\hrule
\vspace{1mm}
{\tt \scriptsize\phantom{~}[1]~~~~}
Set $y_{i1}^1 = c_i$ for all $i \in \Lcal$ and $D^0 = \Dhat^0$. 
\\
{\tt \scriptsize\phantom{~}[2]~~~~}
For $k = 1,\ldots,K$, do the following steps.
\\
{\tt \scriptsize\phantom{~}[3]~~~~} \hspace{5mm}
Set $t=1$, $S^k = \text{False}$.
\\
{\tt \scriptsize\phantom{~}[4]~~~~} \hspace{5mm}
While $S^k = \text{False}$ do the following steps.
\\
{\tt \scriptsize\phantom{~}[5]~~~~} \hspace{10mm}
Sample $J_t^k \in \Jcal$ such that $\mathbb P \{ J_t^k = j \} = \lambda_{jt}^k$.
\\
{\tt \scriptsize\phantom{~}[6]~~~~} \hspace{10mm}
Sample $A_t^k \in \{0,1\}$ such that $\mathbb P \{ A_t^k = 1 \} = \gamma \ts \frac{\xbar_{J_t^k,t}^k(D^{k-1})}{\lambda_{J_t^k,t}^k}$.
\\
{\tt \scriptsize\phantom{~}[7]~~~~} \hspace{10mm} 
Set $C_t^k = {\bf 1}(y_{it}^k \geq a_{i,J_t^k}~\forall \ts i \in \Lcal)$.
\\
{\tt \scriptsize\phantom{~}[8]~~~~} \hspace{10mm}
If $A_t^k = 1$ and $C_t^k = 1$, then 
\\
{\tt \scriptsize\phantom{~}[9]~~~~} \hspace{15mm}
Set $y_{i,t+1}^k = y_{it}^k - a_{i,J_t^k}$ for all $i \in \Lcal$
\\
{\tt \scriptsize[10]~~~~} \hspace{10mm}
End if.
\\
{\tt \scriptsize[11]~~~~} \hspace{10mm}
Sample $E^k \in \{0,1\}$ such that $\mathbb P\{ E^k = 1\} = \theta_t^k(D^{k-1})$.
\\
{\tt \scriptsize[12]~~~~} \hspace{10mm}
If $E^k = 1$, then 
\\
{\tt \scriptsize[13]~~~~} \hspace{15mm} 
Increment $t$ by one
\\ 
{\tt \scriptsize[14]~~~~} \hspace{10mm} 
Else
\\
{\tt \scriptsize[15]~~~~} \hspace{15mm}
Set $S^k = \text{True}$, $D^k = t$ and $y_{i1}^{k+1} = y_{i,t+1}^k$ for all $i \in \Lcal$
\\
{\tt \scriptsize[16]~~~~} \hspace{10mm}
End if.
\\
{\tt \scriptsize[17]~~~~} \hspace{5mm} 
End While.
\\
{\tt \scriptsize[18]~~~~}
End For.
\\
~ \vspace{-5mm}
\hrule
\caption{Detailed description of the approximate policy from the fluid approximation.}
\label{tab:policy_steps}
\vspace{-6mm}
\end{table}

\vspace{-1mm}

\section{Upper Bound on the Optimality Gap of the Approximate Policy}
\label{sec:asymp_ub}

\vspace{-1mm}

We give a problem instance such that if we scale the number of stages and resource capacities with the same rate $\theta$, then the ratio between the optimal total expected revenue and the optimal objective value of (\ref{eqn:lp}) is at most $O\Big(1 - \frac{1}{\sqrt \theta} \Big)$. There are $K$ stages with $K \geq 72 \pi$. The demand in each stage is always one. There is one resource and one product. The capacity of the resource is $\frac 12 \ts K$.~The revenue of the product is one. At each time period, we have a demand for the product  with probability~$\frac 12$. Because there is one product, it is optimal to accept the product requests as much as the capacity of the resource allows. Thus, letting $B^1,\ldots,B^K$ be independent Bernoulli random variables with parameter $\frac 12$, the optimal total expected revenue is $\opt = \mathbb E\{ \min\{ \sum_{k=1}^K B^k , \frac 12 \ts K \} \} = \mathbb E\{ \min\{ \sum_{k=1}^K(B^k - \frac 12 ) , 0\}\} + \frac 12 K$. The function $h(x) = \min\{x,0\}$ is Lipschitz continuous with Lipschitz constant one. We have $\sum_{k=1}^K\text{Var}(B^k - \frac 12 ) = \frac 14 \ts K$ and $\mathbb E\{ |B^k - \frac 12|^3 \} = \frac 18$. In this case, letting $N$ be the standard normal random variable, by Theorem 3.1 in \cite{ChGoSh10}, we obtain 
\begin{align*}
\frac{\mathbb E\{ \min\{ \sum_{k=1}^K(B^k - \frac 12 ) , 0\}\}}{\sqrt {\frac 14 K}}
~\leq~
\mathbb E\{ \min\{N,0\} \} + \frac{3}{\sqrt K}.
\end{align*}
Arranging the terms in the inequality above and using the identity $\mathbb E\{ \min\{N,0\}\} = - \frac{1}{\sqrt{2 \pi}}$, we have $\mathbb E\{  \min\{ \sum_{k=1}^K(B^k - \frac 12 ) , 0\}\} \leq - \frac{1}{2} \sqrt{\frac{K}{2\pi}} + \frac{3}{2}$. Because \mbox{$\opt = \mathbb E\{  \min\{ \sum_{k=1}^K(B^k - \frac 12 ) , 0\}\} + \frac 12 K$}, we get $\opt \leq \frac 12 \ts K - \frac{1}{2} \sqrt{\frac{K}{2\pi}} + \frac{3}{2}$. Noting that the demand in each stage is always one and there is one product, using the decision variables $\xvec = (x_{11}^k(1) : k \in \Kcal)$, problem (\ref{eqn:lp}) for our problem instance is given by $\Zbar_\lp = \max_{\xvec \in \mathbb R_+^K} \{\sum_{k \in \Kcal} x_{11}^k(1) : \sum_{k \in \Kcal} x_{11}^k(1) \leq \frac 12 K,~ x_{11}^k(1) \leq \frac 12 ~\forall \ts k \in \Kcal\}$. The optimal objective value of this problem is $\Zbar_\lp = \frac 12 \ts K$. Therefore, we obtain $\frac{\opt}{\Zbar_\lp} \leq 1 - \frac{1}{\sqrt{2\pi K}} + \frac{3}{K} = 1 - \frac{1}{\sqrt{2\pi K}} + \frac{3}{\sqrt K \sqrt K} \leq 1 - \frac{1}{2\sqrt{2\pi K}} \leq 1 - \frac{1}{\sqrt{26 K}}$, where the second inequality uses the fact that $K \geq 72 \ts \pi$. In our problem instance, the number of stages is $K$ and the capacity of the resource is $\frac 12 K$. Thus, if we scale the number of stages and capacity of the resource with rate $\theta$, then the optimal total expected revenue and the optimal objective value of  (\ref{eqn:lp}) is separated by  $1 - \frac{1}{\sqrt{26 \ts \theta}}$.

\section{Problem Instance with Increasing Support for Demands}
\label{sec:large_capacity}

When there is a single stage, so the demand model is not aware of the calendar and the dependence between the demands in different stages is not an issue, we can construct fluid approximations such that the relative gap between the optimal objective value of the fluid approximation and the optimal total expected revenue vanishes as the capacities of the resources get large, irrespective of how the demand is scaled; see, for example, \cite{BaHo23}. We give a counterexample to demonstrate that if we scale the support of the demand along with the capacities of the resources, then the relative gap between the optimal objective value of problem (\ref{eqn:lp}) and the optimal total expected revenue does not necessarily vanish, even when the demands at different stages are independent of each other.~We consider a problem instance with $K=3$ stages.~There is one resource with a capacity of $C+1$ for an even integer $C$. There are two products indexed by $\{1,2\}$.~The revenue associated with the two products are $f_1 = 1$ and $f_2 = C/4$. The demand in the first stage can take two values with $\mathbb P \{ D^1 = 0 \} = 1/2$ and $\mathbb P \{ D^1 = C \} = 1/2$. The demand in the second and third stages have the distributions $\mathbb P \{ D^2 = 1 \} = 1$ and $\mathbb P \{ D^3 = C/2\} = 1$. The probabilities of getting requests for the different products are given by $\lambda_{1t}^1 = 1$ for all $t=1,\ldots,C$, $\lambda_{11}^2=1$, $\lambda_{1t}^3 = 1$ for all $t=1,\ldots,C/2-1$ and $\lambda_{2, C/2}^3=1$. All other request probabilities are zero. Thus, we have a request for the second product only at the last time period in the last stage. 

The optimal policy for this problem instance is to accept all product requests until there is one unit of remaining capacity and save that unit of capacity for the request for the second product.~In this case, if the demand in the first stage is $C$, then we obtain a total revenue of \mbox{$C + \frac14 C = \frac 54 C$}. If the demand in the first stage is zero, then we obtain a total revenue of $\frac 12 C+ \frac 14 C = \frac 34 C$. Therefore, the optimal total expected revenue is $C$, so $\opt = C$. In the linear program in (\ref{eqn:lp}), noting the last constraint, if $\lambda_{jt}^k = 0$, then we can drop the decision variables $\{x_{jt}^k(q) : q \in \Tcal\}$. In this case, using $\Dhat^0$ to denote the demand right before the beginning of the selling horizon, for the problem instance in the previous paragraph, the linear program in (\ref{eqn:lp}) is given by
\begin{align*}
\Zbar _\lp ~=~\max~~~ & \frac 12 \sum_{t=1}^C x_{1t}^1(\Dhat^0) + \frac 12 x_{11}^2(0) + \frac12 x_{11}^2(C) +\!\!\! \sum_{t=1}^{C/2-1} x_{1t}^3(1) + \frac 14 C \ts x_{2,C/2}^3(1)
\\
\mbox{st} ~~~ & \sum_{t=1}^C x_{1t}^1(\Dhat^0) \leq C+1
\\
& x_{11}^2(0) \leq C+1 \phantom{\sum^c}
\\
& \sum_{t=1}^C x_{1t}^1(\Dhat^0) + x_{11}^2(C) \leq C+1
\\
& \frac 12 \sum_{t=1}^C x_{1t}^1(\Dhat^0) + \frac 12 \ts x_{11}^2(0) + \frac 12 \ts x_{11}^2(C) + \!\!\! \sum_{t=1}^{C/2-1} x_{1t}^3(1) + x_{2, C/2}^3(1) \leq C+1,
\end{align*}
where we understand that all decision variables are restricted to be in $[0,1]$. Setting all decision variables to one yields a feasible solution, so $\Zbar_\lp = \frac 12 C + 1 + \frac 12 C -1 + \frac 14 C = \frac 54 C$.

Thus, for this problem instance, we have $\frac{\opt}{\Zbar_\lp} = \frac{4}{5}$. Therefore, the ratio $\frac{\opt}{\Zbar_\lp}$ stays away from one for this problem instance as the resource capacity gets large.

\section{Performance Guarantee for the Approximate Policy}
\label{sec:constant_factor}

We show that the total expected revenue of the approximate policy satisfies $\frac{\apx}{\opt} \geq \frac{\apx}{\Zbar_\lp} \geq \frac{1}{4L}$.~We use the random variables $\Psi_t^k(q)$, $G_{jt}^k$ and $N_{it}^k(q)$ as defined at the beginning of Section~\ref{sec:perf_asymp}.~By the discussion right after (\ref{eqn:apx_rev}), if we can show that $\mathbb P \{ G_{jt}^k = 1 \ts | \ts \Psi_t^k(q) = 1\}  \geq \alpha$, then  \mbox{$\apx \geq \gamma \ts \alpha \ts \Zbar_\lp$}.~By (\ref{eqn:union}), to lower bound the probability $\mathbb P \{ G_{jt}^k = 1 \ts | \ts \Psi_t^k(q) = 1\}$, it is enough to upper bound the probability $ \mathbb P \{\sum_{\ell=1}^{k-1} \sum_{s \in \Tcal} \sum_{p \in \Tcal} \Psi_s^\ell(p) \ts N_{is}^\ell(p) + \sum_{s=1}^t N_{is}^k(q) \geq c_i \ts | \ts D^{k-1} = q  \}$.~Recall that the random variables $\{ N_{it}^k(q) : i \in \Lcal,~t \in \Tcal,~k \in \Kcal\}$ are independent of demands, so we get
\begin{align*}
& \mathbb E \bigg\{ \sum_{\ell=1}^{k-1} \sum_{s \in \Tcal} \sum_{p \in \Tcal} \Psi_s^\ell(p) \ts N_{is}^\ell(p) + \sum_{s=1}^t N_{is}^k(q) \ts \Big| \ts D^{k-1} = q \bigg\}
\nonumber
\\
~& \qquad \quad \stackrel{(a)}=\ts
\sum_{\ell=1}^{k-1} \sum_{s \in \Tcal} \sum_{p \in \Tcal} \mathbb P \{ D^\ell \geq s,~ D^{\ell-1} = p \ts | \ts D^{k-1} = q \}  \ts \sum_{j \in \Jcal} a_{ij} \ts \gamma \ts  \xbar_{js}^\ell(p)
+
\sum_{s=1}^t \sum_{j \in \Jcal} a_{ij} \ts \gamma \ts \xbar_{js}^k(q) 
~\stackrel{(b)}\leq~ 
\gamma \ts c_i,
\end{align*}
where $(a)$ holds by using the definition of $\Psi_t^k(p)$ and  $\mathbb E\{ N_{it}^k(q) \} = \sum_{j \in \Jcal} a_{ij} \ts \xbar_{jt}^k(q)$, whereas $(b)$ holds by noting that $\xvecbar$ satisfies the first constraint in problem (\ref{eqn:lp}), as well as using the fact that conditional on $D^{k-1}$, $D^k$ is independent of $D^1,\ldots,D^{k-1}$, in which case, for $\ell \leq k-1$, we have $\mathbb P \{ D^\ell \geq s,~ D^{\ell-1} = p \ts | \ts D^{k-1} = q \} = \mathbb P \{ D^\ell \geq s,~ D^{\ell-1} = p \ts | \ts D^k \geq t,~ D^{k-1} = q \}$. Therefore, we get
\begin{align*}
& \mathbb P \bigg\{\sum_{\ell=1}^{k-1} \sum_{s \in \Tcal} \sum_{p \in \Tcal} \Psi_s^\ell(p) \ts N_{is}^\ell(p) + \sum_{s=1}^t N_{is}^k(q) \ts \geq \ts c_i \ts | \ts D^{k-1} = q  \bigg\}
\\
& \qquad \qquad \qquad \qquad \qquad \qquad \stackrel{(c)}\leq
\frac{1}{c_i} \ts 
\mathbb E \bigg\{\sum_{\ell=1}^{k-1} \sum_{s \in \Tcal} \sum_{p \in \Tcal} \Psi_s^\ell(p) \ts N_{is}^\ell(p) + \sum_{s=1}^t N_{is}^k(q) \geq c_i \ts | \ts D^{k-1} = q  \bigg\} 
\leq
\gamma,
\end{align*}
where $(c)$ is the Markov inequality. Thus, by (\ref{eqn:union}), we have $\mathbb P \{ G_{jt}^k = 1 \ts | \ts D^{k-1} = q\} \geq 1 - L \ts \gamma$, yielding $\apx \geq \gamma \ts (1 - L \ts \gamma) \ts \Zbar_\lp$. Setting $\gamma = \frac{1}{2L}$ and using Theorem \ref{thm:bound}, we get $\frac{1}{4L} \leq \frac{\apx}{\Zbar_\lp} \leq \frac{\apx}{\opt}$. \qed

\section{Auxiliary Results for Concentration Inequalities}
\label{sec:auxiliary}

We give proofs for two results used in Section \ref{sec:perf_asymp}. First, we show that $|M_i^k(\ell) - M_i^k(\ell+1)| \leq \frac{2}{\epsilon^3}$ with probability one.  Letting $V^k(q) = \mathbb P \{ D^k = q \}$ and $\theta = 1 - \epsilon$ for notational brevity, we define 
\begin{align}
Q^k(p,q) = \frac{1}{\theta}  \big[ \mathbb P \{ D^k = q \ts | \ts D^{k-1} = p \} - (1-\theta) \ts V^k(q) \big].
\label{eqn:aux_trans}
\end{align}
By the assumption that $\mathbb P \{ D^k =q \ts | \ts D^{k-1} = p \} \geq \epsilon$, we have $Q^k(p,q) \geq \frac{1}{\theta} \ts (\epsilon - (1-\theta)) = 0$.~Also, we have $\sum_{q \in \Tcal} Q^k(p,q) = \frac{1}{\theta} \big[ \sum_{q \in \Tcal} \mathbb P \{ D^k = q \ts | \ts D^{k-1} = p\} - (1-\theta) \ts \sum_{q \in \Tcal} V^k(q) ] = 1$. Thus, we can use $Q^k(p,q)$ to characterize the transition probabilities of a non-stationary Markov chain. Consider the non-stationary Markov chain $Y^0,Y^1,Y^2,\ldots$ over the state space $\Tcal$ characterized by the transition probabilities $\mathbb P \{ Y^k = q \ts | \ts Y^{k-1} = p\} = Q^k(p,q)$ for all $p,q \in \Tcal$ with $\mathbb P\{ Y^0 = q \} = {\bf 1}(D^0 = q)$. To show that $|M_i^{k-1}(\ell) - M_i^{k-1}(\ell+1)| \leq \frac{2}{\epsilon^3}$ with probability one, we will use two preliminary lemmas. In the next lemma, we slightly extend Theorem 4.9 in \cite{LePe17}, which characterizes the mixing times of Markov chains, to non-stationary Markov chains.

\begin{lem}
\label{lem:cond_dev}
For all $p,q \in \Tcal$ and $k,\ell \in \Kcal$ with $\ell \geq k+1$, we have 
\begin{align*}
\mathbb P \{ D^\ell = q \ts | \ts D^k = p \} - \mathbb P \{ D^\ell = q\} ~=~ (1-\epsilon)^{(\ell-k)} \ts \Big[ \mathbb P\{ Y^\ell = q \ts | \ts Y^k = p\} - \mathbb P \{ D^\ell = q \} \Big].
\end{align*}
\end{lem}

\noindent{\it Proof:} Letting $\theta  =1-\epsilon$ for notational brevity, for all $p,q \in \Tcal$ and $k,\ell \in \Kcal$ with $\ell \geq k+1$, we claim that $\mbox{$\mathbb P \{ D^\ell = q \ts | \ts D^k = p\}$} = (1-\theta^{\ell-k}) \ts \mathbb P \{ D^\ell = q\} + \theta^{\ell-k} \ts \mathbb P \{ Y^\ell =q \ts | \ts Y^k = p\}$. We show the claim by using induction over $\ell = k+1,\ldots,K$. Consider the case $\ell = k+1$. By the definition of $Q^{k+1}(p,q)$, we have $\mbox{$\mathbb P \{ Y^{k+1} =q \ts | \ts Y^k = p\}$} = Q^{k+1}(p,q) = \mbox{$\frac{1}{\theta}  \big[ \mathbb P \{ D^{k+1} = q \ts | \ts D^k = p \}$} - \mbox{$(1-\theta) \ts \mbox{$\mathbb P \{ D^{k+1} = q\}$} \big]$}$, so arranging the terms, we get $\mbox{$(1-\theta)$} \ts \mathbb P \{ D^{k+1} = q\} + \theta \ts \mbox{$\mathbb P \{ Y^{k+1} =q \ts | \ts Y^k = p\}$} = \mbox{$\mathbb P \{ D^{k+1} = q \ts | \ts D^k = p\}$}$, establishing the claim for $\ell = k+1$. Assuming that the claim holds for $\ell \geq k +1$, we show that the claim holds for $\ell+1 \geq k+1$ as well. Arranging the terms in the definition of $Q^{\ell+1}(s,q)$ in (\ref{eqn:aux_trans}), we have $\mbox{$\mathbb P \{ D^{\ell+1} = q \ts | \ts D^\ell = s\}$} = \theta \ts Q^{\ell+1}(s,q) + (1- \theta) \ts V^{\ell+1}(q)$. In this case, noting the identity $\mbox{$\mathbb P \{ Y^{\ell+1} = q \ts | \ts Y^k = p \}$} = \sum_{s \in \Tcal} \mathbb P \{ Y^{\ell+1} = q \ts | \ts Y^\ell = s\} \ts \mathbb P \{ Y^\ell = s \ts | \ts Y^k = p \}$, it follows that we~have~the~chain~of equalities 
\begin{align}
& \sum_{s \in \Tcal} \mathbb P \{ D^{\ell+1} =q \ts | \ts D^\ell =s \} \ts\ts \mathbb P \{ Y^\ell =s \ts | \ts Y^k = p\}
\nonumber
\\
& \qquad \quad=~
\sum_{s \in \Tcal} 
\Big[\theta \ts Q^{\ell+1}(s,q) + (1- \theta) \ts V^{\ell+1}(q) \Big]
\ts \mathbb P \{ Y^\ell =s \ts | \ts Y^k = p\}
\nonumber
\\
& \qquad \quad\stackrel{(a)}=~
\theta \ts \sum_{s \in \Tcal} \mathbb P \{ Y^{\ell+1} = q \ts | \ts Y^\ell = s\} \ts \mathbb P \{ Y^\ell =s \ts | \ts Y^k = p\}
+
(1-\theta) \ts V^{\ell+1}(q)
\nonumber
\\
& \qquad \quad  \stackrel{(b)}=~
\theta \ts \mathbb P \{ Y^{\ell+1} = q \ts | \ts Y^k = p\} + (1-\theta) \ts \mathbb P \{ D^{\ell+1} = q\},
\label{eqn:cond_dev_inner}
\end{align}
where $(a)$ holds because  $Q^{\ell+1}(s,q) = \mathbb P \{ Y^{\ell+1} = q \ts | \ts Y^\ell = s\}$ and $\sum_{s \in \Tcal} \mathbb P \{ Y^\ell = s \ts | \ts Y^k = p\} = 1$, whereas $(b)$ uses the fact that we have $V^{\ell+1}(q) = \mathbb P \{ D^{\ell+1} = q\}$ by its definition.

We have $\mathbb P \{ D^{\ell} = s \ts | \ts D^k = p\} = (1-\theta^{\ell-k}) \ts \mathbb P \{ D^\ell = s\} + \theta^{\ell-k} \ts \mathbb P \{ Y^\ell =s \ts | \ts Y^k = p\} $ by the induction assumption, so noting that $\mathbb P \{ D^{\ell+1} = q\} = \sum_{s \in \Tcal} \mathbb P \{ D^{\ell+1} =q \ts | \ts D^\ell = s\} \ts \mathbb P \{ D^\ell = s\}$, we have 
\begin{align*}
\quad~~ & \mathbb P \{ D^{\ell+1} = q \ts | \ts D^k = p\} 
~=~
\sum_{s \in \Tcal} \mathbb P \{ D^{\ell+1} = q \ts | \ts D^\ell = s\} \ts\ts \mathbb P \{ D^{\ell} = s \ts | \ts D^k = p\}
\\
&\quad \qquad  =~
\sum_{s \in \Tcal} \mathbb P \{ D^{\ell+1} =q \ts | \ts D^\ell =s \} \ts \Big[ (1-\theta^{\ell-k}) \ts \mathbb P \{ D^\ell = s\} + \theta^{\ell-k} \ts \mathbb P \{ Y^\ell =s \ts | \ts Y^k = p\} \Big]
\\
& \quad \qquad \stackrel{(c)}=~
(1 - \theta^{\ell-k}) \ts \mathbb P \{ D^{\ell+1} = q \}
+\theta^{\ell-k} \ts 
\Big[\theta \ts \mathbb P \{ Y^{\ell+1} = q \ts | \ts Y^k = p\} + (1-\theta) \ts \mathbb P \{ D^{\ell+1} = q\}
\Big] \phantom{\sum_{s \in \Tcal}}
\\
& \quad \qquad =~
(1 - \theta^{\ell+1-k}) \ts \mathbb P \{ D^{\ell+1} = q \}
+\theta^{\ell+1-k} \ts 
\mathbb P \{ Y^{\ell+1} = q \ts | \ts Y^k = p\}, 
\end{align*}
where $(c)$ is by (\ref{eqn:cond_dev_inner}). Thus, the claim holds for $\ell + 1 \geq k$ as well. By the claim, $\mathbb P \{ D^\ell = q \ts | \ts D^k = p\} = \mbox{$(1-\theta^{\ell-k}) \ts \mathbb P \{ D^\ell = q\}$} + \theta^{\ell-k} \ts \mathbb P \{ Y^\ell =q \ts | \ts Y^k = p\}$, so the lemma follows by arranging the terms. \qed

In the next lemma, we build on Lemma \ref{lem:cond_dev} to bound the gap between the probabilities \mbox{$\mathbb P \{ D^k \geq s,~D^{k-1} = p \ts | \ts D^\ell = q\}$} and  $\mathbb P \{ D^k \geq s,~D^{k-1} = p\}$ for $\ell \geq k+1$.

\begin{lem}
\label{lem:cond_gap}
For all $p,q,s \in \Tcal$ and $k,\ell \in \Kcal$ with $\ell \geq k+1$, we have 
\begin{align*}
\sum_{p \in \Tcal} \Big| \mathbb P \{ D^k \geq s,~D^{k-1} = p \ts | \ts D^\ell = q\} - \mathbb P \{ D^k \geq s,~D^{k-1} = p\} \Big|
~\leq~
\frac{1}{\epsilon} \ts (1-\epsilon)^{\ell-k}.
\end{align*}
\end{lem}

\vspace{-0.25mm}

\noindent{\it Proof:} Using the fact that $\mathbb P \{ D^{k+1} = q \ts | \ts D^k = p\} \geq \epsilon$ for all $p,q \in \Tcal$ and $k \in \Kcal$, it is simple to check that $\mathbb P \{ D^k = q \} \geq \epsilon$ for all $q \in \Tcal$ and $k \in \Kcal$. Also, for $\ell \geq k+1$, by the Bayes rule, we have
\begin{align}
& \frac{\mathbb P \{ D^k = p,~ D^{k-1} = r\}}{\mathbb P \{ D^\ell = q\}} \ts  \mathbb P \{ D^\ell = q \ts | \ts D^k = p \} 
~ = ~
\frac{\mathbb P \{ D^k = p,~ D^{k-1} = r\}}{\mathbb P \{ D^\ell = q\}} \ts  \frac{\mathbb P \{ D^\ell = q ,~ D^k = p \}}{\mathbb P \{ D^k = p\}}
\nonumber
\\
~& \qquad \qquad \quad =~
\mathbb P \{ D^{k-1} = r \ts | \ts D^k = p\} \ts \mathbb P \{ D^k = p \ts | \ts D^\ell = q\}
~=~
\mathbb P \{ D^k = p,~D^{k-1} = r \ts | \ts D^\ell = q\},
\label{eqn:flip_order}
\end{align}
where the last equality uses the fact that given $D^k$, $D^{k-1}$ is independent of $D^\ell$. Noting that we have $\mathbb P \{ Y^\ell = q \ts | \ts Y^k = p\} - \mathbb P \{ D^\ell = q\} \in [-1,1]$, by Lemma \ref{lem:cond_dev}, we obtain the chain of inequalities $- (1-\epsilon)^{\ell-k} \leq \mathbb P \{ D^\ell = q \ts | \ts D^k = p \} - \mathbb P \{ D^\ell = q\} \leq (1-\epsilon)^{\ell-k}$. In this case, multiplying this chain of inequalities with $\frac{\mathbb P \{ D^k = p,~ D^{k-1} = r\}}{\mathbb P \{ D^\ell = q\}}$, as well as using (\ref{eqn:flip_order}), we have  
\begin{align*}
& -(1-\epsilon)^{\ell-k} \ts \frac{\mathbb P \{ D^k = p,~ D^{k-1} = r\}}{\mathbb P \{ D^\ell = q\}}
\\
& \qquad \qquad \qquad \qquad \leq~
\mathbb P \{ D^k = p,~D^{k-1} = r \ts | \ts D^\ell = q\}
-
\mathbb P \{ D^k = p,~D^{k-1} = r \}
\\
& \qquad \qquad \qquad \qquad \leq~
(1-\epsilon)^{\ell-k} \ts \frac{\mathbb P \{ D^k = p,~ D^{k-1} = r\}}{\mathbb P \{ D^\ell = q\}}.
\end{align*}
Because $\mathbb P \{ D^\ell = q\} \geq \epsilon$, the chain of inequalities above yields $\mbox{$-\frac{1}{\epsilon} \ts (1-\epsilon)^{\ell-k} \ts \mathbb P \{ D^k = p,~D^{k-1}=r\}$} \leq \mathbb P \{ D^k = p,~D^{k-1} = r \ts | \ts D^\ell = q\}
-
\mathbb P \{ D^k = p,~D^{k-1} = r \} \leq \frac{1}{\epsilon} \ts (1-\epsilon)^{\ell-k} \ts \mathbb P \{ D^k = p, ~ D^{k-1} = r \}$.

To conclude the proof, if we add the last chain of inequalities in the previous paragraph over all $p \geq s$, then we obtain the chain of inequalities  
\begin{align*}
& -\frac{1}{\epsilon} \ts (1-\epsilon)^{\ell-k} \ts \ts \mathbb P \{ D^k \geq s,~ D^{k-1} = r\}
\\
~& \qquad \qquad \qquad \qquad \leq~
\mathbb P \{ D^k \geq s ,~D^{k-1} = r \ts | \ts D^\ell = q\}
-
\mathbb P \{ D^k \geq s ,~D^{k-1} = r \}
\\
& \qquad \qquad \qquad \qquad  \leq~
\frac{1}{\epsilon} \ts (1-\epsilon)^{\ell-k} \ts \mathbb P \{ D^k \geq s ,~ D^{k-1} = r\}.
\end{align*}
Using the fact that $\sum_{r \in \Tcal} \mathbb P \{ D^k\geq s,~D^{k-1} = r\} = \mathbb P \{ D^k \geq s\} \leq 1$, the chain of inequalities above yields $\sum_{r \in \Tcal} | \mathbb P \{ D^k \geq s ,~D^{k-1} = r \ts | \ts D^\ell = q\}
-
\mathbb P \{ D^k \geq s ,~D^{k-1} = r \}| \leq \frac{1}{\epsilon} \ts (1-\epsilon)^{\ell-k}$. \qed

Using the lemma above, we can show the first result that we are interested in. In the next lemma, we show that $|M_i^k(\ell) - M_i^k(\ell+1)| \leq \frac{2}{\epsilon^3}$ with probability one.

\begin{lem}
\label{lem:martingale_diff}
For all $i \in \Lcal$ and $k,\ell \in \Kcal$ and $\ell+1 \leq k$, with probability one, we have 
\begin{align*}
|M_i^k(\ell) - M_i^k(\ell+1)| ~\leq~ \frac{2}{\epsilon^3}.
\end{align*}
\end{lem}

\noindent{\it Proof:} By the definition of $\Psi_t^k(q)$, we have $\mathbb E\{ \Psi_t^k(q) \ts | \ts \Dvec^{[v,\ell]} \} = \mathbb P \{ D^k \geq t,~D^{k-1} = q \ts | \ts \Dvec^{[v,\ell]}\}$ for any $v \in \Kcal$. Noting that $V_i^k = \sum_{v=1}^k \sum_{s \in \Tcal} \sum_{p \in \Tcal} \Psi_s^v(p) \ts n_{is}^v(p)$, we obtain 
\begin{align}
& | \ts M_i^k(\ell) - M_i^k(\ell+1) \ts | ~=~ |\ts \mathbb E\{ V_i^k \ts | \ts \Dvec^{[\ell,k]}\} - \mathbb E\{ V_i^k \ts | \ts \Dvec^{[\ell+1,k]}\} \ts | 
\nonumber
\\
& \qquad \leq 
\sum_{v=1}^k \sum_{s \in \Tcal} \sum_{p \in \Tcal}n_{is}^v(p) \ts \Big| \ts 
\mathbb P \{ D^v \geq s,~D^{v-1} = p \ts | \ts \Dvec^{[\ell,k]} \} 
-
\mathbb P \{ D^v \geq s,~D^{v-1} = p \ts | \ts \Dvec^{[\ell+1,k]} \} \ts 
\Big|
\nonumber
\\
& \qquad \stackrel{(a)}= 
\sum_{v=1}^{\ell+1} \sum_{s \in \Tcal} \sum_{p \in \Tcal}n_{is}^v(p) \ts \Big| \ts 
\mathbb P \{ D^v \geq s,~D^{v-1} = p \ts | \ts \Dvec^{[\ell,k]} \} 
-
\mathbb P \{ D^v \geq s,~D^{v-1} = p \ts | \ts \Dvec^{[\ell+1,k]} \} \ts 
\Big|,
\label{eqn:cond_gap_first}
\end{align}
where $(a)$ follows from the fact that if $v \geq \ell +2$, then both $D^v$ and $D^{v-1}$ are deterministic functions of $\Dvec^{[\ell+1,k]}$. Therefore, if $v \geq \ell +2$, then the probabilities $\mathbb P \{ D^v \geq s,~D^{v-1} = p \ts | \ts \Dvec^{[\ell,k]} \}$ and \mbox{$\mathbb P \{ D^v \geq s,~D^{v-1} = p \ts | \ts \Dvec^{[\ell+1,k]} \}$} take the same value of zero or one. For $v \leq \ell$, given $D^\ell$, both $D^v$ and $D^{v-1}$ are independent of $D^{\ell+1},\ldots,D^k$. Thus, if $v \leq \ell$, then $\mathbb P \{ D^v \geq s,~D^{v-1} = p \ts | \ts \Dvec^{[\ell,k]} \} = \mathbb P \{ D^v \geq s,~D^{v-1} = p \ts | \ts D^\ell \}$. Similarly, for $v \leq \ell$, given $D^{\ell+1}$, both $D^v$ and $D^{v-1}$ are independent of $D^{\ell+2},\ldots,D^k$. Thus, if $v\leq \ell$, then $\mathbb P \{ D^v \geq s,~D^{v-1} = p \ts | \ts \Dvec^{[\ell+1,k]} \} = \mathbb P \{ D^v \geq s,~D^{v-1} = p \ts | \ts D^{\ell+1} \}$. Lastly, by the definition of $n_{it}^k(q)$, we have $n_{it}^k(q) = \mathbb E\{ N_{it}^k(q) \} = \sum_{j \in \Jcal} a_{ij} \ts \gamma \ts \xbar_{jt}^k(q) \leq \gamma \sum_{j \in \Jcal} \lambda_{jt}^k \leq 1$, where the first inequality uses the fact that $\xvecbar$ satisfies the second constraint in problem (\ref{eqn:lp}). In this case, we can upper bound the expression on the right side of (\ref{eqn:cond_gap_first}) as 
\begin{align}
& \sum_{v=1}^{\ell+1} \sum_{s \in \Tcal} \sum_{p \in \Tcal}n_{is}^v(p) \ts \Big| \ts 
\mathbb P \{ D^v \geq s,~D^{v-1} = p \ts | \ts \Dvec^{[\ell,k]} \} 
-
\mathbb P \{ D^v \geq s,~D^{v-1} = p \ts | \ts \Dvec^{[\ell+1,k]} \} \ts 
\Big|
\nonumber
\\
& \qquad =~\sum_{v=1}^\ell \sum_{s \in \Tcal} \sum_{p \in \Tcal}n_{is}^v(p) \ts \Big| \ts 
\mathbb P \{ D^v \geq s,~D^{v-1} = p \ts | \ts D^\ell \} 
-
\mathbb P \{ D^v \geq s,~D^{v-1} = p \ts | \ts D^{\ell+1} \} \ts 
\Big|
\nonumber
\\
& \qquad \qquad +
\sum_{s \in \Tcal} \sum_{p \in \Tcal}n_{is}^{\ell+1} (p) \ts \Big| \ts 
\mathbb P \{ D^{\ell+1} \geq s,~D^\ell = p \ts | \ts \Dvec^{[\ell,k]} \} 
-
\mathbb P \{ D^{\ell+1} \geq s,~D^\ell = p \ts | \ts \Dvec^{[\ell+1,k]} \} \ts 
\Big|
\nonumber
\\
& \qquad \stackrel{(b)}\leq~
\sum_{v=1}^\ell \sum_{s \in \Tcal} \sum_{p \in \Tcal} \ts \Big| \ts 
\mathbb P \{ D^v \geq s,~D^{v-1} = p \ts | \ts D^\ell \} 
-
\mathbb P \{ D^v \geq s,~D^{v-1} = p \} \ts 
\Big|
\nonumber
\\
&\qquad \qquad + \sum_{v=1}^\ell \sum_{s \in \Tcal} \sum_{p \in \Tcal} \Big| \ts 
\mathbb P \{ D^v \geq s,~D^{v-1} = p \ts | \ts D^{\ell+1} \} 
-
\mathbb P \{ D^v \geq s,~D^{v-1} = p \} \ts 
\Big|
+ T^2
\nonumber
\\
& \qquad \stackrel{(c)}\leq~
\frac{1}{\epsilon} \ts T \sum_{v=1}^\ell (1-\epsilon)^{\ell-v}
+
\frac{1}{\epsilon} \ts T \sum_{v=1}^\ell (1-\epsilon)^{\ell+1-v}
+ T^2
\nonumber
\\
& \qquad \leq~
\frac{1}{\epsilon} \ts T \ts \Big[ \frac{1}{\epsilon} + \frac{1-\epsilon}{\epsilon} \Big] + T^2,
\label{eqn:cond_gap_second}
\end{align}
where $(b)$ holds by noting that $n_{it}^k(q) \in [0,1]$ and $|a - b| \leq |a-c| + |c-b|$ in the first sum on the left side of the inequality, as well as bounding the second sum on the left side of the inequality by $T^2$, whereas $(c)$ uses Lemma \ref{lem:cond_gap}. Because $1 =  \sum_{p \in \Tcal} \mathbb P \{ D^{k+1} = p \ts | \ts D^k = q\} \geq T \ts \epsilon$, we get $T \leq 1/\epsilon$, so $\frac{1}{\epsilon} \ts T \ts \Big[ \frac{1}{\epsilon} + \frac{1-\epsilon}{\epsilon} \Big] + T^2 \leq \frac{1}{\epsilon^2} \ts \Big[ \frac{1}{\epsilon} + \frac{1-\epsilon}{\epsilon} \Big] + \frac{1}{\epsilon^2} = \frac{2}{\epsilon^3}$. Thus, the result follows from (\ref{eqn:cond_gap_first}) and (\ref{eqn:cond_gap_second}).  \qed

In the next lemma, which is the second result that we are interested in, we give an upper bound on the moment generating function of a Bernoulli random variable.

\begin{lem}
\label{lem:mgf_bern}
If $Z$ is Bernoulli with mean $\mu$, then we have $
\mathbb E\{ e^{\lambda  (Z - \mu)} \} \leq e^{\mu \lambda^2}
$ for all $\lambda \in [0,1]$. 

\end{lem}

\noindent{\it Proof:} Because $e^x$ is convex in $x$, over the interval [0,1], the function $e^x$ lies below the line segment that connects the points $(0,1)$ and $(1,e)$. Thus, we have $e^x \leq 1 + x\ts ( e-1)$ for all $x \in [0,1]$. In this case, noting that $e^\lambda = \int_0^\lambda e^x \ts dx + 1$, for all $\lambda \in [0,1]$, we get $e^\lambda \leq \int_0^\lambda (1+x \ts (e-1)) \ts dx + 1 = \lambda + \frac{1}{2} (e-1) \lambda^2 + 1 \leq 1 + \lambda + \lambda^2$, yielding $e^\lambda - 1 - \lambda \leq \lambda^2$ for all $\lambda \in [0,1]$. Because $Z$ is Bernoulli with mean $\mu$, we have $\mathbb P \{ Z = 1 \} = \mu$, so $
\mathbb E\{ e^{\lambda \ts Z } \} = (1-\mu) +  \mu \ts e^{\lambda} = 1 + \mu \ts (e^\lambda - 1) \leq e^{\mu \ts (e^\lambda - 1)}$, where the last inequality holds because $1+x \leq e^x$ for all $x \in \mathbb R$.  By the last chain of inequalities, we obtain $\mathbb E\{ e^{\lambda(Z-\mu)} \} \leq e^{\mu(e^\lambda -1-\lambda)} \leq e^{\mu \lambda^2}$, where we use the fact that  $e^\lambda -1 - \lambda \leq \lambda^2$. \qed

\section{Choice of the Tuning Parameters}
\label{sec:tuning}

%\vspace{-1mm}

In the proof of Theorem~\ref{thm:perf}, our specific choice of the values for $\gamma$ and $\lambda$ is driven by the form of the availability probability bound in Lemma \ref{lem:avail}. To give intuition into these choices, we consider the case with one stage in the selling horizon so that $K=1$. In this case, the probability in  Lemma \ref{lem:avail} is upper bounded by $e^{c_i  \lambda^2 - (1-\gamma)  c_i  \lambda}$. To make this bound as tight as possible, we find the value of $\lambda$ that minimizes $c_i \ts \lambda^2 - (1-\gamma) \ts c_i \ts \lambda$, which is given by $\frac{1-\gamma}{2}$.~Fixing $\lambda = \frac{1-\gamma}{2}$ and noting that $c_i \geq c_{\min}$, we get $c_i \ts \lambda^2 - (1-\gamma) \ts c_i \ts \lambda = - \frac 14 c_i \ts (1-\gamma)^2 \leq - \frac 14 c_{\min} (1-\gamma)^2$, so the probability in Lemma~\ref{lem:avail} is upper bounded by $e^{-\frac 14 \ts c_{\min} \ts (1-\gamma)^2}$. In this case, by (\ref{eqn:union}), the probability \mbox{$\mathbb P \{ G_{jt}^k = 1 \ts | \ts \Psi_t^k(q) = 1\}$} is lower bounded by $1 - L \ts  e^{-\frac 14 \ts c_{\min} \ts (1-\gamma)^2}$. By the discussion just after (\ref{eqn:apx_rev}), if~$\mathbb P \{ G_{jt}^k = 1 \ts | \ts \Psi_t^k(q) = 1\}  \geq \alpha$, then we have $\frac{\apx}{\Zbar_\lp} \geq \gamma \ts \alpha$. Therefore, we get $\frac{\apx}{\Zbar_\lp} \geq \gamma \ts (1 - L \ts  e^{-\frac 14 \ts c_{\min} \ts (1-\gamma)^2})$, which implies that we can use $\gamma \ts (1 - L \ts  e^{-\frac 14 \ts c_{\min} \ts (1-\gamma)^2})$ as a performance guarantee for our approximate policy. To obtain the best performance guarantee, we can find the value of $\gamma$ that maximizes $\gamma \ts (1 - L \ts  e^{-\frac 14 \ts c_{\min} \ts (1-\gamma)^2})$, but finding the value of $\gamma$ that maximizes $\gamma \ts (1 - L \ts  e^{-\frac 14 \ts c_{\min} \ts (1-\gamma)^2})$ explicitly is difficult. 

%\vspace{-0.5mm}

To get around this difficulty, we make two observations. First, we parameterize the value of $\gamma$ as $\gamma = 1 - 2 \ts c_{\min}^\alpha$ for some choice of $\alpha \in \mathbb R$. In this case, the performance guarantee for our approximate policy takes the form  $\gamma \ts (1 - L \ts  e^{-\frac 14 \ts c_{\min} \ts (1-\gamma)^2}) = (1- 2 \ts c_{\min}^\alpha) \ts (1 - L \ts e^{-c_{\min}^{2\alpha+1}})$. As $c_{\min}$ becomes arbitrarily large, we would like the performance guarantee to be arbitrarily close to one.~If \mbox{$\alpha \in (-\frac 12, 0)$}, then we have $\lim_{c_{\min} \rightarrow \infty} (1- 2 \ts c_{\min}^\alpha) \ts (1 - L \ts e^{-c_{\min}^{2\alpha+1}})  =1$. Thus, we choose the value of $\alpha$ in the interval $(-\frac12 ,0)$.~Second, because $c_{\min} \geq 1$ and $\alpha > -\frac 12$, for large enough values of $c_{\min}$, we have \mbox{$(1- 2 \ts c_{\min}^\alpha) \ts (1 - L \ts e^{- c_{\min}^{2\alpha+1}}) \leq 1 - 2 \ts c_{\min}^\alpha \leq 1 - 2 \ts c_{\min}^{-\frac 12}$}, so the best possible performance guarantee that we can hope to obtain is $1 - 2 \ts \frac{1}{\sqrt{c_{\min}}}$. If we fix $\alpha = \frac{\log \log c_{\min}}{2 \log c_{\min}} - \frac 12$, then $2\alpha + 1 = \frac{\log \log c_{\min}}{\log c_{\min}}$,  so using the identity $x^{(\log \log x) / (2\log x)} = \sqrt{\log x}$, the performance guarantee for the approximate policy becomes $(1- 2 \ts c_{\min}^\alpha)  (1 - L \ts e^{-c_{\min}^{2\alpha+1}}) = \Big(1 - 2 \ts \sqrt{\frac{\log c_{\min}}{c_{\min}}} \Big)  \Big( 1 - \frac{L}{c_{\min}} \Big)$. 
Therefore, our choice of $\gamma$ ensures that our approximate policy gets to have a performance guarantee of  $\Big(1 - 2 \ts \sqrt{\frac{\log c_{\min}}{c_{\min}}} \Big)  \Big( 1 - \frac{L}{c_{\min}} \Big)$, whereas the best possible performance guarantee that we can hope for is $1 - \Ommega\Big(\frac{1}{\sqrt{c_{\min}}}\Big)$.

%\vspace{-0.5mm}

Note that $\Big(1 - 2 \ts \sqrt{\frac{\log c_{\min}}{c_{\min}}} \Big)  \Big( 1 - \frac{L}{c_{\min}} \Big) \geq 1 - 2 \ts \sqrt{\frac{\log c_{\min}}{c_{\min}}} - \frac{L}{c_{\min}} = 1 - \Ommega\Big(\sqrt{\frac{\log c_{\min}}{c_{\min}}} \Big)$, so our choice of $\gamma$ ensures that our approximate policy has a performance guarantee of \mbox{$1 - \Ommega\Big(\sqrt{\frac{\log c_{\min}}{c_{\min}}} \Big)$} and the best possible performance guarantee that we can hope for is \mbox{$1 - \Ommega\Big(\frac{1}{\sqrt{c_{\min}}} \Big)$}. Therefore, our choice of $\gamma$ ensures that we match the best possible performance guarantee up to a logarithmic factor. Using the identity $x^{(\log \log x) / (2\log x)} = \sqrt{\log x}$, if we fix  \mbox{$\alpha = \frac{\log \log c_{\min}}{2 \log c_{\min}} - \frac 12$}, then our choice of $\gamma$ is given by $\gamma = 1 - 2 \ts c_{\min}^\alpha = 1 - 2 \sqrt{ \frac{\log c_{\min}}{c_{\min}}} = 1 - \frac{\sqrt{4 \ts c_{\min} \ts \log c_{\min}}}{c_{\min}}$, which matches the choice of $\gamma$ in the proof of Theorem \ref{thm:perf} when there is one stage in the selling horizon so that $K=1$. The preceding discussion sheds light into our choice of the tuning parameter $\gamma$ for the case with one stage in~the selling horizon. We can use precisely the same reasoning for the case with multiple stages in the selling horizon, but the calculations with multiple stages get slightly more tedious.

\section{Independent Demands in Different Stages}
\label{sec:indep}

We study performance guarantees for an approximate policy from our fluid approximation when the demand random variables $D^1,\ldots,D^K$ in different stages are independent of each other.

\subsection{Approximate Policy and Asymptotic Optimality}

In this section, we elaborate on the fact that having $K$ stages with the demand random variables $D^1,\ldots,D^K$ being independent of each other is not equivalent to having one stage with the number of customer arrivals being given by the random variable $D^1 + \ldots + D^K$. Next, when working with independent demand random variables, we assume that the demands in different stages are \mbox{sub-Gaussian}. We explain that this assumption is milder than assuming that $\mathbb P \{ D^{k+1} = p \} \geq \epsilon$ for all $p \in \Tcal$ and $k \in \Kcal$ for some $\epsilon > 0$, which would have been the analogue of the assumption in Section~\ref{sec:form} when the  demand random variables in different stages are independent of each other. We carefully contrast these two assumptions. Following this discussion, we provide an approximate policy that uses our fluid approximation under independent demands. We give a performance guarantee for this policy under sub-Gaussian demand random variables that are independent of each other.~Lastly, we demonstrate that our fluid approximation has a simplified form under independent demand random variables and other intuitive fluid approximations that make seemingly harmless modifications in our fluid approximation do not yield asymptotically tight upper bounds. %We proceed to discussing how having the demand random variables  $D^1,\ldots,D^K$ independent of each other is not equivalent to having one stage with the number of customer arrivals given by the random variable $D^1 + \ldots + D^K$.

Having $K$ stages with the demand random variables $D^1,\ldots,D^K$ being independent of each other is not equivalent to having one stage with the number of customer arrivals being given by the random variable $D^1 + \ldots + D^K$. In particular, the resolution of the demand uncertainty in the two cases are different. Consider having $K$ stages with the demand random variables $D^1,\ldots,D^K$ being independent of each other. If we are in stage $k$ and there have been $t$ customer arrivals in the current stage, then the probability of having one more customer arrival is $\mbox{$\mathbb P \{ D^k \geq t+1 \ts | \ts D^k \geq t\}$} + \mathbb P \{ D^k = t \ts | \ts D^k \geq t\} \ts\ts \mathbb P \{ D^{k+1} + \ldots + D^K \geq 1 \}$, which corresponds to the probability of having one more customer arrival in the current stage plus the probability of having no more customer arrivals in the current stage and one more customer arrival in the remaining stages. If the stages correspond to, for example, weeks, then this first demand model is aware of the calendar in the sense that it distinguishes between having one more customer arrival in the current week and finishing the current week with no more customer arrivals. On the other hand, consider having one stage with the number of customer arrivals being given by the random variable $D^1 + \ldots + D^K$. If there have been $v$ customer arrivals so far, then the probability of having one more customer arrival is $\mathbb P \{ D^1 + \ldots + D^K \geq v + 1 \ts | \ts D^1 + \ldots + D^K \geq v\}$, which does not explicitly consider the distribution of the random demand over different stages. If the stages correspond to weeks, then this second demand model is not aware of the calendar in the sense that it does not pay attention to the beginning and end of each week. To our knowledge, demand models with uncertainty in the demand resolving sequentially over multiple stages have not been considered. Thus, even when the demands in different stages are independent, our fluid approximation yields new results.

When the demand random variables in different stages are independent, we drop the assumption that $\mathbb P \{ D^{k+1} = p \ts | \ts D^k = q\} \geq \epsilon$ for all $p,q \in \Tcal$ and $k \in \Kcal$ for some $\epsilon > 0$. Instead, we assume that the demand random variable in each stage is sub-Gaussian, so we have $\mathbb E \{ e^{\lambda ( D^k - \mathbb E\{D^k\})} \} \leq e^{\frac{\sigma^2}{200} \lambda^2}$ for all $k \in \Kcal$ and $\lambda \geq 0$ for some $\sigma^2 > 0$. The assumption of sub-Gaussian demand is mild, as the class of sub-Gaussian random variables is large; see, for example, Section 2.1.2 in \cite{Wa19}. The parameter $\frac{\sigma^2}{200}$ is known as the variance proxy because if $D^k$ satisfies the \mbox{sub-Gaussian} assumption, then the variance of $D^k$ is at most $\frac{\sigma^2}{200}$. We scale the variance proxy by $1/200$ for notational uniformity in our proofs.~Under independent demands, we continue using the fluid approximation in (\ref{eqn:lp}), but this linear program admits obvious simplifications under independent demands. For example, we can write the probability \mbox{$\mathbb P \{ D^\ell \geq s, \ts D^{\ell-1} =p \ts|\ts D^k \geq t, \ts D^{k-1} = q \}$} in the first constraint as \mbox{$\mathbb P \{ D^\ell \geq s\} \ts\ts \mathbb P\{ D^{\ell-1} =p\}$} for $\ell < k-1$. To construct an approximate policy from the fluid approximation, we solve problem (\ref{eqn:lp}) once at the beginning of the selling horizon. Letting \mbox{$\xvecbar = (\xbar_{jt}^k(q) : j \in \Jcal,~t,q \in \Tcal,~k \in \Kcal)$}~be an optimal solution, we make the decisions as follows.

{\bf \underline{Approximate Policy under Independent Demands}:}
\\
\indent Using $\gamma \in [0,1]$ to denote a tuning parameter, if we have a request for product $j$ at time period~$t$ in stage $k$, then we are willing to accept the request with probability $\gamma \ts \sum_{q \in \Tcal} \mathbb P \{ D^{k-1} = q\} \ts \frac{\xbar_{jt}^k(q)}{\lambda_{jt}^k}$. If we are willing to accept the request and there are enough resource capacities to accept the request, then we accept the request. Otherwise, we reject. 

In the next theorem, we give a performance guarantee for the approximate policy. We continue using $\apx$ to denote the total expected revenue of the approximate policy. 

\begin{thm}[Independent Demands] 
\label{thm:perf_indep}
Under independent demands, there exists a choice of the tuning parameter $\gamma$ such that we have 
\begin{align*}
\frac{\apx}{\opt} 
~\geq~ 
\frac{\apx}{\Zbar_\lp} 
~\geq~ 
\max\Bigg\{ \frac{1}{4L} , \Bigg( 1 - 4 \ts \frac{\sqrt{(c_{\min} + \sigma^2 \ts  (K-1) ) \log c_{\min}}}{c_{\min}} - \frac{L}{c_{\min}} \Bigg) \Bigg\}.
\end{align*}
\end{thm}

We give the proof of the theorem shortly in Extended Results \ref{sec:perf_indep}. Although the demand random variables in the theorem above are independent, we are unable to use concentration inequalities for sums of independent random variables in the proof. In particular, using $N_{it}^k$ to capture the Bernoulli random variable that takes value one if there is a request for a product at time period $t$ in stage $k$ that uses resource $i$ and the approximate policy is willing to accept the request, we can use $\sum_{\ell=1}^{k-1} \sum_{s \in \Tcal} {\bf 1}(D^\ell \geq s) \ts N_{is}^\ell + \sum_{s=1}^t {\bf 1}(D^k \geq s)  N_{is}^k$ to upper bound the capacity consumption of resource $i$ up to and including time period $t$ in stage $k$. This random variable is the analogue of $\sum_{\ell=1}^{k-1} \sum_{s \in \Tcal} \sum_{q \in \Tcal} \Psi_s^\ell(q) \ts N_{is}^\ell(q) + \sum_{s=1}^t \sum_{q \in \Tcal} \Psi_s^k(q) \ts N_{is}^k(q)$ used right after (\ref{eqn:apx_rev}) in Section~\ref{sec:perf_asymp}. Considering the sum $\sum_{\ell=1}^{k-1} \sum_{s \in \Tcal} {\bf 1}(D^\ell \geq s) \ts N_{is}^\ell + \sum_{s=1}^t {\bf 1}(D^k \geq s)  N_{is}^k$, even if the demands are independent, for fixed stage $\ell$, the terms \mbox{$\{ {\bf 1}(D^\ell \geq s) \ts N_{is}^\ell : s \in \Tcal\}$} in this sum are not independent due to the random variable ${\bf 1}(D^\ell \geq s)$. Thus, we cannot bound the tail probabilities of the capacity consumptions by using concentration inequalities  for sums of independent random variables. We derive our concentration inequalities through moment generating function bounds. 

\vspace{-0.15mm}

When the demands in different stages are dependent, to establish our performance guarantee, we assume that  $\mathbb P \{ D^k = t\} \geq \epsilon$ for all $t \in \Tcal$ and $k \in \Kcal$ for some $\epsilon > 0$. If $\mathbb P \{ D^k = t\} \geq \epsilon$ for all $t \in \Tcal$, then $D^k \leq 1/\epsilon$ with probability one. When the demands in different stages are independent, to establish our performance guarantee, we assume that the demands are sub-Gaussian. By Lemma 5.1 in \cite{DuPa09}, bounded random variables are sub-Gaussian, so the latter assumption is milder, but the milder assumption comes at the expense of having to assume independence between the demands. The interpretation Theorem~\ref{thm:perf_indep} is similar to that of Theorem~\ref{thm:perf}. The approximate policy has a constant-factor performance guarantee, as long as $L$ is fixed. Furthermore, if we scale both the number of stages and resource capacities with the same rate $\theta$, then the relative gap between the total expected revenue of the approximate policy and the optimal total expected revenue converges to one with rate $1 - \frac{1}{\sqrt \theta}$.

\vspace{-0.15mm}

There are important differences in the analyses of our approximate policies under dependent and independent demands. Under independent demands, we construct concentration inequalities for a random variable of the form $\sum_{k=1}^K \sum_{t=1}^{D^k} B_t^k$, where $\{ B_t^k :t \in \Tcal,~k \in \Kcal\}$ are independent Bernoulli random variables. The main difficulty in constructing these concentration inequalities is that the quantity $D^k$ is random. Thus, we deal with sums of random numbers of Bernoulli random variables, while the Bernoulli random variables in the sum being independent of each other. On the other hand, under dependent demands, we construct concentration inequalities for a random variable of the form $\sum_{k=1}^K \sum_{t=1}^{D^k} B_t^k(D^{k-1})$, where \mbox{$\{ B_t^k(D^{k-1}) :t \in \Tcal,~k \in \Kcal\}$} are dependent Bernoulli random variables, because the distribution of $B_t^k(D^{k-1})$ depends on $D^{k-1}$ and $\{D^1,\ldots,D^K\}$ are dependent on each other. Thus, we deal with sums of random numbers of Bernoulli random variables, while the Bernoulli random variables in the sum being dependent on each other.

\subsection{Performance Guarantee}
\label{sec:perf_indep}

We give a proof for Theorem \ref{thm:perf_indep}. Our starting point is similar to the outline of Section \ref{sec:perf_asymp}, but we will construct different concentration inequalities to exploit the sub-Gaussian property.

{\bf \underline{Preliminary Random Variables and Availability Probabilities}:}
\\
\indent We define four classes of Bernoulli random variables for each $k \in \Kcal$ and $t \in \Tcal$. We worked with the analogues of these random variables in Section \ref{sec:perf_asymp}.

\indent {$\bullet$ \it {Demand in Each Stage}.} The random variable $\Psi_t^k$ takes value one if we reach time period $t$ in stage $k$ before this stage is over. In other words, recalling that we use ${\bf 1}(\cdot)$ to denote the indicator function, we simply have $\Psi_t^k = {\bf 1}(D^k \geq t)$. 

\indent {$\bullet$ \it {Product Request}.} For each $j \in \Jcal$, the random variable $A_{jt}^k$ takes value one if the customer arriving at time period $t$ in stage $k$ requests product $j$. We have $\mathbb P \{ A_{jt}^k  =1 \} = \lambda_{jt}^k$. The random variables $\{ A_{jt}^k : t \in \Tcal,~k \in \Kcal\}$ are independent of each other. 

\indent {$\bullet$ \it {Policy Decision}.} For each $j \in \Jcal$, the random variable $X_{jt}^k$ takes value one if the approximate policy is willing to accept a request for product $j$ at time period in stage $k$. By the definition of the approximate policy under independent demands, $\mathbb P \{ X_{jt}^k = 1\} = \gamma \ts \sum_{q \in \Tcal} \mathbb P\{ D^{k-1} = q\} \ts \frac{\xbar_{jt}^k(q)}{\lambda_{jt}^k}$.

\indent {$\bullet$ \it {Availability}.} For each $j \in \Jcal$, the random variable $G_{jt}^k$ takes value one if we have enough capacity to accept a request for product $j$ at time period $t$ in stage $k$ under the approximate policy. Instead of calculating the probability $\mathbb P \{ G_{jt}^k =1\}$, we will lower bound $\mathbb P \{ G_{jt}^k =1 \ts | \ts \Psi_t^k = 1\}$.

Under independent demands, we drop the assumption that $\mathbb P \{ D^{k+1} = p \ts | \ts D^k = q\} \geq \epsilon$ for all $p,q \in \Tcal$ and $k \in \Kcal$ for some $\epsilon > 0$. Instead, we assume that $\mathbb E \{ e^{\lambda ( D^k - \mathbb E\{D^k\})} \} \leq e^{\frac{\sigma^2}{200}  \lambda^2}$ for all $k \in \Kcal$ and $\lambda \geq 0$ for some $\sigma^2 > 0$. Throughout this section, when we refer to the approximate policy, we mean the approximate policy under independent demands as described in Section \ref{sec:indep}. Under the approximate policy, the sales for product $j$ at time period $t$ in stage $k$ is given by the random variable $\Psi_t^k \ts G_{jt}^k \ts A_{jt}^k \ts X_{jt}^k$. Taking expectations, the expected sales for product $j$ at time period $t$ in stage $k$ is $\mathbb P\{ \Psi_t^k = 1\} \ts\ts \mathbb P \{ G_{jt}^k = 1 \ts | \ts \Psi_t^k = 1\} \ts\ts \mathbb P \{ A_{jt}^k = 1\} \ts\ts \mathbb P \{ X_{jt}^k = 1\}$, where we use the fact that the random variables $A_{jt}^k$ and $X_{jt}^k$ are independent of the demands and remaining capacities of the resources. Noting that $\mathbb P \{ A_{jt}^k = 1\} = \lambda_{jt}$ and $\mathbb P \{ X_{jt}^k = 1\} = \gamma \ts \sum_{q \in \Tcal} \mathbb P\{ D^{k-1} = q\} \ts \frac{\xbar_{jt}^k(q)}{\lambda_{jt}^k}$, the last expectation is given by $\mathbb P \{ D^k \geq t\} \ts\ts \mathbb P \{ G_{jt}^k = 1 \ts | \ts \Psi_t^k = 1\} \ts \sum_{q \in \Tcal} \gamma \ts\ts \mathbb P \{ D^{k-1}  =q \} \ts\ts \xbar_{jt}^k(q)$. Therefore, the total expected revenue of the approximate policy is
\begin{align*}
\apx
~=~ 
\sum_{k \in \Kcal} \sum_{t \in \Tcal} \sum_{j \in \Jcal} \sum_{q \in \Tcal} f_j \ts \mathbb P \{ D^k \geq t\} \ts\ts \mathbb P \{ G_{jt}^k = 1 \ts | \ts \Psi_t^k = 1\} \ts\ts \gamma \ts \mathbb P\{ D^{k-1} = q \} \ts \xbar_{jt}^k(q).
\end{align*}
Because $\xvecbar$ is an optimal solution to problem (\ref{eqn:lp}), the optimal objective value of problem (\ref{eqn:lp}) is \mbox{$\Zbar_\lp = \sum_{k \in \Kcal} \sum_{t \in \Tcal} \sum_{j \in \Jcal} \sum_{q \in \Tcal} f_j \ts \mathbb P \{ D^k \geq t\} \ts  \mathbb P\{ D^{k-1} = q \} \ts \xbar_{jt}^k(q)$}, where we use the fact that the demands in different stages are independent, so $\mathbb P \{ D^k \geq t,~D^{k-1} = q\} = \mbox{$\mathbb P \{ D^k \geq t\}$} \ts \mbox{$\mathbb P \{D^{k-1} = q\}$}$. Thus, comparing the expressions for the total expected revenue $\apx$ of the approximate policy and the optimal objective value $\Zbar_\lp$ of the linear program in (\ref{eqn:lp}), it follows that if we can show that $\mathbb P\{ G_{jt}^k = 1 \ts | \ts \Psi_t^k = 1\} \geq \alpha$, then we obtain $\apx \geq \gamma \ts \alpha \Zbar_\lp$. We focus on lower bounding the availability probability $\mathbb P\{ G_{jt}^k = 1 \ts | \ts \Psi_t^k = 1\}$. Under the approximate policy, the capacity consumption of resource $i$ at time period $t$ in stage $k$ is given by the random variable $\sum_{j \in \Jcal} a_{ij} \ts \Psi_t^k \ts G_{jt}^k \ts A_{jt}^k \ts X_{jt}^k$, which implies that the capacity consumption of resource $i$ at time period $t$ in stage $k$ is upper bounded by $\sum_{j \in \Jcal} a_{ij} \ts \Psi_t^k \ts A_{jt}^k \ts X_{jt}^k$.~Letting $N_{it}^k = \sum_{j \in \Jcal} a_{ij} \ts A_{jt}^k \ts X_{jt}^k$, we write the upper bound on the  capacity consumption of resource $i$ at time period $t$ in stage $k$ as $\Psi_t^k \ts N_{it}^k$. In this case, if we have~$\sum_{\ell = 1}^{k-1} \sum_{s \in \Tcal} \Psi_s^\ell \ts N_{is}^\ell + \sum_{s=1}^t \Psi_s^k \ts N_{is}^k < c_i$, then the total capacity consumption of resource $i$ up to and including time period $t$ in stage $k$ does not exceed the capacity of the resource, in which case, we must have capacity available for resource $i$ at time period $t$ in stage $k$. Therefore, using the same line of reasoning that we followed to obtain the chain of inequalities in (\ref{eqn:union}), we can lower bound the availability probability as 
\begin{align}
& \mathbb P \{ G_{jt} = 1 \ts | \ts \Psi_t^k = 1\} 
~\geq~ 
1 - \sum_{i \in \Lcal_j} \mathbb P \bigg\{\sum_{\ell=1}^{k-1} \sum_{s \in \Tcal}  \Psi_s^\ell \ts N_{is}^\ell + \sum_{s=1}^t N_{is}^k \ts \geq \ts c_i \ts \Big| \ts D^k \geq t  \bigg\}.
\label{eqn:avail_indep}
\end{align}
The discussion so far closely followed the one at the beginning of Section \ref{sec:perf_asymp}, but we need to deviate from that discussion to exploit the sub-Gaussian assumption.

{\bf \underline{Moment Generating Function Bounds}:}
\\
\indent Note that $N_{it}^k$ is a Bernoulli random variable. We define $n_{it}^k = \mathbb E\{ N_{it}^k\}$. In the next lemma, we bound the moment generating function of the squared deviation of the demand around its mean. 

\begin{lem}
\label{lem:mgf_squared}
Letting $\mu^k = \mathbb E\{ D^k\}$, for all $k \in \Kcal$, if $|\lambda| \leq \frac{2}{\sigma}$, then we have
\begin{align*}
\mathbb E\{ e^{\lambda^2 (D^k - \mu^k)^2} \} ~\leq~ e^{\frac{1}{4} \ts \sigma^2 \lambda^2}.
\end{align*}
\end{lem}

\noindent{\it Proof:} In an auxiliary lemma, labeled as Lemma \ref{lem:one_to_square}, given at the end of this section, we show that if $Z$ is a mean-zero sub-Gaussian random variable with variance proxy $M^2$ so that $\mathbb E\{ e^{\lambda Z} \} \leq e^{M^2\lambda^2 }$ for all $\lambda \in \mathbb R$, then we have $\mathbb E\{ e^{\theta^2 Z^2} \} \leq e^{(7M\theta)^2}$ for all $|\theta| \leq \frac{1}{7M}$. Because $D^k - \mu^k$ is a \mbox{mean-zero} sub-Gaussian random variable with variance proxy $\frac{\sigma^2}{200}$, using Lemma \ref{lem:one_to_square} with $M = \frac{\sigma}{\sqrt{200}}$, we have $\mathbb E\{ e^{\theta^2 \ts (D^k - \mu^k)^2} \} \leq e^{(7 \frac{\sigma}{\sqrt{200}} \theta)^2}$ for all $|\theta| \leq \frac{1}{7 \ts \sigma / \sqrt{200}}$. We have $\frac{1}{7 \ts \sigma / \sqrt{200}} \geq \frac{2}{\sigma}$. Thus, if we have $|\theta| \leq \frac{2}{\sigma}$, then we also have $|\theta| \leq \frac{1}{7 \ts \sigma / \sqrt{200}}$, in which case, $\mathbb E\{ e^{\theta^2 \ts (D^k - \mu^k)^2} \} \leq e^{(7 \frac{\sigma}{\sqrt{200}} \theta)^2} = e^{\frac{49}{200} \sigma^2 \theta^2} \leq e^{\frac{1}{4} \sigma^2 \theta^2}$. \qed

In the next lemma, we use the lemma above to bound the moment generating function of $\sum_{s \in \Tcal} {\bf 1}(D^k \geq s) \ts n_{is}^k$, which will be key to bounding the availability probabilities.

\begin{lem}
\label{lem:mgf_exp_cap}
For all $k \in \Kcal$, $i \in \Lcal$ and $\lambda \in \mathbb R$, we have
\begin{align*}
\mathbb E \{ e^{\lambda \sum_{s \in \Tcal} n_{is}^k \ts {\bf 1}(D^k \geq s)}\} \ts \leq \ts  e^{\sigma^2 \ts \lambda^2 + \lambda \sum_{s \in \Tcal} n_{is}^k \ts \mathbb P \{ D^k \geq s \}}.
\end{align*}
\end{lem}

\noindent{\it Proof:} For the moment, we assume that $|\lambda| \leq \frac{1}{\sigma}$. We will need two inequalities. 
First, we use the random variable $R^k$ to denote an independent and identically distributed copy of $D^k$. Noting that $\sum_{s \in \Tcal} n_{is}^k \ts {\bf 1}(D^k \geq s) = \sum_{s=1}^{D^k} n_{is}^k$, using the fact that $n_{is}^k \in [0,1]$, we get the chain of inequalities $|\sum_{s \in \Tcal} n_{is}^k \ts ({\bf 1}(D^k \geq s) - {\bf 1}(R^k \geq s)) \ts | \ts = \ts | \sum_{s=1}^{D^k} n_{is}^k - \sum_{s=1}^{R^k} n_{is}^k \ts | \ts \leq \ts |D^k - R^k|$. It is simple to verify the inequality $(a-b)^2 \ts \leq \ts 2(a-\delta)^2 + 2(b-\delta)^2$. Thus, letting $\mu^k = \mathbb E\{ D^k\} = \mathbb E \{ R^k\}$, we obtain $[ \ts \! \sum_{s \in \Tcal} n_{is}^k \ts ({\bf 1}(D^k \geq s) - {\bf 1}(R^k \geq s)) \ts ]^2 \ts \leq \ts (D^k - R^k)^2 \ts \leq \ts  2 \ts ( D^k - \mu^k)^2 + 2 \ts (R^k - \mu^k)^2$. Second, because $|\lambda| \leq \frac{1}{\sigma}$,  we have $|\sqrt 2 \ts \lambda| \leq \frac{\sqrt 2}{\sigma} \leq \frac2\sigma$, so by Lemma \ref{lem:mgf_squared}, $\mathbb E\{ e^{2\lambda^2(D^k - \mu^k)^2} \} \leq e^{\frac{1}{2} \ts \sigma^2 \lambda^2}$. %By the same reasoning, we have $\mathbb E\{ e^{2\lambda^2(R^k - \mu^k)^2} \} \leq e^{\frac{1}{2} \ts \sigma^2 \lambda^2}$ as well. 
Thus, we have
\begin{align}
& \mathbb E \{ e^{\lambda^2 [\ts\!\sum_{s \in \Tcal} n_{is}^k \ts ({\bf 1}(D^k \geq s) - \mathbb P \{ D^k \geq s \}) ]^2}\} 
~=~
\mathbb E \{ e^{\lambda^2 [\ts\!\sum_{s \in \Tcal} n_{is}^k \ts ({\bf 1}(D^k \geq s) - \mathbb E \{ {\bf 1}(R^k \geq s) \}) ]^2}\} 
\nonumber
\\
~& \qquad \quad \stackrel{(a)}\leq~
\mathbb E \{ e^{\lambda^2 [\ts\!\sum_{s \in \Tcal} n_{is}^k \ts ({\bf 1}(D^k \geq s) - {\bf 1}(R^k \geq s)) ]^2}\}
~\stackrel{(b)}\leq~
\mathbb E \{ e^{2 \lambda^2 (D^k - \mu^k)^2}\} \ts \ts
\mathbb E \{ e^{2 \lambda^2 (R^k - \mu^k)^2}\} 
~\stackrel{(c)}\leq~
e^{\sigma^2 \ts \lambda^2},
\label{eqn:mgf_exp_cap}
\end{align}
where $(a)$ follows by the Jensen inequality and the fact that $e^{\lambda^2 x^2}$ is convex in $x$, $(b)$ uses the first inequality in this paragraph, as well as the fact that the random variables $D^k$ and $R^k$ are independent and $(c)$ uses the second inequality in this paragraph. In an auxiliary lemma, labeled as Lemma \ref{lem:square_to_one}, given at the end of this section, we show that if $Z$ is a mean-zero random variable that satisfies $\mathbb E\{ e^{\lambda^2 Z^2} \} \leq e^{M^2  \lambda^2}$ for all \mbox{$|\lambda| \leq \frac{1}{M}$}, then we have \mbox{$\mathbb E\{ e^{\theta Z} \} \leq e^{M^2 \theta^2}$} for all $\theta \in \mathbb R$. By (\ref{eqn:mgf_exp_cap}), we have \mbox{$\mathbb E \{ e^{\lambda^2 [\ts\!\sum_{s \in \Tcal} n_{is}^k \ts ({\bf 1}(D^k \geq s) - \mathbb P \{ D^k \geq s \}) ]^2}\} \ts \leq \ts e^{\sigma^2 \lambda^2}$} for all \mbox{$|\lambda| \leq \frac{1}{\sigma}$}. Thus, using Lemma \ref{lem:square_to_one} with \mbox{$M = \sigma$}, $\mathbb E \{ e^{\theta \sum_{s \in \Tcal} n_{is}^k \ts ({\bf 1}(D^k \geq s) - \mathbb P \{ D^k \geq s \})}\} \ts \leq \ts e^{\sigma^2 \theta^2}$ for all $\theta \in \mathbb R$, which is the desired result. 
\qed

{\bf \underline{Performance Guarantee for the Approximate Policy}:}
\\
\indent In the next lemma, we use the moment generating function bound in Lemma \ref{lem:mgf_exp_cap} to give a lower bound on the availability probabilities on the right side of (\ref{eqn:avail_indep}).

\begin{lem}
\label{lem:avail_indep}
For all $k \in \Kcal$,  $i \in \Lcal$, $t \in \Tcal$ and $\lambda \in [0,1]$, we have 
\begin{align*}
\mathbb P \bigg\{ \sum_{\ell=1}^{k-1} \sum_{s \in \Tcal} {\bf 1}(D^\ell \geq s) \ts N_{is}^\ell + \sum_{s=1}^t N_{is}^k \ts \geq \ts c_i \ts \Big| \ts D^k \geq t \bigg\}
~\leq~ 
e^{(c_i + 4 \ts \sigma^2 (k-1) )\lambda^2 - (1-\gamma) \ts c_i \lambda}.
\end{align*}
\end{lem}

\noindent{\it Proof:} Given $D^k$, ${\bf 1}(D^k \geq s) \ts N_{is}^k$ is a Bernoulli random variable with mean ${\bf 1}(D^k \geq s) \ts n_{is}^k$, so by Lemma \ref{lem:mgf_bern}, we have $\mathbb E\{ e^{\lambda \ts {\bf 1}(D^k \geq s) \ts (N_{is}^k - n_{is}^k)} \ts | \ts D^k \} \leq e^{\lambda^2 \ts {\bf 1}(D^k \geq s) \ts n_{is}^k}$. In this case, we obtain
\begin{align*}
\mathbb E\{ e^{\lambda \sum_{s \in \Tcal} {\bf 1}(D^k \geq s) \ts N_{is}^k} \ts | \ts D^k \}
=
e^{\lambda \sum_{s \in \Tcal} {\bf 1}(D^k \geq s) \ts n_{is}^k} \ts 
\mathbb E\{ e^{\lambda \sum_{s \in \Tcal} {\bf 1}(D^k \geq s) \ts (N_{is}^k-n_{is}^k)} \ts | \ts D^k \}
\leq
e^{(\lambda+ \lambda^2) \sum_{s \in \Tcal} {\bf 1}(D^k \geq s) \ts n_{is}^k},
\end{align*}
where we use the fact that $\{ N_{is}^k : s \in \Tcal\}$ are independent of each other. Taking expectations in the chain of inequalities above, we obtain $\mathbb E\{ e^{\lambda \sum_{s \in \Tcal} {\bf 1}(D^k \geq s) \ts N_{is}^k}\} \leq \mathbb E\{ e^{(\lambda+ \lambda^2) \sum_{s \in \Tcal} {\bf 1}(D^k \geq s) \ts n_{is}^k} \}$. Using the same argument, we have $\mathbb E\{ e^{\lambda \ts \sum_{s=1}^t N_{is}^k} \} \leq e^{(\lambda^2 + \lambda) \sum_{s=1}^t n_{is}^k}$. On the other hand, by the discussion in Section \ref{sec:lp}, we can replace the probability $\mathbb P \{ D^\ell \geq s, \ts D^{\ell-1} =p \ts|\ts D^k \geq t, \ts D^{k-1} = q \}$ in the first constraint in problem (\ref{eqn:lp}) with $\mathbb P \{ D^\ell \geq s, \ts D^{\ell-1} =p \ts|\ts \ts D^{k-1} = q \}$. Furthermore, $\xvecbar$ is a feasible solution to problem (\ref{eqn:lp}), so it satisfies the first constraint. Multiplying this constraint with \mbox{$\mathbb P \{ D^{k-1} = q\}$} and adding over all $q \in \Tcal$, we get $\sum_{\ell=1}^{k-1} \sum_{s \in \Tcal} \sum_{p \in \Tcal} \sum_{j \in \Jcal} a_{ij} \mbox{$\mathbb P \{ D^\ell \geq s,~D^{\ell-1} = p\}$} \ts \xbar_{js}^\ell(p) + \sum_{s=1}^t \sum_{j \in \Jcal} \sum_{q \in \Tcal} a_{ij} \ts \mathbb P \{ D^{k-1} = q\} \ts \xbar_{js}^k(q) \leq c_i$. By the definition of $N_{it}^k$, we also have the identity \mbox{$n_{it}^k = \mathbb E\{ N_{it}^k \} = \gamma \sum_{j \in \Jcal} a_{ij} \sum_{q \in \Tcal} \mathbb P \{ D^{k-1} = q \} \ts \xbar_{jt}^k(q)$}. Thus, noting that the demands in different stages are independent of each other, so that $\mathbb P\{ D^\ell \geq s,~ D^{\ell-1} = p \} = \mathbb P\{ D^\ell \geq s\} \ts\ts \mbox{$\mathbb P\{ D^{\ell-1} = p \}$}$, the last inequality yields  $\sum_{\ell=1}^{k-1} \sum_{s \in \Tcal} \mathbb P \{ D^\ell \geq s\} \ts n_{is}^\ell + \sum_{s=1}^t n_{is}^k \leq \gamma \ts c_i$. Noting that the random variables $\{ N_{it}^k: t \in \Tcal,~k \in \Kcal\}$ and $\{ D^k : k \in \Kcal\}$ are all independent of each other, we get
\begin{align*}
& \mathbb P \bigg\{ \sum_{\ell=1}^{k-1} \sum_{s \in \Tcal} {\bf 1}(D^\ell \geq s) \ts N_{is}^\ell + \sum_{s=1}^t N_{is}^k \ts \geq \ts c_i \ts \Big| \ts D^k \geq t \bigg\}
~=~ 
\mathbb P \{ e^{\lambda \ts ( \sum_{\ell=1}^{k-1} \sum_{s \in \Tcal} {\bf 1}(D^\ell \geq s) \ts N_{is}^\ell + \sum_{s=1}^t N_{is}^k )} \ts \geq \ts e^{\lambda \ts c_i} \}
\\
~& \qquad \qquad \qquad \qquad \qquad \qquad \stackrel{(a)}\leq~
\frac{1}{e^{\lambda \ts c_i}} \ts 
\mathbb E \{ e^{\lambda \ts ( \sum_{\ell=1}^{k-1} \sum_{s \in \Tcal} {\bf 1}(D^\ell \geq s) \ts N_{is}^\ell + \sum_{s=1}^t N_{is}^k )} \}
\\
~& \qquad \qquad \qquad \qquad \qquad \qquad =~
\frac{1}{e^{\lambda \ts c_i}} \ts
\ts \prod_{\ell=1}^{k-1} \mathbb E \{ e^{\lambda \ts ( \sum_{s \in \Tcal} {\bf 1}(D^\ell \geq s) \ts N_{is}^\ell} \}  \ts\ts \mathbb E\{ e^{\lambda \sum_{s=1}^t N_{is}^k } \}
\\
~& \qquad \qquad \qquad \qquad \qquad \qquad \stackrel{(b)}\leq~
\frac{1}{e^{\lambda \ts c_i}} \ts
\ts \prod_{\ell=1}^{k-1} \mathbb E \{ e^{(\lambda^2+\lambda) \sum_{s \in \Tcal} {\bf 1}(D^\ell \geq s) \ts n_{is}^\ell} \}  \ts\ts e^{(\lambda^2 + \lambda) \sum_{s=1}^t n_{is}^k}
\\
~& \qquad \qquad \qquad \qquad \qquad \qquad \stackrel{(c)}\leq~
\frac{1}{e^{\lambda \ts c_i}} \ts
\ts \prod_{\ell=1}^{k-1}
e^{4 \sigma^2  \lambda^2 + (\lambda^2 + \lambda) \sum_{s \in \Tcal} \mathbb P \{ D^\ell \geq s \} \ts n_{is}^\ell}  \ts\ts e^{(\lambda^2 + \lambda) \sum_{s=1}^t n_{is}^k}
\\
~& \qquad \qquad \qquad \qquad \qquad \qquad =~
\frac{1}{e^{\lambda \ts c_i}} \ts
e^{4 \sigma^2 (k-1) \lambda^2 } \ts\ts e^{(\lambda^2 + \lambda) (\sum_{\ell = 1}^{k-1} \sum_{s \in \Tcal} \mathbb P \{ D^\ell \geq s \} \ts n_{is}^\ell + \sum_{s=1}^t n_{is}^k)}
\\
~& \qquad \qquad \qquad \qquad \qquad \qquad \stackrel{(d)}\leq~ 
e^{(c_i + 4 \ts \sigma^2 (k-1) )\lambda^2 - (1-\gamma) \ts c_i \lambda},
\end{align*}
where $(a)$ is the Markov inequality, $(b)$ is by the discussion earlier in the proof, $(c)$ uses Lemma~\ref{lem:mgf_exp_cap} and $4 \lambda^2 \geq (\lambda^2 + \lambda)^2$ for $\lambda \in [0,1]$ and $(d)$ holds by $\sum_{\ell=1}^{k-1} \sum_{s \in \Tcal} \mathbb P \{ D^\ell \geq s\} \ts n_{is}^\ell + \sum_{s=1}^t n_{is}^k \leq \gamma \ts c_i$. \qed

Using specific values for $\gamma$ and $\lambda$ in Lemma \ref{lem:avail_indep} will provide a lower bound for the availability probabilities, which, in turn, will yield a performance guarantee for the approximate policy. 

{\it \underline{Proof of Theorem \ref{thm:perf_indep}}:}
\\
\indent Identifying $4 \ts \sigma^2$ with $\frac{3}{\epsilon^6}$, the bounds in Lemmas \ref{lem:avail} and  \ref{lem:avail_indep} have the same form. Thus, choosing $\gamma = 1 - \frac{\sqrt{4 \ts (c_{\min} + 4 \ts \sigma^2 \ts (K-1)) \log c_{\min}}}{c_{\min}}$ and $\lambda = \frac{(1-\gamma) \ts c_i}{2 (c_i + 4 \ts \sigma^2 \ts (K-1))}$, following the proof of Theorem \ref{thm:perf}, we get
\begin{align*}
& \frac{\apx}{\Zbar_\lp} 
~\geq~ 
\Bigg( 1 - 4 \frac{\sqrt{(c_{\min} + \sigma^2 \ts (K-1)) \log c_{\min}}}{c_{\min}}  - \frac{L}{c_{\min}} \Bigg).
\end{align*}
The proof of the inequality above precisely follows the argument at the end of Section \ref{sec:perf_asymp} line by line. On the other hand, we can precisely follow the argument in Extended Results \ref{sec:constant_factor} line by line to show that $\frac{\apx}{\Zbar_\lp} \geq \frac{1}{4L}$. In this case, the desired result follows by collecting these two inequalities together, as well as noting that we have $\frac{\apx}{\opt} \geq \frac{\apx}{\Zbar_\lp}$ by Theorem \ref{thm:bound}. \qed

In the proofs of Lemmas \ref{lem:mgf_squared} and \ref{lem:mgf_exp_cap}, we used two results on sub-Gaussian random variables. We used Lemma \ref{lem:one_to_square} in the proof of Lemma \ref{lem:mgf_squared} and Lemma \ref{lem:square_to_one} in the proof of Lemma \ref{lem:mgf_exp_cap}. 

\begin{lem}
\label{lem:one_to_square}
If $Z$ is a mean-zero sub-Gaussian random variable such that $\mathbb E\{ e^{\lambda Z} \} \leq e^{M^2\lambda^2 }$ for all $\lambda \in \mathbb R$, then we have $\mathbb E\{ e^{\theta^2 Z^2} \} \leq e^{(7M\theta)^2}$ for all $|\theta| \leq \frac{1}{7M}$.
\end{lem}

\noindent{\it Proof:} By the discussion that follows Definition 2.2 in \cite{Wa19}, if $\mathbb E\{ e^{\lambda Z} \} \leq e^{M^2\lambda^2 }$ for all $\lambda \in \mathbb R$, then $\mathbb P\{ |Z| \geq t\} \leq 2 e^{-\frac{t^2}{4M^2}}$ for all $t \geq 0$. Letting $\| Z\|_{L^p} = (\mathbb E\{ |Z|^p\})^{1/p}$ for $p \in \mathbb Z_{++}$, we get
\begin{align*}
& \mathbb E\{ |Z|^p \} 
~=~
\int_0^\infty \mathbb P \{ |Z|^p \geq u \} \ts du 
~\stackrel{(a)}=~
\int_0^\infty \mathbb P \{ |Z| \geq t\} \ts p \ts t^{p-1} \ts dt
~\leq~
2 \!\! \int_0^\infty e^{-\frac{t^2}{4M^2}}\ts p \ts t^{p-1} \ts dt
\\
& ~~ \stackrel{(b)}=~ 
2p \!\! \int_0^\infty \!\!\!\!\! e^{-z} (2M\sqrt{z})^{p-1} \ts \frac{M}{\sqrt z} \ts dz
~=~
p \ts (2M)^p \!\! \int_0^\infty \!\!\!\!\! e^{-z} z^{\frac{p}{2}-1} \ts  dz
~\stackrel{(c)}=~
p \ts (2M)^p \Gamma(p/2) 
~\stackrel{(d)}\leq~
p \ts (2M)^p \ts (p/2)^{p/2},
\end{align*}
where $(a)$ uses the change of variables $u = t^p$, $(b)$  uses the change of variables $z = t^2 / (4M^2)$, $(c)$ is the definition of the gamma function and $(d)$ holds because the gamma function satisfies \mbox{$\Gamma(x) \leq x^x$}.~By~the chain of inequalities above, we obtain $\| Z\|_{L^p} \leq \ts p^{1/p} \ts M \sqrt {2p}$. It is simple to verify that $\max_{p \in \mathbb Z_{++}} p^{1/p} = 3^{1/3}$, so the last inequality yields $\| Z\|_{L^p} \leq \ts 3^{1/3} \ts \ts M \sqrt {2p} \leq 2.04 \ts M \sqrt p$ for all $p \in \mathbb Z_{++}$. In this case, we have $\mathbb E\{ |Z|^{2p}\} \leq (2.04 \ts M \sqrt{2p})^{2p} = (2.04^2 \ts M^2 \ts 2 \ts p)^p \leq (8.33 M^2  p)^p$. Also, by our assumption that $|\theta| \leq \frac{1}{7M}$, we get $8.33 \ts e \ts M^2 \theta^2 \leq \frac{8.33\ts e}{49} \leq \frac 12$. Lastly, we can verify the inequality $e^{2x} \geq \frac{1}{1-x}$ for all $x \in [0,\frac 12]$. Thus, using the Taylor series expansion of $e^x$, we have
\begin{multline*}
\mathbb E\{ e^{\theta^2  Z^2} \}
~=~ 
1 + \sum_{p=1}^\infty \frac{\theta^{2p} \ts \mathbb E\{ |Z|^{2p} \}}{p!}
~\leq~
1 + \sum_{p=1}^\infty \frac{(\theta^2 \ts 8.33 M^2  p )^p}{p!}
~\stackrel{(e)}\leq~
1 + \sum_{p=1}^\infty \frac{(\theta^2 \ts 8.33 M^2  p )^p}{(p/e)^p}
\\
~=~
1 + \sum_{p=1}^\infty (\theta^2 \ts 8.33 \ts e M^2 )^p
~\stackrel{(f)}=~
\frac{1}{1-\theta^2 \ts 8.33 \ts e M^2}
~\stackrel{(g)}\leq~
e^{(7M\theta)^2},
\end{multline*}
where $(e)$ holds because $p! \geq (p/e)^p$ by the Sterling approximation, $(f)$ follows by $8.33 \ts e \ts M^2 \theta^2 \leq \frac{1}{2}$ and $(g)$ holds because $8.33 \ts e \ts M^2 \theta^2 \leq \frac 12$, so $\frac{1}{1-\theta^2 \ts 8.33 \ts e M^2} \leq e^{2 \ts \theta^2 8.33 e M^2} \leq e^{49 \ts \theta^2 M^2}$. \qed

In the proof, we round the fractions up to two decimals to strike a balance between clarity and tightness. We can round the fractions up to integers, in which case, the bounds would be looser.

\begin{lem}
\label{lem:square_to_one}
If $Z$ is a mean-zero random variable that satisfies $\mathbb E\{ e^{\lambda^2 Z^2} \} \leq e^{M^2  \lambda^2}$ for all \mbox{$|\lambda| \leq \frac{1}{M}$}, then we have \mbox{$\mathbb E\{ e^{\theta Z} \} \leq e^{M^2 \theta^2}$} for all $\theta \in \mathbb R$.
\end{lem}

\noindent{\it Proof:} Consider any $\theta \in \mathbb R$. We claim that $\mathbb E\{ e^{\frac{\theta}{M} Z}\} \leq e^{\theta^2}$. First, assume that $|\theta| \leq 1$. Because $|\frac{\theta}{M}| \leq \frac{1}{M}$,  using the assumption of the lemma with $\lambda = \frac{\theta}{M}$, we get $\mathbb E\{ e^{\frac{\theta^2}{M^2}Z^2} \} \leq e^{\theta^2}$. For any $x \in \mathbb R$, we have the inequality $e^x \leq x + e^{x^2}$, in which case, noting that $\mathbb E\{ Z\} = 0$, we obtain  $\mathbb E\{ e^{\frac{\theta}{M} Z}\} \leq \mathbb E\{\frac{\theta}{M} Z\} + \mathbb E\{ e^{\frac{\theta^2}{M^2} Z^2} \} \leq e^{\theta^2}$, establishing the claim when $|\theta| \leq 1$. 
Second, assume that $|\theta| \geq 1$. Because $\frac{1}{\sqrt{2} M} \leq \frac{1}{M}$, using the assumption of the lemma with $\lambda = \frac{1}{\sqrt{2} M}$, we get $\mathbb E\{ e^{\frac{1}{2M^2}Z^2}\} \leq e^{\frac{1}{2}}$. We have the standard inequality $\frac{1}{2} a^2 + \frac{1}{2} b \geq ab$, which implies that $\frac{\theta^2}{2} + \frac{x^2}{2M^2} \geq \frac{\theta \ts x}{M}$ for all $x \in \mathbb R$. In this case, we obtain the chain of inequalities $\mathbb E\{ e^{\frac{\theta}{M} Z}\} \leq e^{\frac{\theta^2}{2}} \ts \mathbb E\{ e^{\frac{1}{2M^2}Z^2}\} \leq e^{\frac{\theta^2}{2}} \ts e^{\frac{1}{2}} \leq e^{\theta^2}$, where the last inequality holds because we have $|\theta| \geq 1$. Therefore, the claim holds when $|\theta| \geq 1$ as well. 

By the claim established in the previous paragraph, we have $\mathbb E\{ e^{\frac{\theta}{M} Z}\} \leq e^{\theta^2}$ for all $\theta \in \mathbb R$, which is equivalent to having \mbox{$\mathbb E\{ e^{\theta Z} \} \leq e^{M^2 \theta^2}$} for all $\theta \in \mathbb R$. \qed

\subsection{Simplifying the Fluid Approximation Under Independent Demands}

Problem (\ref{eqn:lp}) has a simplified form under independent demands. Using the decision variables \mbox{$\xvec = (x_{jt}^k : j \in \Jcal,~t \in \Tcal,~k \in \Kcal)$}, we consider the problem
\begin{align}
\max_{\xvec \in \mathbb R_+^{|\Jcal| TK}} ~~&  \sum_{k \in \Kcal} \sum_{t \in \Tcal} \sum_{j \in \Jcal} f_j \ts \mathbb P \{  D^k \geq t \} \ts x_{jt}^k 
\label{eqn:lp_indep}
\\
\mbox{st}~~&\sum_{\ell =1}^{k-1} \sum_{t \in \Tcal}  \sum_{j \in \Jcal} a_{ij} \ts \mathbb P \{ D^\ell \geq t \} \ts x_{jt}^\ell  + \sum_{t \in \Tcal} \sum_{j \in \Jcal} a_{ij} \ts x_{jt}^k \leq c_i \quad \forall \ts i \in \Lcal,~k \in \Kcal
\nonumber
\\
& x_{jt}^k \leq \lambda_{jt}^k \quad \forall \ts j \in \Jcal,~t \in \Tcal,~k \in \Kcal. \phantom{\Bigg\}}
\nonumber
\end{align}
We shortly show that if the demands are independent, then problem (\ref{eqn:lp}) is equivalent to problem (\ref{eqn:lp_indep}). Furthermore, problem (\ref{eqn:lp_indep}) is useful to see that naive fluid approximations are inadequate to deal with multiple stages even under independent demands. In particular, the first constraint above uses the distribution of the demand for stages $1,\ldots,k-1$, but it does not use the distribution of the demand for stage $k$. The specific form of this constraint is critical.~A naive approach may use a constraint of the form $\sum_{\ell=1}^k \sum_{t \in \Tcal} \sum_{j \in \Jcal} a_{ij} \ts \mathbb P \{ D^\ell \geq t\} \ts x_{jt}^\ell \leq c_i$ or $\sum_{\ell=1}^k \sum_{t \in \Tcal} \sum_{j \in \Jcal} a_{ij} \ts x_{jt}^\ell \leq c_i$. We also shortly demonstrate that the former option does not yield an asymptotically tight fluid approximation, whereas the latter option does not even yield an upper bound on the optimal total expected revenue. \cite{BaHo23} give fluid approximations under random number of customer arrivals that occur in one stage. Because they have one stage, their fluid approximations do not capture dependence between the demands in different stages. Even under independent demands, using one of the two options earlier in this paragraph to extend their work to multiple stages yields fluid approximations that are not asymptotically tight. 

{\bf \underline{Equivalence of Problems (\ref{eqn:lp}) and (\ref{eqn:lp_indep})}:}
\\
\indent We show that if the demands in different stages are independent, then the optimal objective values of (\ref{eqn:lp}) and (\ref{eqn:lp_indep}) are equal to each other and we can construct an optimal solution to one problem by using an optimal solution to the other. First, letting \mbox{$\xvecbar = (\xbar_{jt}^k(q) : j \in \Jcal,~t,q \in \Tcal,~k \in \Kcal)$} be an optimal solution to (\ref{eqn:lp}), we define the solution $\yvecbar = (\ybar_{jt}^k : j \in \Jcal,~t \in \Tcal,~k \in \Kcal)$ to (\ref{eqn:lp_indep}) as $\ybar_{jt}^k = \sum_{q \in \Tcal} \mathbb P\{ D^{k-1} = q \} \ts \xbar_{jt}^k(q)$.  We claim that the solution $\yvecbar$ is feasible to problem (\ref{eqn:lp_indep}) and provides an objective value that is equal to the optimal objective value of problem (\ref{eqn:lp}). For problem (\ref{eqn:lp_indep}), the solution $\yvecbar$ provides an objective value of  $\sum_{k \in \Kcal} \sum_{t \in \Tcal} \sum_{j \in \Jcal} f_j \ts \mathbb P \{ D^k \geq t\} \ts \ybar_{jt}^k = \sum_{k \in \Kcal} \sum_{t \in \Tcal} \sum_{j \in \Jcal} f_j \ts \mathbb P \{ D^k \geq t\} \sum_{q \in \Tcal} \mathbb P \{ D^{k-1} = q\} \ts \xbar_{jt}^k(q)$. Because the demands are independent, we have \mbox{$\mathbb P \{ D^k \geq t\} \ts \mathbb P \{ D^{k-1} = q\} = \mathbb P\{ D^k \geq t,~D^{k-1} = q\}$}, in which case, the last sum corresponds to the optimal objective value of problem (\ref{eqn:lp}). Thus, the solution $\yvecbar$ provides an objective value for (\ref{eqn:lp_indep}) that is equal to the optimal objective value of (\ref{eqn:lp}). 
Because $\xvecbar$ is an optimal solution to problem (\ref{eqn:lp}), it satisfies the second constraint in this problem. Multiplying this constraint with $\mathbb P \{ D^{k-1} = q\}$ and adding over all $q \in \Tcal$, we obtain $\ybar_{jt}^k = \sum_{q \in \Tcal} \mathbb P \{ D^{k-1} = q\} \ts \xbar_{jt}^k(q) \leq \lambda_{jt}^k$, in which case, it follows that the solution $\yvecbar$ satisfies the second constraint in  problem (\ref{eqn:lp_indep}).
Under independent demands, for $\ell = 1,\ldots,k-1$, we have the identity $\sum_{q \in \Tcal} \mathbb P \{ D^\ell \geq s \ts | \ts D^k \geq t,~ D^{k-1} = q \} \ts\ts \mathbb P \{ D^{k-1} = q \} = \mathbb P \{ D^\ell \geq s\}$.  In particular, if $\ell \leq k-2$, then this identity holds because $D^\ell$, $D^{k-1}$ and $D^k$ are independent. If, on the other hand, $\ell = k-1$, then we have $\sum_{q \in \Tcal} \mathbb P\{ D^{k-1} \geq s \ts | \ts D^k \geq t,~D^{k-1} = q\} \ts \mathbb P \{ D^{k-1} =q\} = \sum_{q \in \Tcal} \mathbb P \{ D^{k-1} \geq s,~D^{k-1} = q\}$, where we use the fact that $D^k$ and $D^{k-1}$ are independent. The identity follows by noting that the last sum is equal to $\mathbb P \{ D^{k-1} \geq s\}$. Because $\xvecbar$ is an optimal solution to problem (\ref{eqn:lp}), it satisfies the first constraint in this problem  with $t = T$. Multiplying this constraint with $\mathbb P \{ D^{k-1} = q\}$ and adding over all $q \in \Tcal$, as well as using the fact that the demands are independent, we obtain the inequality 
\begin{multline*}
\sum_{q \in \Tcal} \sum_{\ell =1}^{k-1} \sum_{s \in \Tcal} \sum_{p \in \Tcal} \sum_{j \in \Jcal} a_{ij} \ts  \mathbb P\{ D^{\ell-1} =p \} \ts \mathbb P \{ D^\ell \geq s \ts|\ts D^k \geq t, \ts D^{k-1} = q \} \ts \mathbb P \{ D^{k-1} = q\} \ts  \xbar_{js}^\ell(p) 
\\
  + \sum_{q \in \Tcal} \sum_{s \in \Tcal} \sum_{j \in \Jcal} a_{ij} \ts \mathbb P \{ D^{k-1} = q \} \ts \xbar_{js}^k(q) ~\leq~ c_i. 
\end{multline*}
Arranging the terms on the left side of the inequality above, as well as noting the fact that we have $\ybar_{jt}^k = \sum_{q \in \Tcal} \mathbb P\{ D^{k-1} = q \} \ts \xbar_{jt}^k(q)$ by definition of $\ybar_{jt}^k$ and the identity established just before the inequality above, the left side of the inequality above is $\sum_{\ell = 1}^{k-1} \sum_{s \in \Tcal} \sum_{j \in \Jcal} a_{ij} \ts \mathbb P \{ D^\ell \geq s\} \ts \ybar_{js}^\ell + \sum_{s \in \Tcal} \sum_{j \in \Jcal} a_{ij} \ts \ybar_{js}^k$. Thus, the solution $\yvecbar$ satisfies the first constraint in problem (\ref{eqn:lp_indep}) as well. In this case, the solution $\yvecbar$ is feasible to problem (\ref{eqn:lp_indep}) and the desired claim follows.  

Second, letting $\yvecbar = (\ybar_{jt}^k : j \in \Jcal,~t \in \Tcal,~k \in \Kcal)$ be an optimal solution to problem (\ref{eqn:lp_indep}), we define the solution \mbox{$\xvecbar = (\xbar_{jt}^k(q) : j \in \Jcal,~t,q \in \Tcal,~k \in \Kcal)$} to problem (\ref{eqn:lp}) as $\xbar_{jt}^k(q) = \ybar_{jt}^k$. In this case, we can follow an argument analogous to,  but slightly simpler than, the one in the previous paragraph to show that $\xvecbar$ is a feasible solution to problem (\ref{eqn:lp}) and provides an objective value that is equal to the optimal objective value of problem (\ref{eqn:lp_indep}).~Therefore, by the discussion so far in this section, if the demands in different stages are independent, then the optimal objective values of problems (\ref{eqn:lp}) and (\ref{eqn:lp_indep}) are equal to each other and we can construct an optimal solution to one problem by using an optimal solution to the other. 
In the rest of this section, we give two problem instances. We use the first problem instance to demonstrate that if we replace the first constraint in problem (\ref{eqn:lp_indep}) with  \mbox{$\sum_{\ell=1}^k \sum_{t \in \Tcal} \sum_{j \in \Jcal} a_{ij} \ts \mathbb P \{ D^\ell \geq t\} \ts x_{jt}^\ell \leq c_i$}, then the fluid approximation in (\ref{eqn:lp_indep}) is not asymptotically tight. We use the second problem instance to demonstrate that if we replace the first constraint in problem (\ref{eqn:lp_indep}) with $\sum_{\ell=1}^k \sum_{t \in \Tcal} \sum_{j \in \Jcal} a_{ij} \ts x_{jt}^\ell \leq c_i$, then the optimal objective value of problem (\ref{eqn:lp_indep}) does not provide an upper bound on the optimal total expected revenue.

{\bf \underline{Losing the Asymptotic Tightness of the Fluid Approximation}:}
\\
\indent
We give a problem instance with independent demands to demonstrate that if we replace the first constraint in (\ref{eqn:lp_indep}) with  $\sum_{\ell=1}^k \sum_{t \in \Tcal} \sum_{j \in \Jcal} a_{ij} \ts \mathbb P \{ D^\ell \geq t\} \ts x_{jt}^\ell \leq c_i$, then the fluid approximation in (\ref{eqn:lp_indep}) is not asymptotically tight. In our problem instance, there is one stage, one resource and one product. Thus, we have $K=1$, $\Lcal = \{1\}$ and $\Jcal = \{1\}$. The capacity of the resource is $c_1 = \alpha^3$ for some $\alpha \in \mathbb Z_+$. The revenue of the product is $f_1 = 1$. The demand in the single stage has two possible values with $\mathbb P \{ D^1 = \alpha \} = 1 - \frac{1}{\alpha^2}$ and $\mathbb P \{ D^1 = \alpha^4 \} = \frac{1}{\alpha^2}$. At each time period, there is a request for the product with probability one.
We have $K=1$ and $c_{\min} = \alpha^3$.~Because there is one resource, we have $L=1$. As a function of the parameter $\alpha$, using $\opt^\alpha$ to denote the optimal total expected revenue, $\apx^\alpha$ to denote the total expected revenue of the approximate policy and $\Zbar_\lp^\alpha$ to denote the optimal objective value of  (\ref{eqn:lp_indep}), by Theorem \ref{thm:perf_indep}, we have $\frac{\opt^\alpha}{\Zbar_\lp^\alpha} \geq \frac{\apx^\alpha}{\Zbar_\lp^\alpha} \geq  1 - 4 \ts \frac{\sqrt{\alpha^3 \log \alpha^3}}{\alpha^3} - \frac{1}{\alpha^3}$. By Theorem \ref{thm:bound}, we have $1 \geq \frac{\opt^\alpha}{\Zbar_\lp^\alpha}$. Thus, as $\alpha$ becomes arbitrarily large, the relative gap between the optimal total expected revenue and the optimal objective value of problem (\ref{eqn:lp_indep}) vanishes. 

We proceed to demonstrating that if we replace the first constraint in problem (\ref{eqn:lp_indep})~with  $\sum_{\ell=1}^k \sum_{t \in \Tcal} \sum_{j \in \Jcal} a_{ij} \ts \mathbb P \{ D^\ell \geq t\} \ts x_{jt}^\ell \leq c_i$, then the optimal total expected revenue and the optimal objective value of problem (\ref{eqn:lp_indep}) satisfy $\frac{\opt^\alpha}{\Zbar_\lp^\alpha} = \frac{2 \ts \alpha^2 -1}{\alpha^3 + \alpha^2-1}$. Therefore, as $\alpha$ becomes arbitrarily large, the relative gap between the optimal total expected revenue and the optimal objective value of problem (\ref{eqn:lp_indep}) becomes arbitrarily large, so we lose the asymptotic tightness of the fluid approximation. We compute the optimal total expected revenue. Because there is a single product, it is optimal to accept the requests for the product as much as the capacity allows.~Noting that there is a product request at each time period, if $D^1 = \alpha$, then the capacity of $\alpha^3$ allows us to accept all product requests at $\alpha$ time periods and we obtain a total expected revenue of $\alpha$. If $D^1 = \alpha^4$, then we can sell all of the capacity by using the product requests at $\alpha^4$ time periods, so noting that the capacity of the resource is $\alpha^3$, we obtain a total expected revenue of $\alpha^3$. Thus, the optimal total expected revenue is $\opt^\alpha = \mathbb P \{ D^1 = \alpha\} \ts \alpha  + \mathbb P \{ D^1 = \alpha^4\} \ts \alpha^3 = ( 1 - \frac{1}{\alpha^2} ) \ts \alpha + \frac{1}{\alpha^2} \ts \alpha^3 = \frac{2 \ts \alpha^2 -1}{\alpha}$. For our problem instance, if we replace the first constraint in problem (\ref{eqn:lp_indep}) with   \mbox{$\sum_{\ell=1}^k \sum_{t \in \Tcal} \sum_{j \in \Jcal} a_{ij} \ts \mathbb P \{ D^\ell \geq t\} \ts x_{jt}^\ell \leq c_i$}, then this problem is given by 
\begin{align*}
\Zbar_\lp^\alpha = \max ~~~&\sum_{t=1}^\alpha x_{1t}^1 + \sum_{t=\alpha+1}^{\alpha^4} \frac{1}{\alpha^2} \ts x_{1t}^1
\\
\mbox{st}~~~&\sum_{t=1}^\alpha x_{1t}^1 + \sum_{t=\alpha+1}^{\alpha^4} \frac{1}{\alpha^2} \ts x_{1t}^1 \leq \alpha^3
\\
& 
x_{1t}^1 \in [0,1] ~~\forall \ts t =1 ,\ldots,\alpha^4. \phantom{\sum_{\alpha}^{\alpha}}
\end{align*}
Because $\alpha + (\alpha^4 - \alpha) \ts \frac{1}{\alpha^2} \leq \alpha^3$ for all $\alpha \in \mathbb Z_+$, setting all of the decision variables to one yields a feasible, as well as an optimal, solution, so $\Zbar_\lp^\alpha = \alpha + (\alpha^4 - \alpha) \ts \frac{1}{\alpha^2}  = \frac{\alpha^3 + \alpha^2-1}{\alpha}$. 

\vspace{-0.5mm}

By the discussion in the previous paragraph, we have $\opt^\alpha = \frac{2 \ts \alpha^2 -1}{\alpha}$ and $\Zbar_\lp^\alpha = \frac{\alpha^3 + \alpha^2-1}{\alpha}$ for our problem instance, which implies that $\frac{\opt^\alpha}{\Zbar_\lp^\alpha} = \frac{2 \ts \alpha^2 -1}{\alpha^3 + \alpha^2-1}$. 

\vspace{-0.5mm}

{\bf \underline{Losing the Upper Bound from the Fluid Approximation}:}
\\
\indent
We give a problem instance with independent demands to demonstrate that if we replace the first constraint in (\ref{eqn:lp_indep}) with $\sum_{\ell=1}^k \sum_{t \in \Tcal} \sum_{j \in \Jcal} a_{ij} \ts x_{jt}^\ell \leq c_i$, then the optimal objective value of problem (\ref{eqn:lp_indep}) is not an upper bound on the optimal total expected revenue. In our problem instance, there are two stages, one product and one resource, so $K = 2$, $\Lcal = \{1\}$ and $\Jcal = \{1\}$.~The capacity of the resource is~$3$.~The revenue of the product is $f_1 = 1$.  The demand in each of the two stages has two possible values with $\mathbb P \{ D^1 = 1 \} = \mathbb P \{ D^1 = 2 \} = \mathbb P \{ D^2 = 1\} = \mathbb P \{ D^2 = 2 \} = \frac 12$. At each time period in each stage, we have a request for the product with probability one. We proceed to computing the optimal total expected revenue. Because there is a single product, it is optimal to accept the requests for the product as much as the capacity allows. With probability $\frac 12 \times \frac 12 = \frac 14$, the total number of requests for the product over the selling horizon is two, in which case, we make a revenue of two. With probability $1 - \frac 14 = \frac 34$, the total number of requests for the product over the selling horizon is three or more, in which case, noting that the capacity of the resource is three, we make a revenue of three.~Therefore, the optimal total expected revenue is given by $\opt = \frac 14 \times 2 + \frac 34 \times 3 = \frac{11}{4}$. For our problem instance, if we replace the first constraint in problem (\ref{eqn:lp_indep}) with $\sum_{\ell=1}^k \sum_{t \in \Tcal} \sum_{j \in \Jcal} a_{ij} \ts x_{jt}^\ell \leq c_i$, then this problem is given by
\begin{align*}
\Zbar_\lp ~=~\max~~~&x_{11}^1 + \frac 12 \ts x_{12}^1 + x_{11}^2  + \frac12 \ts x_{12}^2
\\
\mbox{st}~~~ & x_{11}^1 +  x_{12}^1 \leq 3
\\
& x_{11}^1 + x_{12}^1 + x_{11}^2  +  x_{12}^2 \leq 3
\\
& x_{11}^1, x_{12}^1,  x_{11}^2,  x_{12}^2 \in [0,1].
\end{align*}
Setting all of the decision variables other than $x_{12}^2$ to one yields an optimal solution to the problem above, so  $\Zbar_\lp = 1 + \frac 12 + 1 = \frac 52$. Therefore, for our problem instance, we have $\Zbar_\lp = \frac 52 < \frac{11}{4} =  \opt $.

\section{Extension to Assortment Offering and Pricing Decisions}
\label{sec:assort}

\vspace{-2mm}

We give the extension of our results to the case where we make assortment offering or pricing decisions. For economy of space, we only point out the adjustments needed in our development.

\vspace{-1mm}

{\bf \underline{Problem Formulation\phantom{p}\!\!\!}:}
\\
\indent We give the pieces of notation that we need in addition to those given in Section \ref{sec:form}. We pick an assortment of products to offer to the customer arriving at each time period. The customer either makes a choice within the offered assortment or leaves without making a purchase. Our goal is to find a policy to decide which assortment to offer to each customer so that we maximize the total expected revenue over the selling horizon. Given that we offer the assortment of products $S \subseteq \Jcal$, we use $\phi_{jt}^k(S)$ to denote the probability that the customer arriving at time period $t$ in stage~$k$ chooses product $j$. We have $\phi_{jt}^k(S) = 0$ for all $j \in \Jcal \setminus S$.~With probability $1 - \sum_{j \in \Jcal}\phi_{jt}^k(S)$, the customer leaves without making a purchase. As a function of the remaining capacities $\yvec = (y_i : i \in \Lcal) \in \mathbb Z_+^{|\Lcal|}$, the set of feasible decisions is $\Fcal(\yvec) = \{ S \subseteq \Jcal : a_{ij} \ts {\bf 1}(j \in S) \leq y_i ~\forall \ts i \in \Lcal,~j \in \Jcal\}$, ensuring that if we want to offer an assortment that includes product $j$ and the product uses resource $i$, then we need to have at least one unit of remaining capacity for resource $i$. We can find the optimal policy by computing the value functions $\{ J_t^k : t \in \Tcal,~k \in \Kcal\}$ through the dynamic program
\begin{align*}
J_t^k(\yvec,q)
& = \max_{S \in \Fcal(\yvec)} \Bigg\{ \sum_{j \in \Jcal} \phi_{jt}^k(S) \ts \bigg\{ f_j  + \theta_t^k(q) \ts J_{t+1}^k(\yvec -\avec_j ,q) + ( 1 - \theta_t^k(q)) \ts J_1^{k+1}(\yvec - \avec_j ,t) \bigg\} 
\\
& \quad \qquad \qquad \qquad \qquad \qquad + \bigg( 1 - \sum_{j \in \Jcal} \phi_{jt}^k(S) \ts \bigg) \bigg\{\theta_t^k(q) \ts J_{t+1}^k(\yvec ,q) + ( 1 - \theta_t^k(q)) \ts J_1^{k+1}(\yvec ,t ) \bigg\}
\Bigg\},
\end{align*}
with the boundary condition that $J_1^{K+1} = 0$. The demand in the stage just before the beginning of the selling horizon is fixed at $\Dhat^0$, so the optimal total expected revenue is $\opt = J_1^1(\cvec , \Dhat^0)$. 

Our formulation with the assortment offering decisions in the dynamic program above can capture pricing decisions as a special case. To use the formulation above to capture pricing decisions, we create multiple copies of each product, where a copy of a product corresponds to offering the product at a particular price. We refer to each copy of a product as a virtual product. If the set of possible price levels for a product is $\Pcal$, then we have $|\Jcal| \times |\Pcal|$ virtual products and we index the virtual products by $\{(j,p) \in \Jcal \times \Pcal\}$. In the pricing problem, we pick a subset of virtual products to offer to a customer with the constraint that if we offer the virtual product $(j,p)$, then we cannot offer another virtual product of the form $(j,q)$ with $q \in \Pcal \setminus \{p\}$. This constraint ensures that if we offer a product, then we offer it at only one price level. Given that we offer the assortment of virtual products $S \subseteq \Jcal \times \Pcal$, we use $\phi_{jpt}^k(S)$ to denote probability that the customer arriving at time period $t$ in stage $k$ chooses virtual product $(j,p)$. In this case, we can continue using the model in the previous paragraph with the understanding that virtual products play the role of products.~Thus, if there is a set of discrete possible price levels that we can charge for each product, then our discussion in this section applies to the pricing setting.

{\bf \underline{Fluid Approximation}:}
\\
\indent We use the decision variable $x_t^k(S,q)$ to capture the probability of offering assortment $S$ at time period $t$ in stage $k$ given that the demand in the previous stage was $q$. Consider the problem
\begin{align}
\Zbar_\lp =\!\!\! \max_{\xvec \in \mathbb R_+^{2^{|\Jcal|} T^2K}} ~~&  \sum_{k \in \Kcal} \sum_{t \in \Tcal} \sum_{S \subseteq \Jcal} \sum_{q \in \Tcal} \sum_{j \in \Jcal} f_j \ts \mathbb P \{  D^k \geq t, \ts D^{k-1} = q \} \ts \phi_{jt}^k(S) \ts x_t^k(S,q) 
\label{eqn:lp_choice}
\\
\mbox{st}~~&\sum_{\ell =1}^{k-1} \sum_{r \in \Tcal} \sum_{S \subseteq \Jcal} \sum_{p \in \Tcal} \sum_{j \in \Jcal} a_{ij} \ts \mathbb P \{ D^\ell \geq r, \ts D^{\ell-1} =p \ts|\ts D^k \geq t, \ts D^{k-1} = q \} \ts \phi_{jr}^\ell(S) \ts x_r^\ell(S,p) 
\nonumber
\\
& \qquad \qquad \qquad ~~~  + \sum_{r=1}^t \sum_{S \subseteq \Jcal} \sum_{j \in \Jcal} a_{ij} \ts \phi_{jr}^k(S) \ts x_r^k(S,q) \leq c_i \qquad \forall \ts i \in \Lcal,~t,q \in \Tcal,~k \in \Kcal
\nonumber
\\
& \sum_{S \subseteq \Jcal} x_t^k(S,q) = 1 \qquad \forall \ts t,q \in \Tcal,~k \in \Kcal, \phantom{\Bigg\}}
\nonumber
\end{align}
where we use the vector of decision variables $\xvec = (x_t^k(S,q) : S \subseteq \Jcal,~t,q \in \Tcal,~k \in \Kcal)$. Problem (\ref{eqn:lp_choice}) is the analogue of problem (\ref{eqn:lp})  under assortment offering decisions. The left side of the first constraint corresponds to the total expected capacity consumption of resource $i$ up to and including time period $t$ in stage $k$ conditional on the fact that the demand in stage $k$ is at least $t$ and the demand in stage $k-1$ is $q$. The second constraint ensures we offer some assortment with probability one at each time period in each stage, but we allow offering an empty assortment. We can show that the optimal objective value of the problem above is an upper bound on the optimal total expected revenue, so $\Zbar_\lp \geq \opt$. This result is the analogue of Theorem \ref{thm:bound}. We can show this result using any of the three approaches discussed at the end of Section \ref{sec:lp}, which include  Lagrangian relaxation, linear value function approximations and sample paths of the optimal policy.

{\bf \underline{Approximate Policy and Asymptotic Optimality}:}
\\
\indent We can use problem (\ref{eqn:lp_choice}) to construct an approximate policy. In particular, we solve problem (\ref{eqn:lp_choice}) once at the beginning of the selling horizon. Letting $\xvecbar = (\xbar_t^k(S,q) : S \subseteq \Jcal,~t,q \in \Tcal,~k \in \Kcal)$ be an optimal solution to this problem, the approximate policy makes it decisions as follows. Using \mbox{$\gamma \in [0,1]$} to denote a tuning parameter, if the demand in stage $k-1$ was $q$, then we offer assortment $S$ at time period $t$ in stage $k$ with probability $ \gamma \ts \xbar_t^k(S,q)$. With probability $1 - \gamma \sum_{S \subseteq \Jcal} x_t^k(S,q)$, we offer the empty assortment. If the customer chooses a product within the offered assortment for which we do not have remaining resource capacity to serve, then the customer leaves without a purchase.~Otherwise, the customer purchases the chosen product. We use $c_{\min} = \min_{i \in \Lcal} c_i$ to denote the smallest resource capacity and $L = \max_{j \in \Jcal} \sum_{i \in \Lcal} a_{ij}$ to denote the maximum number of resources used by a product. In this case, letting $\apx$ be the total expected revenue of the approximate policy, recalling that $\opt$ is the optimal total expected revenue and $\Zbar_\lp$ is the optimal objective value of the linear program in (\ref{eqn:lp_choice}), we can show that 
\begin{align*}
\frac{\apx}{\opt} 
~\geq~ 
\frac{\apx}{\Zbar_\lp} 
~\geq~ 
\max\Bigg\{ \frac{1}{4L} , \Bigg( 1 - 4 \ts \frac{\sqrt{(c_{\min} + \frac{1}{\epsilon^6} \ts (K-1) ) \log c_{\min}}}{c_{\min}} - \frac{L}{c_{\min}} \Bigg) \Bigg\}.
\end{align*}
This result is the analogue of Theorem \ref{thm:perf} under assortment offering decisions. The interpretation of the performance guarantee above is identical to that of the performance guarantee in Theorem \ref{thm:perf}. In many network revenue management settings, the numbers of resources and products can be large, but the number of resources used by a particular product is small. Therefore, the first part of the performance guarantee ensures that the approximate policy has a constant-factor performance guarantee when $L$ is uniformly bounded. On the other hand, considering a regime where we scale both the number of stages and resource capacities with the same rate $\theta$, the second part of the performance guarantee ensures that the relative gap between the total expected revenue of the approximate policy and the optimal total expected revenue converges to one with rate $1 - \frac{1}{\sqrt \theta}$. We proceed to giving an outline for the proof of the performance guarantee for our approximate policy under assortment offering decisions.

{\bf \underline{Performance Guarantee\phantom{p}\!\!\!}:}
\\
\indent The performance guarantee under assortment offering decisions follows from the same outline in Section \ref{sec:perf_asymp}. For each $k \in \Kcal$ and $t \in \Tcal$, we define four classes of Bernoulli random variables.

\indent {$\bullet$ \it {Demand in Each Stage}.} For each $q \in \Tcal$, the random variable $\Psi_t^k(q)$ takes value one if we reach time period $t$ in stage $k$ before this stage is over and the demand in stage $k-1$ is $q$. In other words, we have $\Psi_t^k(q) = {\bf 1}(D^k \geq t,~D^{k-1} = q)$. 

\indent {$\bullet$ \it {Policy Decision}.} For each $S \subseteq \Jcal$ and $q \in \Tcal$, the random variable $X_t^k(S,q)$ takes value one if the approximate policy  offers assortment $S$ at time period $t$ in stage $k$ when the demand in stage $k-1$ was $q$. By our approximate policy, $\mathbb P \{ X_t^k(S,q) = 1\} = \gamma \ts \xbar_t^k(S,q)$.

\indent {$\bullet$ \it {Customer Choice}.} For each $S \subseteq \Jcal$ and $j \in \Jcal$, given that the approximate policy offers the assortment $S$, the random variable $\Phi_{jt}^k(S)$ takes value one if the customer arriving at time period $t$ in stage $k$ purchases product $j$. We have $\mathbb P \{ \Phi_{jt}^k(S)  =1 \} = \phi_{jt}^k(S)$.

\indent {$\bullet$ \it {Availability}.} For each $j \in \Jcal$, the random variable $G_{jt}^k$ takes value one if we have enough capacity to accept a request for product $j$ at time period $t$ in stage $k$ under the approximate policy. Instead of calculating the probability $\mathbb P \{ G_{jt}^k =1\}$, we lower bound $\mathbb P \{ G_{jt}^k =1 \ts | \ts \Psi_t^k(q) = 1\}$.

The equality in (\ref{eqn:apx_rev}) continues to hold as long we replace $\xbar_{jt}^k(q)$ in this equality with $\sum_{S \subseteq \Jcal} \phi_{jt}^k(S) \ts \xbar_t^k(S,q)$. Thus, if we can show that $\mathbb P \{ G_{jt}^k = 1 \ts | \ts \Psi_t^k(q) =1\} \geq \alpha$, then $\apx \geq \gamma \ts \alpha \ts \Zbar_\lp$. Under the approximate policy, the sales for product $j$ at time period~$t$ in stage $k$ is given by $\sum_{q \in \Tcal} \sum_{S \subseteq \Jcal} \Psi_t^k(q) \ts G_{jt}^k \ts \Phi_{jt}^k(S) \ts X_t^k(S,q)$, where we use the fact that there is a sale for product $j$ at time period $t$ in stage $k$ if and only if we reach time period $t$ in stage $k$, we have enough capacity to accept a request for product $j$ and the customer chooses product~$j$ out of the offered assortment. Thus, we can use $\sum_{q \in \Tcal} \sum_{S \subseteq \Jcal} \Psi_t^k(q) \ts \Phi_{jt}^k(S) \ts X_t^k(S,q)$ as an upper bound on the sales for product $j$ at time period $t$ in stage $k$. In this case, $\sum_{j \in \Jcal} \sum_{q \in \Tcal} \sum_{S \subseteq \Jcal} a_{ij} \ts \Psi_t^k(q) \ts \Phi_{jt}^k(S) \ts X_t^k(S,q)$ is an upper bound on the capacity consumption of resource $i$ at time period $t$ in stage $k$. Letting \mbox{$N_{it}^k(q) = \sum_{S \subseteq \Jcal} \sum_{j \in \Jcal} a_{ij} \ts \Phi_{jt}^k(S) \ts X_t^k(S,q)$}, we express the last  upper bound as $\sum_{q \in \Tcal} \Psi_t^k(q) \ts N_{it}^k(q)$, in which case, the chain of inequalities in (\ref{eqn:union}) continues to hold. Lemmas \ref{lem:mgf_gap}, \ref{lem:mgf_bound} and \ref{lem:avail} go through with no modifications when we use the definition of $N_{it}^k(q)$ given in this paragraph. We obtain the performance guarantee for our approximate policy under assortment offering decisions by following precisely the same outline in the proof of Theorem \ref{thm:perf}. To sum up, if we decide whether to accept or reject each product request, then we use $\sum_{q \in \Tcal} \Psi_t^k(q) \ts A_{jt}^k \ts X_{jt}^k$ as an upper bound on the sales for product $j$ at time period $t$ in stage $k$. If we make assortment offering decisions, then we use $\sum_{q \in \Tcal} \sum_{S \subseteq \Jcal} \Psi_t^k(q) \ts \Phi_{jt}^k(S) \ts X_t^k(S,q)$ as an upper bound on the sales for product $j$ at time period $t$ in stage $k$. Once we draw this parallel, the performance guarantee under assortment offering decisions follows precisely the same outline in the proof of Theorem \ref{thm:perf}.

At the beginning of this section, we explain that our formulation with assortment offering decision captures pricing decisions as a special case. To elaborate further, if the set of possible price levels for a product is $\Pcal$, then the set of virtual products is $\Jcal \times \Pcal$. Offering virtual product $(j,p)$ implies that we  offer product $j$ at price level $p$. In the pricing problem, we choose an assortment of virtual product to offer at each time period in each stage and the arriving customers make a choice among the offered virtual products. A product, if offered, can be offered at a single price. Therefore, if the remaining capacities of the resources at any time period in any stage are given by the vector $\yvec = (y_i : i \in \Lcal)$, then the set of feasible assortments of virtual products that we can offer corresponds to $\Fcal(\yvec) = \{ S \subseteq \Jcal \times \Pcal : a_{ij} \ts {\bf 1}(S \cap (\{j\} \times \Pcal) ) \leq y_i ~\forall \ts i \in \Lcal,~j \in \Jcal,~~|S \cap (\{j\} \times \Pcal)| \leq 1 ~\forall \ts j \in \Jcal\}$, where the first constraint ensures that if we want to offer an assortment that includes a virtual product corresponding to product $j$ and product $j$ uses resource $i$, then we need to have at least one unit of remaining capacity for resource $i$, whereas the second constraint ensures that we can offer at most one virtual product corresponding to product $j$. To capture the reactions of the customers to prices, given that we offer the assortment of virtual products $S \subseteq \Jcal \times \Pcal$, we use $\phi_{jpt}^k(S)$ to denote the probability that the customer arriving at time period $t$ in stage $k$ choses virtual product $(j,p)$.~In this case, our model under assortment offering decisions  applies to pricing decisions, as long as we view each virtual product as a product under assortment offering decisions.

\section{Computational Experiments}
\label{sec:exp}

We conduct a numerical study to investigate the benefits from our fluid approximation when the demands in different stages have arbitrary distributions and are dependent on each other. 

{\bf \underline{Experimental Setup}:} In our computational experiments, we consider an airline network with one hub and three spokes. We have a flight leg that connects the hub to each spoke and each spoke to the hub. Thus, there are six flight legs. We have a high-fare and a low-fare itinerary that connects each origin-destination pair.~Thus, there are $2 \times 4 \times 3 = 24$ itineraries. In this setting, the resources correspond to the flight legs, whereas the products correspond to the itineraries.~The itineraries between the hub and a spoke are direct, each using a single resource, whereas the itineraries between two spokes connect at the hub, each using two resources.  To generate the revenues associated with the itineraries, we place the hub at the center of a $100\times 100$ square.~We place the spokes over the square uniformly at random. The revenue of a low-fare itinerary is the Euclidean distance between the origin and destination locations of the itinerary. The revenue of a high-fare itinerary is $\kappa$ times the revenue of the corresponding low-fare itinerary. The preceding discussion explains our approach for setting up the resources and products, along with the set of resources used by each product. We proceed to explaining our approach for generating the distribution of the demand random variables in each stage and the correlation structure between the demand random variables. %The experimental setup so far closely follows~\cite{To06}, where the author does not consider calendar-aware or dependent demands. We proceed to discussing how we come up with the demand random variables.

We have $K$ stages in the selling horizon. We will vary the parameter $K$. The demand in each stage has a truncated and discretized log-normal distribution. In particular, given that the demand in stage $k-1$ takes the value $\Dbar^{k-1}$, the mean and standard deviation of the demand in stage $k$ are, respectively, $\mu^k(\Dbar^{k-1}) = \rho \ts \Dbar^{k-1} + (1-\rho) \ts 100$ and $\sigma^k(\Dbar^{k-1}) = 0.3 \ts (\rho \ts \Dbar^{k-1} + (1-\rho) \ts 100)$. In this case, given that the demand in stage $k-1$ takes the value $\Dbar^{k-1}$, the demand in stage $k$ is obtained by rounding a log-normal random variable with mean $\mu^k(\Dbar^{k-1})$ and standard deviation $\sigma^k(\Dbar^{k-1})$ up to the closest integer and truncating at $\lceil \mu^k + 3 \ts \sigma^k \rceil$. We fix $D^0$ at $\Dhat^0 = 100$.~Therefore, the demand in each stage, before rounding and truncating, always has the mean of 100 and coefficient of variation of 0.3. The parameter $\rho$ controls the correlation between the demand random variables in successive stages. Because a linear transformation of a log-normal random variable is not a \mbox{log-normal} random variable, the correlation coefficient between the demands in successive stages may not be necessarily $\rho$. After generating each test problem, however, we computed the correlation coefficient between the demands in successive stages and the correlation coefficient came out to be quite close to $\rho$. We will vary the parameter $\rho$.

To come up with the request arrival probabilities $\{ \lambda_{jt}^k : j \in \Jcal,~t \in \Tcal,~k \in \Kcal\}$, we assume that the~requests for the low-fare itineraries tend to arrive towards the beginning of the selling horizon, whereas the requests for the high-fare itineraries tend to arrive towards the end of the selling horizon. In this way, we generate test problems in which it is important to protect capacities for  the high-fare itinerary requests that tend to arrive later. In particular, the probability of getting a request for a low-fare itinerary linearly decreases over time. There are $K$ stages, each with at most~$T$ time periods. We sample a threshold from the uniform distribution over $\{\lceil \frac{1}{2} KT \rceil, \ldots, \lceil \frac{2}{3} KT \rceil\}$ so that the requests for the high-fare itineraries arrive only after the threshold and the probability of getting a request for a high-fare itinerary linearly increases over time. We use the following approach to come up with such request arrival probabilities.

For each origin-destination pair $(f,g)$, we sample $\zeta_{fg}$ from the uniform distribution over $[0,1]$. Using~$N$ to denote the set of all locations, we  set $\gamma_{fg} = \frac{\zeta_{fg}}{\sum_{k \in N} \sum_{\ell \in N \setminus \{\ell\}} \zeta_{k\ell}}$ so that the parameters \mbox{$\{ \gamma_{fg} : f \in N,~g \in N \setminus \{f\}\}$} are normalized to add up to one. At any time period in any stage, we have a request for an itinerary that connects origin-destination pair $(f,g)$ with probability~$\gamma_{fg}$. There is a \mbox{high-fare} and a low-fare itinerary that connects each origin-destination pair. The probability that we have a request for a low-fare itinerary linearly increases over time.~The probability that we have a request for a high-fare itinerary is zero until the threshold and increases linearly over time after the threshold. Therefore, for each origin-destination pair $(f,g)$, we sample the threshold $\tau_{fg}$ from the uniform distribution over $\{\lceil \frac{1}{2} KT \rceil, \ldots, \lceil \frac{2}{3} KT \rceil\}$. The function $F(t) = \frac{KT+1-t}{KT}$ is decreasing for \mbox{$t \in [1,KT]$}, whereas the function $H_{fg}(t) = \frac{[t-\tau_{fg}]^+}{KT - \tau_{fg}}$ takes the value zero for $t \in[ 1,\tau_{fg}]$ and is increasing for $t \in [\tau_{fg},KT]$. In this case, noting that there are $(k-1)\ts T$ time periods before stage~$k$, if product $j$ corresponds to a low-fare itinerary that connects origin-destination pair $(f,g)$, then we set $\lambda_{jt}^k = \gamma_{fg} \ts \frac{F((k-1)\ts T + t)}{F((k-1)\ts T + t) + H_{fg}((k-1)\ts T + t)}$, whereas if product $j$ corresponds to a high-fare itinerary that connects origin-destination pair $(f,g)$, then we set $\lambda_{jt}^k = \gamma_{fg} \ts \frac{H_{fg}((k-1)\ts T + t)}{F((k-1)\ts T + t) + H_{fg}((k-1)\ts T + t)}$.

Once we generate the request arrival probabilities, the total expected demand for the capacity on fight leg $i$ is $\xi_i = \sum_{k \in \Kcal} \sum_{t \in \Tcal} \sum_{j \in \Jcal} a_{ij} \ts \mathbb P \{ D^k \geq t \ts \} \ts \lambda_{jt}^k$. We set the capacity of flight leg~$i$ as $c_i = \lceil \xi_i / \beta \rceil$, where the parameter $\beta$ controls the tightness of~the capacities.~Thus, the total expected demand for the capacity on a flight leg exceeds the capacity of the flight leg by a factor of $\beta$. 
We fix $\kappa = 8$ and $\beta = 1.6$. We experimented with different values for these two parameters and our results remained qualitatively the same. We build on the log-normal distribution to capture the demand random variables. Using this distribution, we can increase the coefficient of variation of the demand as much as we like and not be concerned about negative values. %Poisson or normal distributions, for example, do not provide such flexibility.

We vary $K \in \{5,10,15,20,25,30\}$ and $\rho \in \{0.2,0.4,0.6,0.8\}$.  For each parameter combination $(K,\rho)$, we generate a test problem using the approach described so far in this section.

{\bf \underline{Benchmark\phantom{p}\!\!\!}:} We use a benchmark fluid approximation that uses only the expected values of the demands. The total expected demand for product $j$ is $\sum_{k \in \Kcal} \sum_{t \in \Tcal} \mathbb P\{ D^k \geq t \} \ts \lambda_{jt}^k$. Using the decision variable $w_j$ to capture the total number of requests for product $j$ that we accept, using the vector $\wvec = (w_j : j \in \Jcal)$, we consider the linear program 
\begin{align}
\max_{\wvec \in \mathbb R_+^{|\Jcal|}} \Bigg\{ 
\sum_{j \in \Jcal} f_j \ts w_j~:~
\sum_{j \in \Jcal} a_{ij} \ts w_j \leq c_i ~~ \forall \ts i \in \Lcal,~~
 w_j \leq \sum_{k \in \Kcal} \sum_{t \in \Tcal} \mathbb P \{ D^k \geq t \} \ts \lambda_{jt}^k ~~ \forall \ts j \in \Jcal \Bigg\}.
 \label{eqn:simple_lp}
\end{align}
By the first constraint, the expected capacity consumption of each resource does not exceed its capacity. By the second constraint, the expected number of accepted requests for each product does not exceed its expected demand. The optimal objective value of (\ref{eqn:simple_lp}) is an upper bound on the optimal total expected revenue. In particular, letting $\xvecbar$ be an optimal solution to (\ref{eqn:lp}), we set $\wbar_j = \sum_{k \in \Kcal} \sum_{t \in \Tcal} \sum_{q \in \Tcal} \mathbb P\{ D^k \geq t,~D^{k-1} = q\} \ts \xbar_{jt}^k(q)$. Considering the first constraint in (\ref{eqn:lp}) for $k=K$, if we multiply this constraint with $\mathbb P \{ D^k = t,~D^{k-1} = q\}$ and add over all $t,q \in \Tcal$, then we can verify that $\wvecbar$ satisfies the first constraint in (\ref{eqn:simple_lp}). Considering the second constraint in (\ref{eqn:lp}), if we multiply this constraint with $\mathbb P \{ D^k \geq t,~D^{k-1} = q\}$ and add over all $t,q \in \Tcal$, then we can verify that $\wvecbar$ satisfies the second constraint in (\ref{eqn:simple_lp}). In this way, the solution $\wvecbar$ is feasible to problem (\ref{eqn:simple_lp}) and provides an objective value that is equal to the optimal objective value of problem (\ref{eqn:lp}). Thus, the optimal objective value of (\ref{eqn:simple_lp}) is an upper bound on that of (\ref{eqn:lp}), which is, in turn, an upper bound on the optimal total expected revenue. However, at the end of this section, we demonstrate that the upper bound from (\ref{eqn:simple_lp}) is not asymptotically tight. Using (\ref{eqn:simple_lp}), we can also obtain an approximate policy. Letting $\wvecbar$ be an optimal solution to (\ref{eqn:simple_lp}), the approximate policy is willing to accept a request for product $j$ at any time period with probability $\wbar_j / \sum_{k \in \Kcal} \sum_{t \in \Tcal} \mathbb P \{ D^k \geq t  \} \ts \lambda_{jt}^k$.

{\bf \underline{Computational Results\phantom{p}\!\!\!}:} We refer to the fluid approximation in (\ref{eqn:lp}) as $\prf$ to indicate that this fluid approximation uses information on the full probability distribution of the demands, whereas we refer to the fluid approximation in (\ref{eqn:simple_lp}) as $\exf$ to indicate that this fluid approximation uses information on only the expected values of the demands. In our computational experiments, we generate 24 test problems using the approach discussed earlier in this section. For each test problem, we solve $\prf$ and $\exf$ to compute the upper bounds on the optimal total  expected revenue provided by the two fluid approximations. Furthermore, we simulate the decisions of the approximate policies provided by $\prf$ and $\exf$ for 1000 sample paths to estimate the total expected revenues of these two approximate policies. We use common random numbers when simulating the decisions of the~two~approximate policies. We give our computational results in Table \ref{tab:exp}. In this table, the first column gives the parameter combinations for our test problems using the pair $(K,\rho)$. The second column gives~the upper bound on the optimal total expected revenue provided by $\prf$, whereas the third column gives the total expected revenue of the approximate policy from $\prf$. The fourth column gives the ratio between the total expected revenue of the approximate policy and the upper bound provided by $\prf$ multiplied by 100. The fifth, sixth and seventh columns give the same statistics for $\exf$. The eighth column gives the percent gap between the upper bounds provided by $\prf$ and $\exf$, whereas the ninth column gives the percent gap between the total expected revenues of the approximate policies provided by the two fluid approximations. 

\begin{table}
\begin{center}
\footnotesize
\begin{tabular}{|c|rrc|rrc|cr|}
\hline
Params. & \multicolumn{3}{c|}{$\prf$} & \multicolumn{3}{c|}{$\exf$} & Bnd. & Plcy.~ \\
\cline{2-7}
$(K,\rho)$ & Bound & Policy & ~Ratio~ & Bound & Policy & ~Ratio~ & ~~Gap~~ & Gap~ \\
\hline
\hline
$(5,0.2)$	&	32,752	&	30,185	&	92.16\%	&	34,533	&	28,206	&	81.68\%	&	5.16\%	&	7.01\%	\\
$(5,0.4)$	&	31,755	&	28,700	&	90.38\%	&	33,529	&	26,603	&	79.34\%	&	5.29\%	&	7.88\%	\\
$(5,0.6)$	&	30,244	&	26,866	&	88.83\%	&	31,983	&	24,404	&	76.30\%	&	5.44\%	&	10.09\%	\\
$(5,0.8)$	&	27,690	&	24,285	&	87.70\%	&	29,449	&	20,492	&	69.58\%	&	5.97\%	&	18.51\% \\
\hline
$(10,0.2)$	&	71,923	&	68,361	&	95.05\%	&	74,333	&	65,324	&	87.88\%	&	3.24\%	&	4.65\%	\\
$(10,0.4)$	&	69,342	&	64,991	&	93.73\%	&	72,099	&	61,155	&	84.82\%	&	3.82\%	&	6.27\%	\\
$(10,0.6)$	&	64,802	&	59,938	&	92.49\%	&	67,861	&	55,599	&	81.93\%	&	4.51\%	&	7.80\%	\\
$(10,0.8)$	&	55,250	&	49,893	&	90.30\%	&	58,482	&	42,911	&	73.37\%	&	5.53\%	&	16.27\%	\\
\hline
$(15,0.2)$	&	111,579	&	106,140	&	95.13\%	&	114,223	&	102,259	&	89.53\%	&	2.31\%	&	3.80\%	\\
$(15,0.4)$	&	107,571	&	102,204	&	95.01\%	&	110,808	&	97,262	&	87.78\%	&	2.92\%	&	5.08\%	\\
$(15,0.6)$	&	100,060	&	93,949	&	93.89\%	&	104,080	&	87,636	&	84.20\%	&	3.86\%	&	7.20\%	\\
$(15,0.8)$	&	82,458	&	76,242	&	92.46\%	&	87,337	&	66,527	&	76.17\%	&	5.59\%	&	14.60\%	\\
\hline
$(20,0.2)$	&	151,270	&	144,609	&	95.60\%	&	154,021	&	140,783	&	91.40\%	&	1.79\%	&	2.72\%	\\
$(20,0.4)$	&	145,995	&	139,229	&	95.37\%	&	149,455	&	133,671	&	89.44\%	&	2.31\%	&	4.16\%	\\
$(20,0.6)$	&	135,759	&	129,531	&	95.41\%	&	140,308	&	122,021	&	86.97\%	&	3.24\%	&	6.15\%	\\
$(20,0.8)$	&	110,234	&	105,274	&	95.50\%	&	116,391	&	92,082	&	79.11\%	&	5.29\%	&	14.33\%	\\
\hline
$(25,0.2)$	&	191,025	&	181,896	&	95.22\%	&	193,835	&	177,921	&	91.79\%	&	1.45\%	&	2.23\%	\\
$(25,0.4)$	&	184,421	&	176,305	&	95.60\%	&	188,005	&	170,534	&	90.71\%	&	1.91\%	&	3.38\%	\\
$(25,0.6)$	&	171,628	&	164,427	&	95.80\%	&	176,483	&	156,205	&	88.51\%	&	2.75\%	&	5.26\%	\\
$(25,0.8)$	&	138,672	&	131,287	&	94.67\%	&	145,754	&	117,589	&	80.68\%	&	4.86\%	&	11.65\%	\\
\hline
$(30,0.2)$	&	230,756	&	220,695	&	95.64\%	&	233,608	&	216,844	&	92.82\%	&	1.22\%	&	1.78\%	\\
$(30,0.4)$	&	222,995	&	214,150	&	96.03\%	&	226,656	&	207,869	&	91.71\%	&	1.62\%	&	3.02\%	\\
$(30,0.6)$	&	207,655	&	196,766	&	94.76\%	&	212,704	&	188,936	&	88.83\%	&	2.37\%	&	4.14\%	\\
$(30,0.8)$	&	167,445	&	158,680	&	94.77\%	&	175,196	&	143,432	&	81.87\%	&	4.42\%	&	10.63\%	\\															
\hline
\hline
Avg.			& & &			93.81\% & & &						84.43\%	&	3.62\%	&	7.44\% \\
\hline
\end{tabular}
\caption{Performance of the fluid approximations.}
\label{tab:exp}
\end{center}
\end{table}

Our results indicate that the approximate policy from $\prf$ performs quite well. Over all of our test problems, the average gap between the total expected revenues of the approximate policy and the upper bounds provided by $\prf$ is 6.19\%. In other words, noting that we compare the performance of the approximate policy with an upper bound on the optimal total expected revenue, rather than the optimal total expected revenue itself, a conservative estimate of the average optimality gap of the approximate policy provided by $\prf$ is 6.19\%. In alignment with Theorem~\ref{thm:perf}, the relative gaps between the upper bounds and the total expected revenues of the approximate policy from $\prf$ tend to vanish as the number of stages increases. Considering groups of test problems with 5, 10, 15, 20, 25 and 30 stages, on average, the ratios between the total expected revenues of the approximate policy and the upper bounds provided by $\prf$ are, respectively, $89.77$, $92.89$, $94.12$, $95.47$, $95.32$ and $95.30$.~The performance of $\exf$, both in terms of the tightness of the upper bounds and the total expected revenues of the approximate policy that it provides, is noticeably inferior to $\prf$. Over all of our test problems, the average gap between the total expected revenues of the approximate policy and the upper bounds provided by $\exf$ is 15.57\%. On average, the upper bounds on the optimal total expected revenue provided by $\prf$ improve those provided by $\exf$ by 3.62\%.~Similarly, the total expected revenues of the approximate policy provided by $\prf$ improve those provided by $\exf$ by 7.44\%. There are test problems where the gap between the two upper bounds can reach 5.97\% and the gap between the total expected revenues of the two approximate policies can reach 18.51\%. Overall, $\exf$ builds on the standard fluid approximation framework by using only the expected values of the demands, but there are significant benefits from using a fluid approximation that handles dependent demands with arbitrary distributions. In our computational experiments, we set the tuning parameter $\gamma$ to one in our implementation of the approximate policy for $\prf$. While using a value of $\gamma$ away from one is necessary to obtain our performance guarantee, our numerical experience indicates that using the value of one provides the strongest performance for our test problems. In Figure \ref{fig:gamma}, we show the total expected revenue obtained by $\prf$  as a function of $\gamma$ for our problem instances with $K = 5$. The performance of our approximate policy is strongest when we set $\gamma =1$. We observed a similar pattern in our other test problems.

\begin{figure}
\begin{center}
%\vspace{2in}
\includegraphics[scale = 0.4 ]{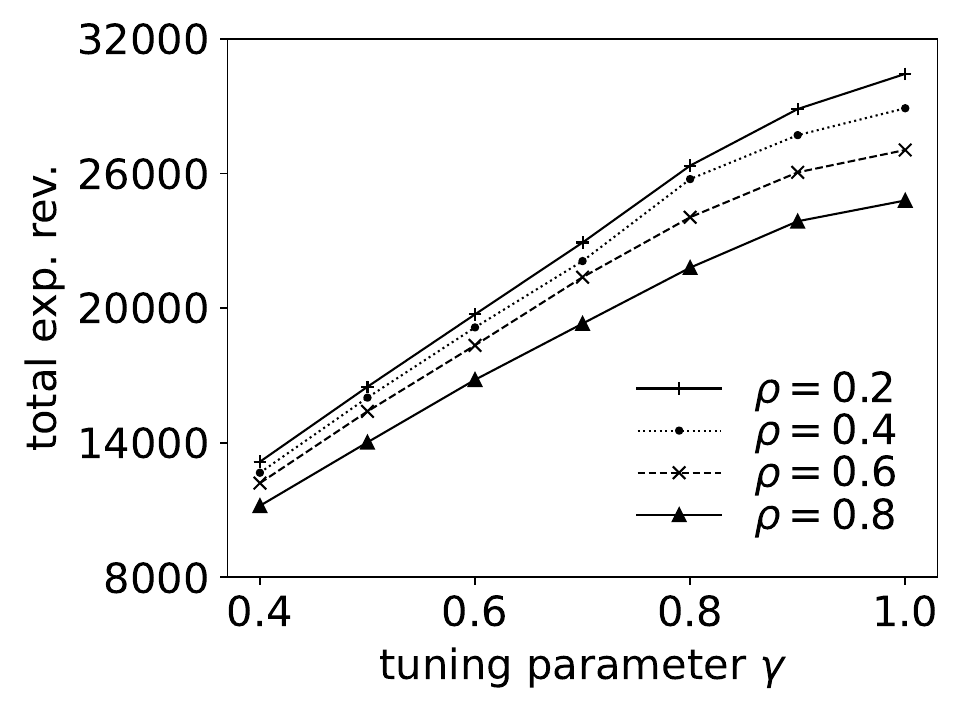}
\caption{Performance of $\prf$ on test problem with $K = 5$ for different values of the tuning parameter.}
\label{fig:gamma}
\vspace{-5mm}
\end{center}
\end{figure}

{\bf \underline{Upper Bound from Problem (\ref{eqn:simple_lp})}:} 
We give a problem instance to demonstrate that the upper bound provided by the optimal objective value of problem (\ref{eqn:simple_lp}) is not asymptotically tight. In our problem instance, there is one stage, one resource and one product, so $K = 1$, $\Lcal = \{1\}$ and \mbox{$\Jcal = \{1\}$}. The capacity of the resource is $\alpha^2$ for some $\alpha \in \mathbb Z_+$. The revenue of the product is \mbox{$f_1 =1$}.~The demand in the single stage has two possible values with $\mbox{$\mathbb P \{ D^1 = \alpha \}$} = 1 - \frac{1}{2 (\alpha-1)}$ and \mbox{$\mathbb P \{ D^1 = \alpha^3 \ts \} = \frac{1}{2(\alpha-1)}$}. At each time period, there is a request for the product with probability one. For this problem instance, we have $K=1$ and \mbox{$c_{\min} = \alpha^2$}.~Furthermore, noting that there is one resource, we have $L=1$. As a function of the parameter $\alpha$ in our problem instance, using $\opt^\alpha$ to denote the optimal total expected revenue, $\apx^\alpha$ to denote the total expected revenue of the approximate policy and $\Zbar_\lp^\alpha$ to denote the optimal objective value of problem (\ref{eqn:lp}), by Theorem~\ref{thm:perf}, we have  \mbox{$\frac{\opt^\alpha}{\Zbar_\lp^\alpha} \geq \frac{\apx^\alpha}{\Zbar_\lp^\alpha} \geq  1 - 4 \ts \frac{\sqrt{\alpha^2 \log \alpha^2}}{\alpha^2} - \frac{1}{\alpha^2}$}. By Theorem \ref{thm:bound}. we also have \mbox{$1 \geq \frac{\opt^\alpha}{\Zbar_\lp^\alpha}$}.~Thus, as $\alpha$ becomes arbitrarily large, the relative gap between the optimal total expected revenue and the optimal objective value of problem (\ref{eqn:lp}) vanishes for our problem instance.

Using $\Ztilde_\lp^\alpha$ to denote the optimal objective value of problem (\ref{eqn:simple_lp}) as a function of the parameter $\alpha$ in our problem instance, we argue that $\frac{\opt^\alpha}{\Ztilde_\lp^\alpha} = \frac{3}{\alpha + 3}$. Therefore, as $\alpha$ becomes arbitrarily large, the relative gap between the optimal total expected revenue and the optimal objective value of problem (\ref{eqn:simple_lp}) becomes arbitrarily large, so the optimal objective value of problem (\ref{eqn:simple_lp}) is not asymptotically tight. Because there is a single product, it is optimal to accept the requests for the product as much as the capacity allows. If $D^1 = \alpha$, then the capacity of $\alpha^2$ allows us to accept all product requests at $\alpha$ time periods and we obtain a total expected revenue of $\alpha$. If $D^1 = \alpha^3$, then we can sell all of the capacity by using the product requests at $\alpha^3$ time periods, so noting that the capacity of the resource is $\alpha^2$, we obtain a total expected revenue of $\alpha^2$. Thus, the optimal total expected revenue is $\opt^\alpha = \mbox{$\mathbb P \{ D^1 = \alpha \} \ts \alpha + \mathbb P \{ D^1 = \alpha^3 \} \ts \alpha^2$} = (1 - \frac{1}{2 \ts (\alpha-1)}) \ts \alpha + \frac{1}{2\ts (\alpha-1)} \ts \alpha^2 = \frac{3\alpha}{2}$. On the other hand, for our problem instance, the linear program in (\ref{eqn:simple_lp}) is given by
\begin{align*}
\Ztilde_\lp^\alpha ~=~~\max_{w_1 \geq 0} ~~ \Bigg\{ w_1  ~:~ w_1 \leq \alpha^2,~~w_1 \leq \alpha + (\alpha^3 - \alpha) \ts \frac{1}{2  (\alpha-1)} \Bigg\}.
\end{align*}
For $\alpha \geq 3$, we have $\alpha + (\alpha^3 - \alpha) \ts \frac{1}{2 (\alpha-1)} = \frac{\alpha^2 + 3 \alpha}{2} \leq \alpha^2$, so $\Ztilde_\lp^\alpha = \frac{\alpha^2 + 3 \alpha}{2}$. Noting that $\opt^\alpha = \frac{3\alpha}{2}$, we obtain $\frac{\opt^\alpha}{\Ztilde_\lp^\alpha} = \frac{3}{\alpha + 3}$ for our problem instance, as long as $\alpha \geq 3$.

\clearpage

\end{APPENDICES}

\end{document}